\documentclass[12pt,leqno]{article}

\def\diam{\mathop{\rm diam}}
\def\dist{\mathop{\rm dist}}

\newtheorem{theorem}{Theorem}
\newtheorem{lemma}[theorem]{Lemma}
\newtheorem{proposition}[theorem]{Proposition}
\newtheorem{sublemma}[theorem]{Sublemma}
\newtheorem{definition}[theorem]{Definition}
\newtheorem{corollary}[theorem]{Corollary}
\newtheorem{problem}[theorem]{Problem}
\newtheorem{remark}[theorem]{Remark}
\newtheorem{claim}[theorem]{Claim}
\newtheorem{assumptions}[theorem]{Assumptions}
\newtheorem{examples}[theorem]{Examples}
\newtheorem{basicfact}[theorem]{Basic Fact}

\newcommand{\begintheorem}{\addtocounter{equation}{1}\begin{theorem}}
\newcommand{\beginlemma}{\addtocounter{equation}{1}\begin{lemma}}
\newcommand{\beginproposition}{\addtocounter{equation}{1}\begin{proposition}}
\newcommand{\beginsublemma}{\addtocounter{equation}{1}\begin{sublemma}}
\newcommand{\begindefinition}{\addtocounter{equation}{1}\begin{definition}}
\newcommand{\begincorollary}{\addtocounter{equation}{1}\begin{corollary}}
\newcommand{\beginproblem}{\addtocounter{equation}{1}\begin{problem}}
\newcommand{\beginremark}{\addtocounter{equation}{1}\begin{remark}}
\newcommand{\beginclaim}{\addtocounter{equation}{1}\begin{claim}}
\newcommand{\beginassumptions}{\addtocounter{equation}{1}\begin{assumptions}}
\newcommand{\beginexamples}{\addtocounter{equation}{1}\begin{examples}}
\newcommand{\beginbasicfact}{\addtocounter{equation}{1}\begin{basicfact}}

\renewcommand{\thetheorem}{\arabic{section}.\arabic{equation}}
\renewcommand{\theequation}{\arabic{section}.\arabic{equation}}

\begin{document}

\title{Real Analysis, Quantitative Topology, and Geometric Complexity}


\author{Stephen Semmes}

\date{}

\maketitle


\renewcommand{\thefootnote}{}   

\footnotetext{The author was partially supported by the National
Science Foundation.}

	This survey originated with the John J. Gergen Memorial
Lectures at Duke University in January, 1998.  The author would like
to thank the Mathematics Department at Duke University for the
opportunity to give these lectures.  See \cite{Gromov-green,
Gromov-translation, Gromov-conf, Seapp} for related topics, in
somewhat different directions.

\tableofcontents

\newpage

\section{Mappings and distortion}
\label{Mappings and distortion}
\setcounter{equation}{0}

	A very basic mechanism for controlling geometric complexity is
to limit the way that distances can be distorted by a mapping.

	If distances are distorted by only a small amount, then one
might think of the mapping as being approximately ``flat''.  Let us
look more closely at this, and see what actually happens.

	Let $\delta$ be a small positive number, and let $f$ be a
mapping from the Euclidean plane ${\bf R}^2$ to itself.  Given two
points $x, y \in {\bf R}^2$, let $|x-y|$ denote the usual Euclidean
distance between them.  We shall assume that
\begin{equation}
\label{1+delta bilipschitz}
	(1 + \delta)^{-1} \, |x-y| \le |f(x) - f(y)| 
					\le (1 + \delta) \, |x-y|
\end{equation}
for all $x, y \in {\bf R}^2$.  This says exactly that $f$ does not
ever shrink or expand distances by more than a factor of $1+\delta$.

	What does this really mean about the behavior of $f$?  A first
point is that if $\delta$ were equal to $0$, so that $f$ does not
distort distances at all, then $f$ would have to be a ``rigid''
mapping.  This means that $f$ could be expressed as
\begin{equation}
\label{special form for isometries}
	f(x) = A(x) + b,
\end{equation}
where $b$ is an element of ${\bf R}^2$ and $A$ is a linear mapping on 
${\bf R}^2$ which is either a rotation or a combination of a rotation
and a reflection.  This is well known, and it is not hard to prove.
For instance, it is not hard to show that the assumption that $f$
preserve distances implies that $f$ takes lines to lines, and that
it preserve angles, and from there it is not hard to see that $f$
must be of the form (\ref{special form for isometries}) as above.

	If $\delta$ is not equal to zero, then one would like to say that
$f$ is approximately equal to a rigid mapping when $\delta$ is small enough.
Here is a precise statement.  Let $D$ be a (closed) disk of radius $r$
in the plane.  This means that there is a point $w \in {\bf R}^2$ such that
\begin{equation}
	D = \{x \in {\bf R}^2 : |x-w| \le r \}.
\end{equation}
Then there is a rigid mapping $T : {\bf R}^2 \to {\bf R}^2$, depending
on $D$ and $f$, such that
\begin{equation}
\label{first approximation}
	r^{-1} \, \sup_{x \in D} |f(x) - T(x)| \le {\rm small}(\delta),
\end{equation}
where ${\rm small}(\delta)$ depends only on $\delta$, and not on $D$ or $f$,
and has the property that
\begin{equation}
\label{small(delta) to 0}
	{\rm small}(\delta) \to 0 \quad\hbox{as}\quad \delta \to 0.
\end{equation}
There are a number of ways to look at this.  One can give direct
constructive arguments, through basic geometric considerations or
computations.  In particular, one can derive explicit bounds for ${\rm
small}(\delta)$ in terms of $\delta$.  Results of this kind are given
in \cite{John}.  There are also abstract and inexplicit methods, in
which one argues by contradiction using compactness and the
Arzela--Ascoli theorem.  (In some related but different contexts, this
can be fairly easy or manageable, while explicit arguments and
estimates are less clear.)

	The presence of the factor of $r^{-1}$ on the left side of
(\ref{first approximation}) may not make sense at first glance, but it
is absolutely on target, and indispensable.  It reflects the natural
\emph{scaling} of the problem, and converts the left-hand side of
(\ref{first approximation}) into a dimensionless quantity, just as
$\delta$ is dimensionless.  One can view this in terms of the natural
\emph{invariances} of the problem.  Nothing changes here if we compose
$f$ (on either side) with a translation, rotation, or reflection, and
the same is true if we make simultaneous dilations on both the domain
and the range of equal amounts.  In other words, if $a$ is any
positive number, and if we define $f_a : {\bf R}^2 \to {\bf R}^2$ by
\begin{equation}
	f_a(x) = a^{-1} f(a x),
\end{equation}
then $f_a$ satisfies (\ref{1+delta bilipschitz}) exactly when $f$
does.  The approximation condition (\ref{first approximation}) is
formulated in such a way as to respect the same kind of invariances as
(\ref{1+delta bilipschitz}) does, and the factor of $r^{-1}$ accounts
for the dilation-invariance.

	This kind of approximation by rigid mappings is pretty good,
but can we do better?  Is it possible that the approximation works at
the level of the \emph{derivatives} of the mappings, rather than just
the mappings themselves?

	Here is another way to think about this, more directly in
terms of distance geometry.  Let us consider a simple mechanism by
which mappings that satisfy (\ref{1+delta bilipschitz}) can be
produced, and ask whether this mechanism gives everything.  Fix a
nonnegative number $k$, and call a mapping $g : {\bf R}^2 \to {\bf
R}^2$ is \emph{$k$-Lipschitz} if
\begin{equation}
\label{k-Lipschitz}
	|g(x) - g(y)| \le k \, |x-y|
\end{equation}
for all $x, y \in {\bf R}^2$.  This condition is roughly equivalent to
saying that the differential of $g$ has norm less than or equal to $k$
everywhere.  Specifically, if $g$ is differentiable at every point in
${\bf R}^2$, and if the norm of its differential is bounded by $k$
everywhere, then (\ref{k-Lipschitz}) holds, and this can be derived
from the mean value theorem.  The converse is not quite true, however,
because Lipschitz mappings need not be differentiable everywhere.
They are differentiable \emph{almost everywhere}, in the sense of
Lebesgue measure.  (See \cite{Fe, St1, Seapp}.)  To get a proper
equivalence one can consider derivatives in the sense of
distributions.

	If $f = S + g$, where $S$ is a rigid mapping and $g$ is
$k$-Lipschitz, and if $k \le 1/2$ (say), then $f$ satisfies
(\ref{1+delta bilipschitz}) with $\delta = 2k$.  (More precisely, one
can take $\delta = (1-k)^{-1} - 1$.)  This is not hard to check.  When
$k$ is small, this is a much stronger kind of approximation of $f$ by
rigid mappings than (\ref{first approximation}) is.  In particular, it
implies that the differential of $f$ is uniformly close to the
differential of $S$.

	To what extent can one go in the opposite direction, and say
that if $f$ satisfies (\ref{1+delta bilipschitz}) with $\delta$ small,
then $f$ can be approximated by rigid mappings in this stronger sense?
Let us begin by looking at what happens with the differential of $f$
at individual points.  Let $x$ be some point in ${\bf R}^2$, and assume
that the differential $df_x$ of $f$ at $x$ exists.  Thus $df_x$ is a linear
mapping from ${\bf R}^2$ to itself, and 
\begin{equation}
	f(x) + df_x(y-x)
\end{equation}
provides a good approximation to $f(y)$ for $y$ near $x$, in the
sense that
\begin{equation}
\label{approximation by the differential}
	|f(y) - \{f(x) + df_x(y-x)\}| = o(|y-x|).
\end{equation}
One can also think of the differential as the map obtained from $f$ by
``blowing up'' at $x$.  This corresponds to the formula
\begin{equation}
	df_x(v) = \lim_{t \to 0} t^{-1} (f(x+tv) - f(x)),
\end{equation}
with $t$ taken from positive real numbers.

	It is not hard to check that $df_x$, as a mapping on ${\bf
R}^2$ (with $x$ fixed), automatically satisfies (\ref{1+delta
bilipschitz}) when $f$ does.  Because the differential is already
linear, standard arguments from linear algebra imply that it is close
to a rotation or to the composition of a rotation and a reflection
when $\delta$ is small, and with easy and explicit estimates for the
degree of approximation.

	This might sound pretty good, but it is actually much weaker
than something like a representation of $f$ as $S + g$, where $S$ is a
rigid mapping and $g$ is $k$-Lipschitz with a reasonably-small value
of $k$.  If there is a representation of this type, then it means that
the differential $df_x$ of $f$ is always close to the differential of
$S$, which is \emph{constant}, i.e., independent of $x$.  The simple
method of the preceding paragraph implies that $df_x$ is always close
to being a rotation or a rotation composed with a reflection, but a
priori the choice of such a linear mapping could depend on $x$ in a
strong way.  That is very different from saying that there is a single
linear mapping that works for every $x$.

	Here is an example which shows how this sort of phenomenon can
happen.  (See also \cite{John}.)  Let us work in polar coordinates, so
that a point $z$ in ${\bf R}^2$ is represented by a radius $r \ge 0$
and an angle $\theta$.  We define $f : {\bf R}^2 \to {\bf R}^2$ by
saying that if $x$ is described by the polar coordinates $(r,\theta)$,
then
\begin{equation}
\label{example in polar coordinates}
	f(x) \enspace \hbox{has polar coordinates} \enspace 
				(r, \theta + \epsilon \, \log r).  
\end{equation}
Here $\epsilon$ is a small positive number that we get to choose.  Of
course $f$ should also take the origin to itself, despite the fact
that the formula for the angle degenerates there.

	Thus $f$ maps each circle centered at the origin to itself,
and on each such circle $f$ acts by a rotation.  We do not use a
single rotation for the whole plane, but instead let the rotation
depend logarithmically on the radius, as above.  This has the effect
of transforming every line through the origin into a logarithmic
spiral.  This spiral is very ``flat'' when $\epsilon$ is small, but
nonetheless it does wrap around the origin infinitely often in every
neighborhood of the origin.

	It is not hard to verify that this construction leads to a
mapping $f$ that satisfies (\ref{1+delta bilipschitz}), with a
$\delta$ that goes to $0$ when $\epsilon$ does, and at an easily
computable (linear) rate.  This gives an example of a mapping that
cannot be represented as $S + g$ with $S$ rigid and $g$ $k$-Lipschitz
for a fairly small value of $k$ (namely, $k < 1$).  For if $f$ did
admit such a representation, it would not be able to transform lines
into curves that spiral around a fixed point infinitely often; instead
it would take a line $L$ to a curve $\Gamma$ which can be realized as
the graph of a function over the line $S(L)$.  The spirals that we get
can never be realized as a graph of a function over any line.  This is
not hard to check.

	This spiralling is not incompatible with the kind of
approximation by rigid mappings in (\ref{first approximation}).  Let
us consider the case where $D$ is a disk centered at the origin, which
is the worst-case scenario anyway.  One might think that (\ref{first
approximation}) fails when we get too close to the origin (as compared
to the radius of $D$), but this is not the case.  Let $T$ be the
rotation on ${\bf R}^2$ that agrees with $f$ on the boundary of $D$.
If $\epsilon$ is small (which is necessary in order for the $\delta$
to be small in (\ref{1+delta bilipschitz})), then $T$ provides a good
approximation to $f$ on $D$ in the sense of (\ref{first
approximation}).  In fact, $T$ provides a good approximation to $f$ at
the level of their derivatives too on most of $D$, i.e., on the
complement of a small neighborhood of the origin.  The approximation of
derivatives breaks down near the origin, but the approximation of
values does not, as in (\ref{first approximation}), because $f$ and
$T$ both take points near the origin to points near the origin.

	This example suggests another kind of approximation by rigid
mappings that might be possible.  Given a disk $D$ of radius $r$ and a
mapping $f$ that satisfies (\ref{1+delta bilipschitz}), one would like
to have a rigid mapping $T$ on ${\bf R}^2$ so that (\ref{first
approximation}) holds, and also so that
\begin{equation}
\label{approximation in average}
	{1 \over \pi r^2} \int_D \|df_x - dT\| \, dx \le {\rm small}'(\delta),
\end{equation}
where ${\rm small}'(\delta)$ is, as before, a positive quantity which
depends only on $\delta$ (and not on $f$ or $D$) and which tends to
$0$ when $\delta$ tends to $0$.  Here $dx$ refers to the ordinary
$2$-dimensional integration against area on ${\bf R}^2$, and we think
of $df_x - dT$ as a matrix-valued function of $x$, with $\|df_x -
dT\|$ denoting its norm (in any reasonable sense).

	In other words, instead of asking that the differential of $f$
be approximated \emph{uniformly} by the differential of a rigid
mapping, which is not true in general, one can ask only that the
differential of $f$ be approximated by the differential of $T$
\emph{on average}.

	This does work, and in fact one can say more.  Consider the
expression
\begin{equation}
\label{definition of P(lambda)}
	P(\lambda) = {\rm Probability}(\{x \in D : \|df_x - dT\| \ge
				{\rm small}'(\delta) \cdot \lambda \}),
\end{equation}
where $\lambda$ is a positive real number.  Here ``probability'' means
Lebesgue measure divided by $\pi r^2$, which is the total measure of
the disk $D$.  The inequality (\ref{approximation in average}) implies
that
\begin{equation}
	P(\lambda) \le {1 \over \lambda}
\end{equation}
for all $\lambda > 0$.  It turns out that there is actually a
universal bound for $P(\lambda)$ with exponential decay for $\lambda
\ge 1$.  This was proved by John \cite{John} (with concrete
estimates).

	Notice that uniform approximation of the differential of $f$
by the differential of $T$ would correspond to a statement like
\begin{equation}
	P(\lambda) = 0
\end{equation}
for all $\lambda$ larger than some fixed (universal) constant.  John's
result of exponential decay is about the next best thing.

	As a technical point, let us mention that one can get
exponential decay conditions concerning the way that $\|df_x - dT\|$
should be small most of the time in a kind of trivial manner, with
constants that are not very good (at all), using the linear decay
conditions with good constants, together with the fact that $df$ is
bounded, so that $\|df_x - dT\|$ is bounded.  In the exponential decay
result mentioned above, the situation is quite different, and one
keeps constants like those from the linear decay condition.  This
comes out clearly in the proof, and we shall see more about it
later.

	This type of exponential decay occurs in a simple way in the
example above, in (\ref{example in polar coordinates}).  (This also
comes up in \cite{John}.)  One can obtain this from the presence of
$\epsilon \log r$ in the angle coordinate in the image.  The use of
the logarithm here is not accidental, but fits exactly with the
requirements on the mapping.  For instance, if one differentiates
$\log r$ in ordinary Cartesian coordinates, then one gets a quantity
of size $1/r$, and this is balanced by the $r$ in the first part of
the polar coordinates in (\ref{example in polar coordinates}), to give
a result which is bounded.

	It may be a bit surprising, or disappointing, that uniform
approximation to the differential of $f$ does not work here.  After
all, we did have ``uniform'' (or ``supremum'') bounds in the
hypothesis (\ref{1+delta bilipschitz}), and so one might hope to have
the same kind of bounds in the conclusion.  This type of failure of
supremum bounds is quite common, and in much the same manner as in the
present case.  We shall return to this in Section \ref{BMO frame of
mind}.

	How might one prove (\ref{approximation in average}), or the
exponential decay bounds for $P(\lambda)$?  Let us start with a
slightly simpler situation.  Imagine that we have a rectifiable curve
$\gamma$ in the plane whose total length is only slightly larger than
the distance between its two endpoints.  If the length of $\gamma$
were \emph{equal} to the distance between the endpoints, then $\gamma$
would have to be a straight line segment, and nothing more.  If the
length is slightly larger, then $\gamma$ has to stay close to the line
segment that joins its endpoints.  In analogy with (\ref{approximation
in average}), we would like to say that the tangents to $\gamma$ are
nearly parallel, on average, to the line that passes through the
endpoints of $\gamma$.

	In order to analyze this further, let $z(t)$, $t \in {\bf R}$,
$a \le t \le b$, be a parameterization of $\gamma$ by arclength.  This
means that $z(t)$ should be $1$-Lipschitz, so that
\begin{equation}
	|z(s) - z(t)| \le |s-t|
\end{equation}
for all $s, t \in [a,b]$, and that $|z'(t)| = 1$ for almost all $t$,
where $z'(t)$ denotes the derivative of $z(t)$.  Set
\begin{equation}
\label{definition of zeta}
	\zeta = {z(b) - z(a) \over b-a} = {1 \over b-a} \int_a^b z'(t) \, dt.
\end{equation}
Let us compute
\begin{equation}
	{1 \over b-a} \int_a^b | z'(s) - \zeta|^2 \, ds,
\end{equation}
which controls the average oscillation of $z'(s)$.  Let $\langle
\cdot, \cdot \rangle$ denote the standard inner product on ${\bf
R}^2$, so that
\begin{eqnarray}
	|x-y|^2 = \langle x-y, x-y \rangle 
		& = & \langle x, x \rangle - 2 \langle x, y \rangle 
						+ \langle y, y \rangle	\\
		& = & |x|^2 - 2 \langle x, y \rangle + |y|^2	\nonumber
\end{eqnarray}
for all $x, y \in {\bf R}^2$.  Applying this with $x = z'(s)$, $y =
\zeta$, we get that
\begin{equation}
\label{computation}
	{1 \over b-a} \int_a^b | z'(s) - \zeta|^2 \, ds
  	  = 1 - 2  {1 \over b-a} \int_a^b \langle z'(s), \zeta \rangle \, ds
			+ |\zeta|^2,
\end{equation}
since $|z'(s)| = 1$ a.e., and $\zeta$ does not depend on $s$.  The
middle term on the right side reduces to
\begin{equation}
	2 \langle \zeta, \zeta \rangle,
\end{equation}
because of (\ref{definition of zeta}).  Thus (\ref{computation}) yields
\begin{equation}
\label{result of computation}
	{1 \over b-a} \int_a^b | z'(s) - \zeta|^2 \, ds
		= 1 - 2 |\zeta|^2 + |\zeta|^2 = 1 - |\zeta|^2.
\end{equation}
On the other hand, $|z(b) - z(a)|$ is the same as the distance between
the endpoints of $\gamma$, and $b-a$ is the same as the length of
$\gamma$, since $z(t)$ is the parameterization of $\gamma$ by
arclength.  Thus $|\zeta|$ is exactly the ratio of the distance
between the endpoints of $\gamma$ to the length of $\gamma$, by
(\ref{definition of zeta}), and $1 - |\zeta|^2$ is a dimensionless
quantity which is small exactly when the length of $\gamma$ and the
distance between its endpoints are close to each other
(proportionately).  In this case (\ref{result of computation})
provides precise information about the way that $z'(s)$ is
approximately a constant on average.  (These computations follow ones
in \cite{CM2}.)

	One can use these results for curves for looking at mappings
from ${\bf R}^2$ (or ${\bf R}^n$) to itself, by considering images of
segments under the mappings.  This does not seem to give the proper
bounds in (\ref{approximation in average}), in terms of dependence on
$\delta$, though.  In this regard, see John's paper \cite{John}.
(Compare also with Appendix \ref{Some Fourier transform
calculations}.)  Note that for curves by themselves, the computations
above are quite sharp, as indicated by the equality in (\ref{result of
computation}).  See also \cite{CM2}.

	The exponential decay of $P(\lambda)$ requires more work.  A
basic point is that exponential decay bounds can be derived in a very
general way once one knows (\ref{approximation in average}) for
\emph{all} disks $D$ in the plane.  This is a famous result of John
and Nirenberg \cite{John-Nirenberg}, which will be discussed further
in Section \ref{BMO frame of mind}.  In the present situation, having
estimates like (\ref{approximation in average}) for all disks $D$ (and
with uniform bounds) is quite natural, and is essentially automatic,
because of the invariances of the condition (\ref{1+delta
bilipschitz}) under translations and dilations.  In other words, once
one has an estimate like (\ref{approximation in average}) for some
fixed disk $D$ and all mappings $f$ which satisfy (\ref{1+delta
bilipschitz}), one can conclude that the same estimate works for all
disks $D$, because of invariance under translations and dilations.

\section{The mathematics of good behavior much of the time, and the 
BMO frame of mind}
\label{BMO frame of mind}
\setcounter{equation}{0}

	Let us start anew for the moment, and consider the following
question in analysis.  Let $h$ be a real-valued function on ${\bf
R}^2$.  Let $\Delta$ denote the Laplace operator, given by
\begin{equation}
	\Delta = {\partial^2 \over \partial x_1^2} + 
			{\partial^2 \over \partial x_2^2},
\end{equation}
where $x_1$, $x_2$ are the standard coordinates on ${\bf R}^2$.  To what
extent does the behavior of $\Delta h$ control the behavior of the
other second derivatives of $h$?

	Of course it is easy to make examples where $\Delta h$
vanishes at a point but the other second derivatives do not vanish at
the same point.  Let us instead look for ways in which the overall
behavior of $\Delta h$ can control the overall behavior of the other
second derivatives.

	Here is a basic example of such a result.  Let us assume (for
simplicity) that $h$ is smooth and that it has compact support, and
let us write $\partial_1$ and $\partial_2$ for $\partial/\partial x_1$
and $\partial/\partial x_2$, respectively.  Then
\begin{equation}
\label{L^2 estimate}
	\int_{{\bf R}^2} |\partial_1 \partial_2 h(x)|^2 \, dx
		\le \int_{{\bf R}^2} |\Delta h(x)|^2 \, dx.
\end{equation}
This is a well-known fact, and it can be derived as follows.  We begin
with the identity
\begin{equation}
\label{identity}
	\int_{{\bf R}^2} \partial_1 \partial_2 h(x) \, 
				\partial_1 \partial_2 h(x) \, dx
	= \int_{{\bf R}^2} \partial_1^2 h(x) \, \partial_2^2 h(x) \, dx,
\end{equation}
which uses two integrations by parts.  On the other hand,
\begin{eqnarray}
\lefteqn{\int_{{\bf R}^2} |\Delta h(x)|^2 \, dx  = 
	  \int_{{\bf R}^2} (\partial_1^2 h(x) + \partial_2^2 h(x))^2 \, dx }
									\\
	&& \qquad
	    =  \int_{{\bf R}^2} (\partial_1^2 h(x))^2 
		+ 2 \, \partial_1^2 h(x) \, \partial_2^2 h(x)
		+ (\partial_2^2 h(x))^2 \, dx.			\nonumber
\end{eqnarray}
Combining this with (\ref{identity}) we get that
\begin{eqnarray}
\label{formula for the difference}
\lefteqn{\int_{{\bf R}^2} |\Delta h(x)|^2 \, dx   
		- 2 \int_{{\bf R}^2} |\partial_1 \partial_2 h(x)|^2 \, dx 
							} \\
	&& \qquad\qquad\qquad\qquad\qquad 
	      = \int_{{\bf R}^2} (\partial_1^2 h(x))^2 
		+ (\partial_2^2 h(x))^2 \, dx.			\nonumber
\end{eqnarray}
This implies (\ref{L^2 estimate}), and with an extra factor of $2$ on
the left-hand side, because the right side of (\ref{formula for the
difference}) is nonnegative.  (One can improve this to get a factor of
$4$ on the left side of (\ref{L^2 estimate}), using the right-hand
side of (\ref{formula for the difference}).)

	In short, the $L^2$ norm of $\Delta h$ always bounds the $L^2$
norm of $\partial_1 \partial_2 h$.  There are similar bounds for $L^p$
norms when $1 < p < \infty$.  Specifically, for each $p$ in
$(1,\infty)$, there is a constant $C(p)$ such that
\begin{equation}
\label{L^p estimate}
	\int_{{\bf R}^2} |\partial_1 \partial_2 h(x)|^p \, dx
		\le C(p) \int_{{\bf R}^2} |\Delta h(x)|^p \, dx
\end{equation}
whenever $h$ is a smooth function with compact support.  This is a
typical example of a ``Calder\'on--Zygmund inequality'', as in
\cite{St1}.  Such inequalities do \emph{not} work for $p = 1$ or
$\infty$, and the $p = \infty$ case is like the question of supremum
estimates in Section \ref{Mappings and distortion}.  Note that the $p
= 1$ and $p = \infty$ cases are closely connected to each other,
because of duality (of spaces and operators); the operators $\Delta$
and $\partial_1 \partial_2$ here are equal to their own transposes,
with respect to the standard bilinear form on functions on ${\bf R}^2$
(defined by taking the integral of the product of two given
functions).  In a modestly different direction, there are classical
results which give bounds in terms of the norm for H\"older continuous
(or Lipschitz) functions of order $\alpha$, for every $\alpha \in
(0,1)$, instead of the $L^p$ norm.  To be explicit, given $\alpha$,
this norm for a function $g$ on ${\bf R}^2$ can be described as the
smallest constant $A$ such that
\begin{equation}
	|g(x) - g(y)| \le A \, |x-y|^\alpha	
\end{equation}
for all $x, y \in {\bf R}^2$.  One can view this as a $p = \infty$
situation, like the $L^\infty$ norm for $g$, but with a positive order
$\alpha$ of smoothness, unlike $L^\infty$.  There is a variety of
other norms and spaces which one can consider, and for which there are
results about estimates along the lines of (\ref{L^p estimate}), but
for the norm in question instead of the $L^p$ norm.

	The $p=\infty$ version of (\ref{L^p estimate}) would say that
there is a constant $C$ such that
\begin{equation}
\label{p = infty version}
	\sup_{x \in {\bf R}^2} |\partial_1 \partial_2 h(x)|
		\le C \, \sup_{x \in {\bf R}^2} |\Delta h(x)|
\end{equation}
whenever $h$ is smooth and has compact support.  In order to see that
this is not the case, consider the function $h(x)$ given by
\begin{equation}
\label{example of a function h}
	h(x) = x_1 x_2 \log(x_1^2 + x_2^2),
\end{equation}
$x = (x_1, x_2)$.  It is not hard to compute $\Delta h$ and
$\partial_1 \partial_2 h$ explicitly, and to see that $\Delta h$ is
bounded while $\partial_1 \partial_2 h$ is not.  Indeed, 
\begin{equation}
\label{log singularity}
	\partial_1 \partial_2 h(x) = \log(x_1^2 + x_2^2) 
				+  \hbox{bounded terms},
\end{equation}
while the logarithm does not survive in $\Delta h$, because $\Delta
(x_1 x_2) \equiv 0$.

	This choice of $h$ is neither smooth nor compactly supported,
but these defects can be corrected easily.  For smoothness we can
consider instead
\begin{equation}
	h_\epsilon(x) = x_1 x_2 \log(x_1^2 + x_2^2 + \epsilon),
\end{equation}
where $\epsilon > 0$, and then look at what happens as $\epsilon \to
0$.  To make the support compact we can simply multiply by a fixed
cut-off function that does not vanish at the origin.  With these
modifications we still get a singularity at the origin as $\epsilon
\to 0$, and we see that (\ref{p = infty version}) cannot be true
(with a fixed constant $C$ that does not depend on $h$).

	This is exactly analogous to what happened in Section
\ref{Mappings and distortion}, i.e., with a uniform bound going in but
not coming out.  Instead of a uniform bound for the output, we also
have a substitute in terms of ``mean oscillation'', just as before.
To be precise, let $D$ be any disk in ${\bf R}^2$ of radius $r$, and
consider the quantity
\begin{equation}
	{1 \over \pi r^2} \int_D 
   |\partial_1 \partial_2 h(x) - {\rm Average}_D (\partial_1 \partial_2 h)|
			\, dx,
\end{equation}
where ``${\rm Average}_D \partial_1 \partial_2 h$'' is the average of
$\partial_1 \partial_2 h$ over the disk $D$, i.e.,
\begin{equation}
	{\rm Average}_D (\partial_1 \partial_2 h) =
		{1 \over \pi r^2} \int_D \partial_1 \partial_2 h(u) \, du.
\end{equation}
Instead of (\ref{p = infty version}), it is true that there is a constant
$C > 0$ so that 
\begin{equation}
\label{mean oscillation bound, first version}
	{1 \over \pi r^2} \int_D 
   |\partial_1 \partial_2 h(x) - {\rm Average}_D (\partial_1 \partial_2 h)|
			\, dx
		\le C \, \sup_{x \in {\bf R}^2} |\Delta h(x)|
\end{equation}
for every disk $D$ in ${\bf R}^2$ of radius $r$ and every smooth
function $h$ with compact support.  This is not too hard to prove;
roughly speaking, the point is to ``localize'' the $L^2$ estimate that
we had before.  (More general results of this nature are discussed
in \cite{GCRF, Journe, St2}.)

	Let us formalize this estimate by defining a new space of
functions, namely the space BMO of functions of \emph{bounded mean
oscillation}, introduced by John and Nirenberg in
\cite{John-Nirenberg}.  A locally-integrable function $g$ on ${\bf
R}^2$ is said to lie in BMO if there is a nonnegative number $k$ such
that
\begin{equation}
\label{def of bmo}
	{1 \over \pi r^2} \int_D |g(x) - {\rm Average}_D (g)| \, dx
		\le k
\end{equation}
for every disk $D$ in ${\bf R}^2$ of radius $r$.  In this case we set
\begin{equation}
	\|g\|_* = \sup_D 
	{1 \over \pi r^2} \int_D |g(x) - {\rm Average}_D (g)| \, dx,
\end{equation}
with the supremum taken over all disks $D$ in ${\bf R}^2$.  This is
the same as saying that $\|g\|_*$ is the \emph{smallest} number $k$
that satisfies (\ref{def of bmo}).  One refers to $\|g\|_*$ as the
``BMO norm of $g$'', but notice that $\|g\|_*=0$ when $g$ is equal to
a constant almost everywhere.  (The converse is also true.)

	This definition may look a little crazy, but it works quite
well in practice.  Let us reformulate (\ref{mean oscillation bound,
first version}) by saying that there is a constant $C$ so that
\begin{equation}
	\|\partial_1 \partial_2 h\|_* \le C \, \| \Delta h \|_\infty,
\end{equation}
where $\|\phi\|_\infty$ denotes the $L^\infty$ norm of a given
function $\phi$.  In other words, although the $L^\infty$ norm of
$\partial_1 \partial_2 h$ is not controlled (for all $h$) by the
$L^\infty$ norm of $\Delta h$, the BMO norm of $\partial_1 \partial_2
h$ is controlled by the $L^\infty$ norm of $\Delta h$.

	Similarly, one of the main points in Section \ref{Mappings and
distortion} can be reformulated as saying that if a mapping $f : {\bf
R}^2 \to {\bf R}^2$ distorts distances by only a small amount, as in
(\ref{1+delta bilipschitz}), then the BMO norm $\|df\|_*$ of the
differential of $f$ is small (and with precise estimates being
available).

	In Section \ref{Mappings and distortion} we mentioned a
stronger estimate with exponential decay in the measure of certain
``bad'' sets.  This works for all BMO functions, and can be given as
follows.  Suppose that $g$ is a BMO function on ${\bf R}^2$ with
$\|g\|_* \le 1$, and let $D$ be a disk in ${\bf R}^2$ with radius $r$.
As in (\ref{definition of P(lambda)}), consider the ``distribution
function'' $P(\lambda)$ defined by
\begin{equation}
\label{new definition of P(lambda)}
	P(\lambda) = {\rm Probability}(\{x \in D : 
			|g(x) - {\rm Average}_D (g)| \ge \lambda \}),
\end{equation}
where ``Probability'' means Lebesgue measure divided by the area $\pi
r^2$ of $D$.  Under these conditions, there is a universal bound for
$P(\lambda)$ with exponential decay, i.e., an inequality of the form
\begin{equation}
\label{exponential decay}
	P(\lambda) \le B^{- \lambda} \qquad\hbox{for } \lambda \ge 1,
\end{equation}
where $B$ is a positive number greater than $1$, and $B$ does not
depend on $g$ or $D$.  This is a theorem of John and Nirenberg
\cite{John-Nirenberg}.

	Although we have restricted ourselves to ${\bf R}^2$ here for
simplicity, everything goes over in a natural way to Euclidean spaces
of arbitrary dimension.  In fact, there is a much more general
framework of ``spaces of homogeneous type'' in which basic properties
of BMO (and other aspects of real-variable harmonic analysis) carry
over.  See \cite{CW1, CW2}, and compare also with \cite{GCRF, St2}.
This framework includes certain Carnot spaces that arise in several
complex variables, like the unit sphere in ${\bf C}^n$ with the
appropriate (noneuclidean) metric.

	The exponential decay bound in (\ref{exponential decay}) helps
to make precise the idea that BMO functions are very close to being
bounded (which would correspond to having $P(\lambda) = 0$ for all
sufficiently large $\lambda$).  The exponential rate of decay implies
that BMO functions lie in $L^p$ locally for all finite $p$, but it
is quite a bit stronger than that.

	A basic example of a BMO function is $\log |x|$.  This is not
hard to check, and it shows that exponential decay in
(\ref{exponential decay}) is sharp, i.e., one does not have
superexponential decay in general.  This example also fits with
(\ref{log singularity}), and with the ``rotational'' part of the
differential of the mapping $f$ in (\ref{example in polar
coordinates}).

	In general, BMO functions can be much more complicated than
the logarithm.  Roughly speaking, the total ``size'' of the
unboundedness is no worse than for the logarithm, as in
(\ref{exponential decay}), but the arrangement of the singularities
can be more intricate, just as one can make much more complex singular
examples than in (\ref{example of a function h}) and (\ref{example in
polar coordinates}).  There are a lot of tools available in harmonic
analysis for understanding exactly how BMO functions behave.  (See
\cite{GCRF, Ga, Journe, St2}, for instance.)

	BMO functions show up all over the place.  One can reformulate
the basic scenario in this section with the Laplacian and $\partial_1
\partial_2$ by saying that the pseudodifferential or singular integral
operator
\begin{equation}
	{\partial_1 \partial_2 \over \Delta}
\end{equation}
maps $L^\infty$ to BMO, and this holds for similar operators (of order
$0$) much more generally (as in \cite{GCRF, Ga, Journe, St2}).  This will
be discussed a bit further in Appendix \ref{Some Fourier transform
calculations}.  Note that the nonlinear problem in Section
\ref{Mappings and distortion} has a natural linearization which falls
into this rubric.  (See Appendix \ref{Some Fourier transform
calculations}.)

	Sobolev embeddings provide another class of linear problems in
which BMO comes up naturally.  One might wish that a function $g$ on
${\bf R}^n$ that satisfies $\nabla g \in L^n({\bf R}^n)$ (in the sense
of distributions) were bounded or continuous, but neither of these are
true in general, when $n > 1$.  However, such a function $g$ is always
in BMO, and in the subspace VMO (``vanishing mean oscillation''), in
which the measurements of mean oscillation (as in the left side of
(\ref{def of bmo}) when $n=2$) tend to $0$ as the radius $r$ goes to
$0$.  This is a well-known analogue of continuity in the context of
BMO.  (See \cite{BN, GCRF, Ga, Seapp, St2}.)

	BMO arises in a lot of nonlinear problems, in addition to the
one in Section \ref{Mappings and distortion}.  For instance, there are
circumstances in which one might wish that the derivative of a
conformal mapping in the complex plane were bounded, and it is not,
but there are natural estimates in terms of BMO.  More precisely, it
is BMO for the \emph{logarithm} of the derivative that comes up most
naturally.  This is closely related to BMO conditions for tangents to
curves under certain geometric conditions.  See \cite{CM1, CM2, CM3,
Da-3rdcycle, JK1, P1, P2, P3}, for instance.  Some basic computations
related to the latter were given in Section \ref{Mappings and
distortion}, near the end.  In general dimensions (larger than $1$),
BMO shows up naturally as the logarithm of the density for harmonic
measure for Lipschitz domains, and for the logarithm of Jacobians of
quasiconformal mappings.  See \cite{Dh1, Dh2, JK2, Ge, Re, St2} and
the references therein.  In all dimensions, there are interesting
classes of ``weights'', positive functions which one can use as
densities for modifications of Lebesgue measure, whose logarithms lie
in BMO, and which in fact correspond to open subsets of BMO (for
real-valued functions).  These weights have good properties concerning
$L^p$ boundedness of singular integral and other operators, and they
also show up in other situations, in connection with conformal
mappings in the plane, harmonic measure, and Jacobians of
quasiconformal mappings in particular, as above.  See \cite{GCRF, Ga,
Journe, St2, ST} for information about these classes of weights.

	There is a simple reason for BMO functions to arise frequently
as some kind of logarithm.  In many nonlinear problems there is a
symmetry which permits one to multiply some quantity by a constant
without changing anything in a significant way.  (E.g., think of
rescaling or rotating a domain, or a mapping, or multiplying a weight
by a positive constant.)  At the level of the logarithm this
invariance is converted into a freedom to add constants, and this is
something that BMO accommodates automatically.

	To summarize a bit, there are a lot of situations in which one
has some function that one would like to be bounded, but it is not,
and for which BMO provides a good substitute.  One may not expect at
first to have to take measure theory into account, but then it comes
up on its own, or works in a natural or reasonable way. 

	Before leaving this section, let us return to the
John--Nirenberg theorem, i.e., the exponential decay estimate in
(\ref{exponential decay}).  How might one try to prove this?  The
first main point is that one cannot prove (\ref{exponential decay})
for a particular disk $D$ using only a bound like (\ref{def of bmo})
for that one disk.  That would only give a rate of decay on the order
of $1/\lambda$.  Instead one uses (\ref{def of bmo}) over and over
again, for many different disks.

	Here is a basic strategy.  Assume that $g$ is a BMO function
with $\|g\|_* \le 1$.  First use (\ref{def of bmo}) for $D$ itself
(with $k = 1$) to obtain that the set of points $x$ in $D$ such that
\begin{equation}
	|g(x) - {\rm Average}_D (g)| \ge 10,
\end{equation}
is pretty small (in terms of probability).  On the bad set where this
happens, try to make a good covering by smaller disks on which one can
apply the same type of argument.  The idea is to then show that the
set of points $x$ in $D$ which satisfy
\begin{equation}
	|g(x) - {\rm Average}_D (g)| \ge 10 + 10
\end{equation}
is significantly smaller still, and by a definite proportion.  If one
can repeat this forever, then one can get exponential decay as in
(\ref{exponential decay}).  More precisely, at each stage the size of
the deviation of $g(x)$ from ${\rm Average}_D (g)$ will increase by
the \emph{addition} of $10$, while the decrease in the measure of the
bad set will decrease \emph{multiplicatively}.

	This strategy is roughly correct in spirit, but to carry it
out one has to be more careful in the choice of ``bad'' set at each
stage, and in the transition from one stage to the next.  In
particular, one should try to control the difference between the
average of $g$ over one disk and over one of the smaller disks created
in the next step of the process.  As a practical matter, it is simpler
to work with cubes instead of disks, for the way that they can be
decomposed evenly into smaller pieces.  The actual construction used
is the ``Calder\'on--Zygmund decomposition'', which itself has a lot
of other applications.  See \cite{John-Nirenberg, GCRF, Ga, Journe,
St2, Seapp} for more information.

\section{Finite polyhedra and combinatorial parameterization problems}
\label{parameterization problems}
\setcounter{equation}{0}

	Let us now forget about measure theory for the time being, and look
at a problem which is, in principle, purely combinatorial.

	Fix a positive integer $d$, and let $P$ be a $d$-dimensional
polyhedron.  We assume that $P$ is a finite union of $d$-dimensional
simplices, so that $P$ has ``pure'' dimension $d$ (i.e., with no
lower-dimensional pieces sticking off on their own).  

\beginproblem
\label{recognition problem}
How can one tell if $P$ is a PL (piecewise-linear) manifold?  In other words,
when is $P$ locally PL-equivalent to ${\bf R}^d$ at each point?
\end{problem}

	To be precise, $P$ is locally PL-equivalent to ${\bf R}^d$ at
a point $x \in P$ if there is a neighborhood of $x$ in $P$ which is
homeomorphic to an open set in ${\bf R}^d$ through a mapping which is
piecewise-linear.

	This is really just a particular example of a general issue,
concerning existence and complexity of parameterizations of a given
set.  Problem \ref{recognition problem} has the nice feature that
finite polyhedra and piecewise-linear mappings between them can, in
principle, be described in finite terms.

	Before we try to address Problem \ref{recognition problem}
directly, let us review some preliminary matters.  It will be
convenient to think of $P$ as being like a simplicial complex, so that
it is made up of simplices which are always either disjoint or meet in
a whole face of some (lower) dimension.  Thus we can speak about the
vertices of $P$, the edges, the $2$-dimensional faces, and so on, up
to the $d$-dimensional faces.

	Since $P$ is a finite polyhedron, its local structure at any
point is pretty simple.  Namely, $P$ looks like a cone over a
$(d-1)$-dimensional polyhedron at every point.  To make this precise,
imagine that $Q$ is some finite polyhedron in some ${\bf R}^n$, and
let $z$ be a point in ${\bf R}^n$ which is affinely-independent of
$Q$, i.e., which lies in the complement of an (affine) plane that
contains $Q$.  (We can always replace ${\bf R}^n$ with ${\bf
R}^{n+1}$, if necessary, to ensure that there is such a point.)  Let
$c(Q)$ denote the set which consists of all rays in ${\bf R}^n$ which
emanate from $z$ and pass through an element of $Q$.  We include $z$
itself in each of these rays.  This defines the ``cone over $Q$
centered at $z$''.  It does not really depend on the choice of $z$, in
the sense that a different choice of $z$ leads to a set which is
equivalent to the one just defined through an invertible affine
transformation.

	If $x$ is a ``vertex'' of $P$, in the sense described above,
then there is a natural way to choose a $(d-1)$-dimensional polyhedron
$Q$ so that $P$ is the same as the cone over $Q$ centered at $x$ in a
neighborhood of $x$.  Let us call $Q$ the \emph{link} of $P$ at $x$.
(Actually, with this description $Q$ is only determined up to
piecewise-linear equivalence, but this is adequate for our purposes.)

	Now suppose that $x$ is not a vertex.  One can still realize
$P$ as a cone over a $(d-1)$-dimensional polyhedron near $x$, but one
can also do something more precise.  If $x$ is not a vertex, then
there is a positive integer $k$ and a $k$-dimensional face $F$ of $P$
such that $x$ lies in the interior of $F$.  In this case there is a
$(d-k-1)$-dimensional polyhedron $Q$ such that $P$ is locally
equivalent to ${\bf R}^k \times c(Q)$ near $x$, with $x$ in $P$
corresponding to a point $(y,z)$ in ${\bf R}^k \times c(Q)$, where $z$
is the center of $c(Q)$.  This same polyhedron $Q$ works for all the
points in the interior of $F$, and we call $Q$ the \emph{link} of $F$.

\beginbasicfact
\label{characterization in terms of spheres}
$P$ is everywhere locally equivalent to ${\bf R}^d$ if and only if
all of the various links of $P$ (of all dimensions) are piecewise-linearly
equivalent to standard spheres (of the same dimension).  
\end{basicfact}

	Here the ``standard sphere of dimension $m$'' can be taken to
be the boundary of the standard $(m+1)$-dimensional simplex.  

	Basic Fact \ref{characterization in terms of spheres} is
standard and not hard to see.  The ``if'' part is immediate, since one
knows exactly what the cone over a standard sphere looks like, but for
the converse there is a bit more to check.  A useful observation is
that if $Q$ is a $j$-dimensional polyhedron whose cone $c(Q)$ is
piecewise-linearly equivalent to ${\bf R}^{j+1}$ in a neighborhood of
the center of $c(Q)$, then $Q$ must be piecewise-linearly equivalent
to a standard $j$-dimensional sphere.  This is pretty easy to verify,
and one can use it repeatedly for the links of $P$ of codimension
larger than $1$.  (A well-known point here is that one should be
careful not to use \emph{radial} projections to investigate links
around vertices, but suitable \emph{pseudo-radial} projections, to fit
with the piecewise-linear structure, and not just the topological
structure.)

	A nice feature of Basic Fact \ref{characterization in terms of
spheres} is that it sets up a natural induction in the dimensions,
since the links of $P$ always have dimension less than $P$.  This
leads to the following question.

\beginproblem
\label{recognizing PL spheres}
If $Q$ is a finite polyhedron which is a $k$-dimensional PL manifold,
how can one tell if $Q$ is a PL sphere of dimension $k$?
\end{problem}

	It is reasonable to assume here that $Q$ is a PL-manifold,
because of the way that one can use Basic Fact \ref{characterization
in terms of spheres} and induction arguments.

	Problem \ref{recognizing PL spheres} is part of the matter of
the Poincar\'e conjecture, which would seek to say that $Q$ is a PL
sphere as soon as it is homotopy-equivalent to a sphere.  This has
been established in all dimensions except $3$ and $4$.  (Compare with
\cite{Rourke-Sanderson}.)  In dimension $4$ the Poincar\'e conjecture
was settled by M. Freedman \cite{Fr} in the ``topological'' category
(with ordinary homeomorphisms (continuous mappings with continuous
inverses) and topological manifolds), but it remains unknown in the PL
case.  The PL case is equivalent to the smooth version in this
dimension, and both are equivalent to the ordinary topological version
in dimension $3$.  (A brief survey related to these statements is
given in Section 8.3 of \cite{FQ}.)  Although the Poincar\'e
conjecture is known to hold in the PL category in all higher
dimensions (than $4$), it does not always work in the smooth category,
because of exotic spheres (as in \cite{Milnor-exotic-spheres,
Kervaire-Milnor}).
	
	If the PL version of the Poincar\'e conjecture is true in all
dimensions, then this would give one answer to the question of
recognizing PL manifolds among finite polyhedra in Problem
\ref{recognition problem}.  Specifically, our polyhedron $P$ would be
a PL manifold if and only if its links are all homotopy-equivalent to
spheres (of the correct dimension).  

	This might seem like a pretty good answer, but there are
strong difficulties concerning complexity for matters of homotopy.  In
order for a $k$-dimensional polyhedron $Q$ to be a homotopy sphere, it
has to be simply connected in particular, at least when $j \ge 2$.  In
other words, it should be possible to continuously deform any loop in
$Q$ to a single point, or, equivalently, to take any continuous
mapping from a circle into $Q$ and extend it to a continuous mapping
from a closed disk into $Q$.  This extension can entail enormous
complexity, in the sense that the filling to the disk might have to be
of much greater complexity than the original loop itself.

	This is an issue whose geometric significance is often
emphasized by Gromov.  To describe it more precisely it is helpful to
begin with some related algebraic problems, concerning
finitely-presented groups.

	Let $G$ be a group.  A finite presentation of $G$ is given by
a finite list $g_1, g_2, \ldots, g_n$ of generators for $G$ together
with a finite set $r_1, r_2, \ldots, r_m$ of ``relations''.  The
latter are (finite) words made out of the $g_i$'s and their inverses.
Let us assume for convenience that the set of relations includes the
inverses of all of its elements, and also the empty word.  The $r_j$'s
are required to be trivial, in the sense that they represent the
identity element of $G$.  This implies that arbitrary products of
conjugates of the $r_j$'s also represent the identity element, and the
final requirement is that if $w$ is any word in the $g_i$'s and their
inverses which represents the identity element in $G$, then it should
be possible to obtain $w$ from some product of conjugates of the
$r_j$'s through cancellations of subwords of the form $g_i^{-1} g_i$
and $g_i g_i^{-1}$.

	For instance, the group ${\bf Z}^2$ can be described by two
generators $a$, $b$ and one relation, $aba^{-1}b^{-1}$.  As another
concrete example, there is the (Baumslag--Solitar) group with two
generators $x$, $y$ and one relation $x^2 y x^{-1} y^{-1}$.

	Suppose that a group $G$ and finite presentation of $G$ are
given and fixed, and let $w$ be a word in the generators of $G$ and
their inverses.  Given this information, how can one decide whether
$w$ represents the identity element in $G$?  This is called ``the word
problem'' (for $G$).  It is a famous result that there exist finite
presentations of groups for which there is no algorithm to solve the
word problem.  (See \cite{Manin}.)

	To understand what this really means, let us first notice that
the set of trivial words for the given presentation is ``recursively
enumerable''.  This means that there is an algorithm for listing all
of the trivial words.  To do this, one simply has to have the
algorithm systematically generate all possible conjugates of the
relations, all possible products of conjugates of relations, and all
possible words derived from these through cancellations as above.  In
this way the algorithm will constantly generate trivial words, and
every trivial word will eventually show up on the list.

	However, this does not give a finite procedure for determining
that a given word is \emph{not} trivial.  A priori one cannot conclude
that a given word is not trivial until one goes through the entire
list of trivial words.

	The real trouble comes from the cancellations.  In order to
establish the triviality of a given word $w$, one might have to make
derivations through words which are enormously larger, with a lot of
collapsing at the end.  If one had a bound for the size of the words
needed for at least one derivation of the triviality of a given word
$w$, a bound in terms of an effectively computable (or ``recursive'')
function of the length of $w$, then the word problem would be
algorithmically solvable.  One could simply search through all
derivations of at most a computable size.

	This would not be very efficient, but it would be an
algorithm.  As it is, even this does not always work, and there are
finitely-presented groups for which the derivations of triviality may
need to involve words of nonrecursive size compared to the given word.

	One should keep in mind that for a given group and a given
presentation there is always \emph{some} function $f(n)$ on the
positive integers so that trivial words of length at most $n$ admit
derivations of their triviality through words of size no greater than
$f(n)$.  This is true simply because there are only finitely many
words of size at most $n$, and so one can take $f(n)$ to be the
maximum size incurred in some finite collection of derivations.  The
point is that such a function $f$ may not be bounded by a recursive
function.  This means that $f$ could be really huge, larger than any
tower of exponentials, for instance.

	The same kind of phenomenon occurs geometrically, for deciding
whether a loop in a given polyhedron can be continuously contracted to
a point.  This is because any finite presentation of a group $G$ can
be coded into a finite polyhedron, in such a way that the group $G$ is
represented by the fundamental group of the polyhedron.  This is a
well-known construction in topology.

	Note that while the fundamental group of a space is normally
defined in terms of \emph{continuous} (based) loops in the space and
the continuous deformations between them, in the case of finite
polyhedra it is enough to consider \emph{polygonal} loops and
deformations which are piecewise-linear (in addition to being
continuous).  This is another standard fact, and it provides a
convenient way to think about complexity for loops and their
deformations.

	Although arbitrary finite presentations can be coded into
finite polyhedra, as mentioned above, this is not the same as saying
that they can be coded into compact \emph{manifolds}.  It turns out
that this does work when the dimension is at least $4$, i.e., for each
$n \ge 4$ it is true that every finite presentation can be coded into
a compact PL manifold of dimension $n$.  This type of coding can be
used to convert algorithmic unsolvability results for problems in
group theory into algorithmic unsolvability statements in topology.
For instance, there does not exist an algorithm to decide when a given
finite presentation for a group actually defines the trivial group,
and, similarly, there does not exist an algorithm for deciding when a
given manifold (of dimension at least $4$) is simply-connected.  See
\cite{BHP, Markov1, Markov2, Markov3} for more information and
results.

	Let us mention that in dimensions $3$ and less, it is
\emph{not} true that arbitrary finitely-presented groups can be
realized as fundamental groups of compact manifolds.  Fundamental
groups of manifolds are very special in dimensions $1$ and $2$, as is
well known.  The situation in dimension $3$ is more complicated, but
there are substantial restrictions on the groups that can arise as
fundamental groups.  As an aspect of this, one can look at
restrictions related to Poincar\'e duality.  In a different vein, the
fundamental group of a $3$-dimensional manifold has the property that
all of its finitely-generated subgroups are finitely-presented.  See
\cite{Scott}, and Theorem 8.2 on p70 of \cite{Hempel1}.  See also
\cite{Jaco}.  In another direction, there are relatively few abelian
groups which can arise as subgroups of fundamental groups of
$3$-dimensional manifolds.  See \cite{Epstein, Evans-Moser}, Theorems
9.13 and 9.14 on p84f of \cite{Hempel1}, and p67-9 of \cite{Jaco}.  At
any rate, it is a large open problem to know exactly what groups arise
as fundamental groups of $3$-dimensional manifolds.

	See also \cite{Thurston-Levy} and Chapter 12 of
\cite{ECHLPT} concerning these groups.  The book \cite{ECHLPT} treats
a number of topics related to computability and groups, and not just
in connection with fundamental groups of $3$-manifolds.  This includes
broad classes of groups for which positive results and methods are
available.  See \cite{Farb} as well in this regard.

	Beginning in dimension $5$, it is known that there is no
algorithm for deciding when a compact PL manifold is
piecewise-linearly equivalent to a standard (PL) sphere.  This is a
result of S. Novikov.  See Section 10 of \cite{VKF}, and also the
appendix to \cite{Na}.  (Note that in dimensions less than or equal to
$3$, such algorithms do exist.  This is classical for dimensions $1$,
$2$; see \cite{Rubinstein1, Rubinstein2, Thompson} concerning
dimension $3$, and related problems and results.)  Imagine that we
have a PL manifold $M$ of some dimension $n$ whose equivalence to a
standard sphere is true but ``hard'' to check.  According to the
solution of the Poincar\'e conjecture in these dimensions, $M$ will be
equivalent to an $n$-sphere if it is homotopy-equivalent to ${\bf
S}^n$.  For standard reasons of algebraic topology, this will happen
exactly when $M$ is simply-connected and has trivial homology in
dimensions $2$ through $n-1$.  (Specifically, this uses Theorem 9 and
Corollary 24 on pages 399 and 405, respectively, of \cite{Sp}.  It
also uses the existence of a degree-$1$ mapping from $M$ to ${\bf
S}^n$ to get started (i.e., to have a mapping to which the
aforementioned results can be applied), and the fact that the homology
of $M$ and ${\bf S}^n$ vanish in dimensions larger than $n$, and are
equal to ${\bf Z}$ in dimension $n$.  To obtain the degree-$1$ mapping
from $M$ to ${\bf S}^n$, one can start with any point in $M$ and a
neighborhood of that point which is homeomorphic to a ball.  One then
collapses the complement of that neighborhood to a point, which gives
rise to the desired mapping.)  The vanishing of homology can be
determined algorithmically, and so if the equivalence of $M$ with an
$n$-sphere is ``hard'' for algorithmic verification, then the problem
must occur already with the simple-connectivity of $M$.  (Concerning
this statement about homology, see Appendix \ref{A few facts about
homology}.)

	To determine whether $M$ is simply-connected it is enough to
check that a finite number of loops in $M$ can be contracted to a
point, i.e., some collection of generators for the fundamental group.
If this is ``hard'', then it means that the complexity of the
contractions should be enormous compared to the complexity of $M$.
For if there were a bound in terms of a recursive function, then one
could reverse the process and use this to get an algorithm which could
decide whether $M$ is PL equivalent to a sphere, and this is not
possible.

	If $M$ is a hard example of a PL manifold which is equivalent
to an $n$-sphere, then any mapping from $M$ to the sphere which
realizes this equivalence must necessarily be of very high complexity
as well.  Because of the preceding discussion, this is also true for
mappings which are homotopy-equivalences, or even which merely induce
isomorphisms on $\pi_1$, if one includes as part of the package of
data enough information to justify the condition that the induced
mapping on $\pi_1$ be an isomorphism.  (For a homotopy equivalence,
for instance, one could include the mapping $f$ from $M$ to the
$n$-sphere, a mapping $g$ from the $n$-sphere to $M$ which is a
homotopy-inverse to $f$, and mappings which give homotopies between $f
\circ g$ and $g \circ f$ to the identity on the $n$-sphere and $M$,
respectively.)  This is because one could use the mapping to reduce
the problem of contracting a loop in $M$ to a point to the
corresponding problem for the $n$-sphere, where the matter of bounds
is straightforward.

	Similar considerations apply to the problem of deciding when a
finite polyhedron $P$ is a PL manifold.  Indeed, given a PL manifold
$M$ whose equivalence to a sphere is in question, one can use it to
make a new polyhedron $P$ by taking the ``suspension'' of $M$.  This
is defined by taking two points $y$ and $z$ which lie outside of a
plane that contains $M$, and then taking the union of all of the
(closed) line segments that go from either of $y$ or $z$ to a point in
$M$.  One should also be careful to choose $y$ and $z$ so that these
line segments never meet, except in the trivial case of line segments
from $y$ and $z$ to the same point $x$ in $M$, with $x$ being the only
point of intersection of the two segments.  (One can imagine $y$ and
$z$ as lying on ``opposite sides'' of an affine plane that contains
$M$.)

	If $M$ is equivalent to a sphere, then this operation of
suspension produces a PL manifold equivalent to the sphere of $1$
larger dimension, as one can easily check.  If $M$ is not PL
equivalent to a sphere, then the suspension $P$ of $M$ is not a PL
manifold at all.  This is because $M$ is the link of $P$ at the
vertices $y$ and $z$, by construction, so that one is back to the
situation of Basic Fact \ref{characterization in terms of spheres}.

	Just as there are PL manifolds $M$ whose equivalence with a
sphere is hard, the use of the suspension shows that there are
polyhedra $P$ for which the property of being a PL manifold is hard to
establish.  Through the type of arguments outlined above, when PL
coordinates exist for a polyhedron $P$, they may have to be of
enormous complexity compared to the complexity of $P$ itself.  This
works more robustly than just for PL coordinates, i.e., for any
information which is strong enough to give the simple-connectivity of
the links of $P$.  Again, this follows the discussion above.

	We have focussed on piecewise-linear coordinates for finite
polyhedra for the sake of simplicity, but similar themes of complexity
come up much more generally, and in a number of different ways.  In
particular, existence and complexity of parameterizations is often
related in a strong manner to the behavior of something like $\pi_1$,
sometimes in a localized form, as with the links of a polyhedron.  For
topology of manifolds in high dimensions, $\pi_1$ and the filling of
loops with disks comes up in the Whitney lemma, for instance.  This
concerns the separation of crossings of submanifolds through the use
of embedded $2$-dimensional disks, and it can be very useful for
making some geometric constructions.  (A very nice brief review of
some of these matters is given in Section 1.2 of \cite{DoK}.)
Localized $\pi_1$-type conditions play a crucial role in taming
theorems in geometric topology.  Some references related to this are
\cite{Bing5, Bing-book, Bing-coll, Bu, BuC, C1, C-bulletin, Dm1, Dm2,
E-demension, Moise, Ru-book, Ru}.

	In Appendix \ref{more on existence and behavior of
homeomorphisms}, we shall review some aspects of geometric topology
and the existence and behavior of parameterizations, and the role of
localized versions of fundamental groups in particular.

	As another type of example, one has the famous ``double
suspension'' results of Edwards and Cannon \cite{C1, C2, Dm2, E}.
Here one starts with a finite polyhedron $H$ which is a manifold with
the same homology as a sphere of the same dimension, and one takes the
suspension (described above) of the suspension of $H$ to get a new
polyhedron $K$.  The result is that $K$ is actually homeomorphic to a
sphere.  A key point is that $H$ is \emph{not} required to be
simply-connected.  When $\pi_1(H) \ne 0$, it is not possible for the
homeomorphism from $K$ to a standard sphere to be piecewise-linear, or
even Lipschitz (as in (\ref{k-Lipschitz})).  Concerning the latter,
see \cite{SS}.  Not much is known about the complexity of the
homeomorphisms in this case.  (We shall say a bit more about this in
Section \ref{Uniform rectifiability} and Subsection \ref{Manifold
factors}.)

	Note that if $J$ is obtained as a single suspension of $H$,
and if $\pi_1(H) \ne 0$, then $J$ cannot be a topological manifold at
all (at least if the dimension of $H$ is at least $2$).  Indeed, if
$M$ is a topological manifold of dimension $n$, then for every point
$p$ in $M$ there are arbitrarily small neighborhoods $U$ of $p$ which
are homeomorphic to an open $n$-ball, and $U \backslash \{p\}$ must
then be simply-connected when $n \ge 3$.  This cannot work for the
suspension $J$ of $H$ when $\pi_1(H) \ne 0$, with $p$ taken to be one
of the two cone points introduced in the suspension construction.

	However, $J$ has the advantage over $H$ that it \emph{is}
simply-connected.  This comes from the process of passing to the
suspension (and the fact that $H$ should be connected, since it has
the same homology as a sphere).  It is for this reason that the cone
points of $K$ do not have the same trouble as in $J$ itself, with no
small deleted neighborhoods which are simply-connected.  The
singularities at the cone points in $J$ lead to trouble with the
codimension-$2$ links in $K$, but this turns out not to be enough to
prevent $K$ from being a topological manifold, or a topological
sphere.  It \emph{does} imply that the homeomorphisms involved have to
distort distances in a very strong way, as in \cite{SS}.

	In other words, local homeomorphic coordinates for $K$ do
exist, but they are necessarily much more complicated than PL
homeomorphisms, even though $K$ is itself a finite polyhedron.  As
above, there is also a global homeomorphism from $K$ to a sphere.  The
first examples of finite polyhedra which are homeomorphic to each
other but not piecewise-linearly equivalent were given by Milnor
\cite{Mi}.  See also \cite{Stallings2}.  This is the ``failure of the
Hauptvermutung'' (in general).  These polyhedra are not PL manifolds,
and it turns out that there are examples of compact PL manifolds which
are homeomorphic but not piecewise-linearly equivalent too.  See
\cite{Sie} for dimensions $5$ and higher, and \cite{DoK, FQ} for
dimension $4$.  In dimensions $3$ and lower, this does not happen
\cite{Moise, Bing-book}.  The examples in \cite{Mi, Stallings2, Sie}
involved non-PL homeomorphisms whose behavior is much milder than in
the case of double-suspension spheres.  There are general results in
this direction for PL manifolds (and more broadly) in dimensions
greater than or equal to $5$.  See \cite{Su1, SS}.  Analogous
statements fail in dimension $4$, by \cite{DoS}.

	Some other examples where homeomorphic coordinates do not
exist, or necessarily have complicated behavior, even though the
geometry behaves well in other ways, are given in \cite{Seqs, Sebil}. 

	See \cite{DS-trans, HS, HY, MuS, Se-casscII, T1, T2} for some
related topics concerning homeomorphisms and bounds for their
behavior.

	One can try to avoid difficulties connected to $\pi_1$ by
using mappings with branching rather than homeomorphisms.  This is
discussed further in Appendix \ref{branching}.  

	Questions of algorithmic undecidability in topology have been
revisited in recent years, in particular by Nabutovsky and Weinberger.
See \cite{NW, NW2}, for instance, and the references therein.

\section{Quantitative topology, and calculus on singular spaces}
\label{Quantitative topology, and calculus on singular spaces}
\setcounter{equation}{0}

	One of the nice features of Euclidean spaces is that it is
easy to work with functions, derivatives, and integrals.  Here is
a basic example of this.  Let $f$ be a real-valued function on 
${\bf R}^n$ which is continuously differentiable and has compact
support, and fix a point $x \in {\bf R}^n$.  Then
\begin{equation}
\label{basic inequality}
	|f(x)| \le {1 \over \nu_n} 
		\int_{{\bf R}^n} {1 \over |x-y|^{n-1}} \, |\nabla f(y)| \, dy,
\end{equation}
where $\nu_n$ denotes the $(n-1)$-dimensional volume of the unit
sphere in ${\bf R}^n$, and $dy$ refers to ordinary $n$-dimensional
volume.

	This inequality provides a way to say that the values of a
function are controlled by \emph{averages} of its derivative.  In this
respect it is like Sobolev and isoperimetric inequalities, to which we
shall return in a moment.

	To prove (\ref{basic inequality}) one can proceed as follows
(as on p125 of \cite{St1}).  Let $v$ be any element of ${\bf R}^n$
with $|v|=1$.  Then
\begin{equation}
	f(x) = - \int_0^\infty {\partial \over \partial t} f(x + t v) \, dt,
\end{equation}
by the fundamental theorem of calculus.  Thus
\begin{equation}
\label{inequality in a line}
	|f(x)| \le \int_0^\infty |\nabla f(x + t v)| \, dt.
\end{equation}
This is true for every $v$ in the unit sphere of ${\bf R}^n$, and by
averaging over these $v$'s one can derive (\ref{basic inequality})
from (\ref{inequality in a line}). 

	To put this into perspective, it is helpful to look at a
situation where analogous inequalities make sense but fail to hold.
Imagine that one is interested in inequalities like (\ref{basic
inequality}), but for $2$-dimensional surfaces in ${\bf R}^3$ instead
of Euclidean spaces themselves.  Let $S$ be a smoothly embedded
$2$-dimensional submanifold of ${\bf R}^3$ which looks like a
$2$-plane with a bubble attached to it.  Specifically, let us start
with the union of a $2$-plane $P$ and a standard (round)
$2$-dimensional sphere $\Sigma$ which is tangent to $P$ at a single
point $z$.  Then cut out a little neighborhood of $z$, and glue in a
small ``neck'' as a bridge between the plane and the sphere to get a
smooth surface $S$.

	If the neck in $S$ is very small compared to the size of
$\Sigma$, then this is bad for an inequality like (\ref{basic
inequality}).  Indeed, let $x$ be the point on $\Sigma$ which is
exactly opposite from $z$, and consider a smooth function $f$ which is
equal to $1$ on most of $\Sigma$ (and at $x$ in particular) and equal
to $0$ on most of $P$.  More precisely, let us choose $f$ so that its
gradient is concentrated near the bridge between $\Sigma$ and $P$.  If
$f$ makes the transition from vanishing to being $1$ in a reasonable
manner, then the integral of $|\nabla f|$ on $S$ will be very small.
This is not hard to check, and it is bad for having an inequality like
(\ref{basic inequality}), since the left-hand side would be $1$ and
the right-hand side would be small.  In particular, one could not have
uniform bounds that would work for arbitrarily small bridges between
$P$ and $\Sigma$.

	The inequality (\ref{basic inequality}) is a relative of the
usual Sobolev and isoperimetric inequalities, which say the following.
Fix a dimension $n$ again, and an exponent $p$ that satisfies $1 \le p
< n$.  Define $q$ by $1/q = 1/p - 1/n$, so that $p < q < \infty$.
The Sobolev inequalities assert the existence of a constant $C(n,p)$
such that
\begin{equation}
\label{sobolev inequality}
	\bigg(\int_{{\bf R}^n} |f(x)|^q \, dx \bigg)^\frac{1}{q}
		\le C(n,p) \,
	\bigg(\int_{{\bf R}^n} |\nabla f(x)|^p \, dx \bigg)^\frac{1}{p}
\end{equation}
for all functions $f$ on ${\bf R}^n$ that are continuously
differentiable and have compact support.  One can also allow more
general functions, with $\nabla f$ interpreted in the sense of
distributions.

	  The \emph{isoperimetric inequality} states that if $D$ is a
domain in ${\bf R}^n$ (which is bounded and has reasonably smooth
boundary, say), then
\begin{eqnarray}
\label{isoperimetric inequality}
\lefteqn{
	n\hbox{-dimensional volume of } D} 			\\ 
		& & \qquad \quad \le C(n) \, 
	((n-1)\hbox{-dimensional volume of } \partial D)^\frac{n}{n-1}.
							\nonumber
\end{eqnarray}
This is really just a special case of (\ref{sobolev inequality}), with $p=1$
and $f$ taken to be the characteristic function of $D$ (i.e., the function
that is equal to $1$ on $D$ and $0$ on the complement of $D$).  In this
case $\nabla f$ is a (vector-valued) measure, and the right-hand side of
(\ref{sobolev inequality}) should be interpreted accordingly.  Conversely,
Sobolev inequalities for all $p$ can be derived from isoperimetric
inequalities, by applying the latter to sets of the form
\begin{equation}
	\{x \in {\bf R}^n : |f(x)| > t\},
\end{equation}
and then making suitable integrations in $t$.
	
	The sharp version of the isoperimetric inequality states that
(\ref{isoperimetric inequality}) holds with the constant $C(n)$ that
gives equality in the case of a ball.  See \cite{Fe}.  One can also
determine sharp constants for (\ref{sobolev inequality}), as on
p39 of \cite{Aubin}.

	Note that the choice of the exponent $n/(n-1)$ in the right
side of (\ref{isoperimetric inequality}) is determined by scaling
considerations, i.e., in looking what happens to the two sides of
(\ref{isoperimetric inequality}) when one dilates the domain $D$
by a positive factor.  The same is true of the relationship between
$p$ and $q$ in (\ref{sobolev inequality}), and the power $n-1$ in
the kernel on the right side of (\ref{basic inequality}).

	The inequality (\ref{basic inequality}) is a basic ingredient
in one of the standard methods for proving Sobolev and isoperimetric
inequalities (but not necessarily with sharp constants).  Roughly
speaking, once one has (\ref{basic inequality}), the rest of the
argument works at a very general level of integral operators on
measure spaces, rather than manifolds and derivatives.  This is not
quite true for the $p=1$ case of (\ref{sobolev inequality}), for which
the general measure-theoretic argument gives a slightly weaker result.
See Chapter V of \cite{St1} for details.  The slightly weaker result
does give an isoperimetric inequality (\ref{isoperimetric
inequality}), and it is not hard to recover the $p=1$ case of
(\ref{sobolev inequality}) from the weaker version using a bit more of
the localization properties of the gradient than are kept in
(\ref{basic inequality}).  (See also Appendix C of \cite{Setop},
especially Proposition C.14.)

	The idea of these inequalities makes sense much more broadly
than just on Euclidean spaces, but they may not always work very well,
as in the earlier example with bubbling.  To consider this further,
let $M$ be a smooth Riemannian manifold of dimension $n$, and let us
assume for simplicity that $M$ is unbounded (like ${\bf R}^n$).  Let
us also think of $M$ as coming equipped with a distance function
$d(x,y)$ with the usual properties ($d(x,y)$ is nonnegative, symmetric
in $x$ and $y$, vanishes exactly when $x=y$, and satisfies the
triangle inequality).  One might take $d(x,y)$ to be the geodesic
distance associated to the Riemannian metric on $M$, but let us not
restrict ourselves to this case.  For instance, imagine that $M$ is a
smooth submanifold of some higher-dimensional ${\bf R}^k$, and that
$d(x,y)$ is simply the ambient Euclidean distance $|x-y|$ inherited
from ${\bf R}^k$.  In general this could be much smaller than the
geodesic distance.

	We shall make the standing assumption that the distance
$d(x,y)$ and the Riemannian geodesic distance are approximately the
same, each bounded by twice the other, on sufficiently small
neighborhoods about any given point in $M$.  This ensures that
$d(x,y)$ is compatible with quantities defined locally on $M$ using
the Riemannian metric, like the volume measure, and the length of the
gradient of a function.  Note that this local compatibility condition
for the distance function $d(x,y)$ and the Riemannian metric is
satisfied automatically in the situation mentioned above, where $M$ is
a submanifold of a larger Euclidean space and $d(x,y)$ is inherited
from the ambient Euclidean distance.  We shall also require that the
distance $d(x,y)$ be compatible with the (manifold) topology on $M$,
and that it be complete.  This prevents things like infinite ends in
$M$ which asymptotically approach finite points in $M$ with respect to
$d(x,y)$.

	The smoothness of $M$ should be taken in the character of an a
priori assumption, with the real point being to have bounds that do
not depend on the presence of the smoothness in a quantitative way.
Indeed, the smoothness of $M$ will not really play an essential role
here, but will be convenient, so that concepts like volume, gradient,
and lengths of gradients are automatically meaningful.

	Suppose for the moment that $M$ is \emph{bilipschitz
equivalent} to ${\bf R}^n$ equipped with the usual Euclidean metric.
This means that there is a mapping $\phi$ from ${\bf R}^n$ onto $M$
and a constant $k$ such that
\begin{equation}
\label{k-bilipschitz}
	k^{-1} |z-w| \le d(\phi(z), \phi(w)) \le k \, |z-w|
			\qquad\hbox{for all } z, w \in {\bf R}^n.
\end{equation}
In other words, $\phi$ should neither expand or shrink distances by
more than a bounded amount.  This implies that $\phi$ does not distort
the corresponding Riemannian metrics or volume by more than bounded
factors either, as one can readily show.  In this case the analogues
of (\ref{basic inequality}), (\ref{sobolev inequality}), and
(\ref{isoperimetric inequality}) all hold for $M$, with constants that
depend only on the constants for ${\bf R}^n$ and the distortion factor
$k$.  This is because any test of these inequalities on $M$ can be
``pulled back'' to ${\bf R}^n$ using $\phi$, with the loss of
information in moving between $M$ and ${\bf R}^n$ limited by the
bilipschitz condition for $\phi$.

	This observation helps to make clear the fact that
inequalities like (\ref{basic inequality}), (\ref{sobolev
inequality}), and (\ref{isoperimetric inequality}) do not really
require much in the way of smoothness for the underlying space.
Bounds on curvature are not preserved by bilipschitz mappings, just as
bounds on higher derivative of functions are not preserved.
Bilipschitz mappings can allow plenty of spiralling and corners in $M$
(or approximate corners, since we are asking that $M$ that be smooth a
priori).

	Although bilipschitz mappings are appropriate here for the
small amount of regularity involved, the idea of a
``parameterization'' is too strong for the purposes of inequalities
like (\ref{basic inequality}), (\ref{sobolev inequality}), and
(\ref{isoperimetric inequality}).  One might say that these
inequalities are like algebraic topology, but more quantitative, while
parameterizations are more like homeomorphisms, which are always more
difficult.  (Some other themes along these lines will be discussed in
Appendix \ref{Working on spaces which may not have nice coordinates}.
Appendix \ref{more on existence and behavior of homeomorphisms} is
related to this as well.  See also \cite{hanson-juha}.)

	I would like to describe now some conditions on $M$ which are
strong enough to give bounds as in (\ref{basic inequality}), but which
are quite a bit weaker than the existence of a bilipschitz
parameterization.  First, let us explicitly write down the analogue of
(\ref{basic inequality}) for $M$.  If $x$ is any element of $M$, this
analogue would say that there is a constant $C$ so that
\begin{equation}
\label{basic inequality, general version}
	|f(x)| \le C \, 
		\int_M {1 \over d(x,y)^{n-1}} \, |\nabla f(y)| \, dVol(y)
\end{equation}
for all continuously differentiable functions $f$ on $M$, where
$|\nabla f(y)|$ and the volume measure $dVol(y)$ are defined in terms
of the Riemannian structure that comes with $M$.

	The next two definitions give the conditions on $M$ that we
shall consider.  These and similar notions have come up many times in
various parts of geometry and analysis, as in \cite{Al, AV1, AV2, A1,
A2, A3, CW1, CW2, Gromov-green, Gromov-translation, HK1, HK2, HK3, HY,
Pe1, Pe2, V3}.

\begindefinition [The doubling condition]
\label{doubling con}
A metric space $(M, d(x,y))$ is said to be \emph{doubling} (with constant
$L_0$) if each ball $B$ in $M$ with respect to $d(x,y)$ can be covered by at
most $L_0$ balls of half the radius of $B$.
\end{definition}

	Notice that Euclidean spaces are automatically doubling, with
a constant $L_0$ that depends only on the dimension.  Similarly, every
subset of a (finite-dimensional) Euclidean space is doubling, with a
uniform bound for its doubling constant.

\begindefinition [Local linear contractability] 
\label{local lin con}
A metric space $(M,d(x,y))$ is said to be \emph{locally linearly
contractable} (with constant $L_1$) if the following is true.  Let $B$
be a ball in $M$ with respect to $d(x,y)$, and with radius no greater
than $L_1^{-1}$ times the diameter of $M$.  (Arbitrary radii are
permitted when $M$ is unbounded, as in the context of the present
general discussion.)  Then (local linear contractability means that)
it should be possible to continuously contract $B$ to a point inside
of $L_1 \, B$, i.e., inside the ball with the same center as $B$ and
$L_1$ times the radius.
\end{definition}

	This is a kind of quantitative and scale-invariant condition
of local contractability.  It prevents certain types of cusps or
bubbling, for instance.  Both this and the doubling condition hold
automatically when $M$ admits a bilipschitz parameterization by ${\bf
R}^n$, with uniform bounds in terms of the bilipschitz constant $k$ in
(\ref{k-bilipschitz}) (and the dimension for the doubling condition).

\begintheorem 
\label{theorem}
If $M$ and $d(x,y)$ are as before, and if $(M,d(x,y))$
satisfies the doubling and local linear contractability conditions with
constants $L_0$ and $L_1$, respectively, then (\ref{basic inequality,
general version}) holds with a constant $C$ that depends only on
$L_0$, $L_1$, and the dimension $n$.
\end{theorem}

	This was proved in \cite{Setop}.  Before we look at some
aspects of the proof, some remarks are in order about what the
conclusions really mean.

	In general one cannot derive bounds for Sobolev and
isoperimetric inequalities for $M$ just using (\ref{basic inequality,
general version}).  One might say that (\ref{basic inequality, general
version}) is only as good as the behavior of the volume measure on
$M$.  If the volume measure on $M$ behaves well, with bounds for the
measure of balls like ones on ${\bf R}^n$, then one can derive
conclusions from (\ref{basic inequality, general version}) in much the
same way as for Euclidean spaces.  See Appendices B and C in
\cite{Setop}.

	The doubling and local linear contractability conditions do
not themselves say anything about the behavior of the volume on $M$,
and indeed they tolerate fractal behavior perfectly well.  To see
this, consider the metric space which is ${\bf R}^n$ as a set, but
with the metric $|x-y|^\alpha$, where $\alpha$ is some fixed number in
$(0,1)$.  This is a kind of abstract and higher-dimensional version of
standard fractal snowflake curves in the plane.  However, the doubling
and local linear contractability conditions work just as well for
$({\bf R}^n, |x-y|^\alpha)$ as for $({\bf R}^n, |x-y|)$, just with
slightly different constants.

	How might one prove Theorem \ref{theorem}?  It would be nice
to be able to mimic the proof of (\ref{basic inequality}), i.e., to
find a family of rectifiable curves in $M$ which go from $x$ to
infinity and whose arclength measures have approximately the same kind
of distribution in $M$ as rays in ${\bf R}^n$ emanating from a given
point.  Such families exist (with suitable bounds) when $M$ admits a
bilipschitz parameterization by ${\bf R}^n$, and they also exist in
more singular circumstances.

	Unfortunately, it is not so clear how to produce families of
curves like these without some explicit information about the space
$M$ in question.  This problem was treated in a special case in
\cite{DS-sob}, with $M$ a certain kind of (nonsmooth) conformal
deformation of ${\bf R}^n$.  The basic idea was to obtain these curves
from level sets of certain mappings with controlled behavior.  When
$n=2$, for instance, imagine a standard square $Q$, with opposing
vertices $p$ and $q$.  The boundary of $Q$ can be thought of as a pair
of paths $\alpha$, $\beta$ from $p$ to $q$, each with two segments, two
sides of $Q$.  If $\tau$ is a function on $Q$ which equals $0$ on
$\alpha$ and $1$ on $\beta$ (and is somewhat singular at $p$ and $q$),
then one can try to extract a family of paths from $p$ to $q$ in $Q$
from the level sets 
\begin{equation}
	\{x \in {\bf R}^2 : \tau(x) = t\}, \qquad 0 < t < 1.
\end{equation}
For the standard geometry on ${\bf R}^2$ one can write down a good
family of curves and a good function $\tau$ explicitly.  For a certain
class of conformal deformations of ${\bf R}^2$ one can make
constructions of functions $\tau$ with approximately the same behavior
as in the case of the standard metric, and from these one can get
controlled families of curves.

	These constructions of functions $\tau$ used the standard
Euclidean geometry in the background in an important way.  For the
more general setting of Theorem \ref{theorem} one needs to proceed
somewhat differently, and it is helpful to begin with a different 
formulation of the kind of auxiliary functions to be used.

	Given a point $x$ in ${\bf R}^n$, there is an associated
spherical projection $\pi_x : {\bf R}^n \backslash \{x\} \to {\bf
S}^{n-1}$ given by
\begin{equation}
\label{standard pi_x}
	\pi_x(u) = {u-x \over |u-x|}.
\end{equation}
This projection is topologically nondegenerate, in the sense that it
has degree equal to $1$.  Here the ``degree'' can be defined by
restricting $\pi_x$ to a sphere around $x$ and taking the degree of
this mapping (from an $(n-1)$-dimensional sphere to another one) in
the usual sense.  (See \cite{Massey, Milnor, Nirenberg} concerning the
notion of degree of a mapping.)  Also, this mapping satisfies the
bound
\begin{equation}
\label{bound for the differential}
	|d \pi_x(u)| \le C \, |u-x|^{-1}
\end{equation}
for all $u \in {\bf R}^n \backslash \{x\}$, where $d \pi_x(u)$ denotes
the differential of $\pi_x$ at $u$, and $C$ is some constant.  One can
write down the differential of $\pi_x$ explicitly, and (\ref{bound for
the differential}) can be replaced by an equality, but this precision
is not needed here, and not available in general.  The rays in ${\bf
R}^n$ that emanate from $x$ are exactly the fibers of the mapping
$\pi_x$, and bounds for the distribution of their arclength measures
can be seen as a consequence of (\ref{bound for the differential}),
using the ``co-area theorem'' \cite{Fe, Morgan, Simon}.

	One can also think of $\pi_x$ as giving (\ref{basic
inequality}) in the following manner.  Let $\omega$ denote the
standard volume form on ${\bf S}^{n-1}$, a differential form
of degree $n-1$, and normalized so that
\begin{equation}
	\int_{{\bf S}^{n-1}} \omega = 1.
\end{equation}
Let $\lambda$ denote the differential form on ${\bf R}^n \backslash
\{0\}$ obtained by pulling $\omega$ back using $\pi_x$.  Then
(\ref{bound for the differential}) yields
\begin{equation}
\label{bound for lambda}
	|\lambda(u)| \le C' \, |u-x|^{-n+1}
\end{equation}
for all $u \in {\bf R}^n \backslash \{x\}$, where $C'$ is a slightly
different constant from before.  In particular, $\lambda$ is locally
integrable across $x$ (and smooth everywhere else).  This permits one
to take the exterior derivative of $\lambda$ on all of ${\bf R}^n$ in
the (distributional) sense of currents \cite{Fe, Morgan}, and the
result is that $d \lambda$ is the current of degree $n$ which is a
Dirac mass at $x$.  More precisely, $d \lambda = 0$ away from $x$
because $\omega$ is automatically closed (being a form of top degree
on ${\bf S}^{n-1}$), and because the pull-back of a closed form is
always closed.  The Dirac mass at $x$ comes from a standard Stokes'
theorem computation, which uses the observation that the integral of
$\lambda$ over any $(n-1)$-sphere in ${\bf R}^n$ around $x$ is equal
to $1$.  (The latter is one way to formulate the fact that the degree
of $\pi_x$ is $1$.)

	This characterization of $d \lambda$ as a current on ${\bf
R}^n$ means that
\begin{equation}
	\int_{{\bf R}^n} df \wedge \lambda = - f(x)
\end{equation}
when $f$ is a smooth function on ${\bf R}^n$ with compact support.
This yields (\ref{basic inequality}), because of (\ref{bound for
lambda}).  (A similar use of differential forms was employed in
\cite{DS-sob}.)

	The general idea of the mapping $\pi_x$ also makes sense in
the context of Theorem \ref{theorem}.  Let $M$, $d(y,z)$ be as before,
and fix a point $x$ in $M$.  One would like to find a mapping $\pi_x :
M \backslash \{x\} \to {\bf S}^{n-1}$ which is topologically
nondegenerate and satisfies
\begin{equation}
\label{bound for the differential on M}
	|d \pi_x(u)| \le K \, d(u,x)^{-1}
\end{equation}
for some constant $K$ and all $u \in M \backslash \{x\}$.  Note that
now the norm of the differential of $\pi_x$ involves the Riemannian
metric on $M$.  For the topological nondegeneracy of $\pi_x$, let us
ask that it have nonzero degree on small spheres in $M$ that surround
$x$ in a standard way.  This makes sense, because of the a priori
assumption that $M$ be smooth.

	If one can produce such a mapping $\pi_x$, then one can derive
(\ref{basic inequality, general version}) as a consequence, using the
same kind of argument with differential forms as above.  One can also
find enough curves in the fibers of $\pi_x$, with control on the way
that their arclength measures are distributed in $M$, through the use
of co-area estimates.  For this the topological nondegeneracy of
$\pi_x$ is needed for showing that the fibers of $\pi_x$ connect
$x$ to infinity in $M$.

	In the context of conformal deformations of ${\bf R}^n$, as in
\cite{DS-sob}, such mappings $\pi_x$ can be obtained as perturbations
to the standard mapping in (\ref{standard pi_x}).  This is described
in \cite{Sebarc}.  For Theorem \ref{theorem}, the method of
\cite{Setop} does not use mappings quite like $\pi_x$, but a
``stabilized'' version from which one can draw similar conclusions.
In this stabilized version one looks for mappings from $M$ to ${\bf
S}^n$ (instead of ${\bf S}^{n-1}$) which are constant outside of a
(given) ball, topologically nontrivial (in the sense of nonzero
degree), and which satisfy suitable bounds on their differentials.
These mappings are like snapshots of pieces of $M$, and one has to
move them around in a controlled manner.  This means moving them both
in terms of location (the center of the supporting ball) and scale
(the radius of the ball).

	At this stage the hypotheses of Theorem \ref{theorem} may make
more sense.  Existence of mappings like the ones described above is a
standard matter in topology, except for the question of uniform
bounds.  The hypotheses of Theorem \ref{theorem} (the doubling
condition and local linear contractability) are also in the nature of
quantitative topology.  Note, however, that the kind of bounds
involved in the hypotheses of the theorem and the construction of
mappings into spheres are somewhat different from each other, with
bounds on the differentials being crucial for the latter, while
control over moduli of continuity does not come up in the former.
(The local linear contractability condition restricts the overall
distances by which points are displaced in the contractions, but not
the sizes of the smaller-scale oscillations, as in a modulus of
continuity.)  In the end the bounds for the differentials come about
because the hypotheses of Theorem \ref{theorem} permit one to reduce
various constructions and comparisons to finite models of controlled
complexity.

	In the proof of Theorem \ref{theorem} there are three related
pieces of information that come out, namely (1) estimates for the
behavior of functions on our space $M$ in terms of their derivatives,
as in (\ref{basic inequality, general version}), (2) families of
curves in $M$ which are well-distributed in terms of arclength
measure, and (3) mappings to spheres with certain estimates and
nondegeneracy properties.  These three kinds of information are
closely linked, through various dualities, but to some extent they
also have their own lives.  Each would be immediate if $M$ had a
bilipschitz parameterization by ${\bf R}^n$, but in fact they are more
robust than that, and much easier to verify.

	Indeed, one of the original motivations for \cite{DS-sob} was
the problem of determining which conformal deformations of ${\bf R}^n$
lead to metric spaces (through the geodesic distance) which are
bilipschitz equivalent to ${\bf R}^n$.  The deformations are allowed
to be nonsmooth here, but this does not matter too much, because of
the natural scale-invariance of the problem, and because one seeks
uniform bounds.  This problem is the same in essence as asking which
(positive) functions on ${\bf R}^n$ arise as the Jacobian of a
quasiconformal mapping, modulo multiplication by a positive function
which is bounded and bounded away from $0$.

	Some natural necessary conditions are known for these
questions, with a principal ingredient coming from \cite{Ge}.  It was
natural to wonder whether the necessary conditions were also
sufficient.  As a test for this, \cite{DS-sob} looked at the Sobolev
and related inequalities that would follow if the necessary conditions
were sufficient.  These inequalities could be stated directly in terms
of the data of the problems, the conformal factor or prospective
Jacobian.  The conclusion of \cite{DS-sob} was that these inequalities
could be derived directly from the conditions on the data,
independently of whether these conditions were sufficient for the
existence of bilipschitz/quasiconformal mappings as above.

	In \cite{Sebil} it was shown that the candidate conditions are
\emph{not} sufficient for the existence of such mappings, at least in
dimensions $3$ and higher.  (Dimension $2$ remains open.)  The
simplest counterexamples involved considerations of localized
fundamental groups, in much the same fashion as in Section
\ref{parameterization problems}.  (Another class of counterexamples
were based on a different mechanism, although these did not start in
dimension $3$.)  These counterexamples are all perfectly well-behaved
in terms of the doubling and local linear contractability properties,
and in fact are much better than that.

	Part of the bottom line here is that spaces can have geometry
which behaves quite well for many purposes even if they do not behave
so well in terms of parameterizations.

	For some other aspects of ``quantitative topology'', see
\cite{Al, AV1, AV2, At1, At2, BW, CF, Ch, Fe1, Fe2, Fe3, Fe4, Ge-link,
Gromov-green, Gromov-translation, HY, HS, Lu, Pe1, Pe2, TV, V1, V2,
V3}.  Related matters of Sobolev and other inequalities on non-smooth
spaces come up in \cite{HK2, HK3, HKST}, in connection with the
behavior of quasiconformal mappings.

\section{Uniform rectifiability}
\label{Uniform rectifiability}
\setcounter{equation}{0}

	A basic fact in topology is that there are spaces which are
\emph{manifold factors} but not manifolds.  That is, there are
topological spaces $M$ such that $M \times {\bf R}$ is a manifold
(locally homeomorphic to a Euclidean space) but $M$ is not.  This can
even happen for finite polyhedra, because of the double-suspension
results of Edwards and Cannon.  See \cite{Dm2, E, K} for more
information and specific examples.  (We shall say a bit more about
this in Appendix \ref{more on existence and behavior of
homeomorphisms}.)

	\emph{Uniform rectifiability} is a notion of controlled
geometry that trades topology for estimates.  It tolerates some amount
of singularities, like holes and crossings, and avoids some common
difficulties with homeomorphisms, such as manifold factors.

	The precise definition is slightly technical, and relies on
measure theory in a crucial way.  In many respects it is analogous to
the notion of BMO from Section \ref{BMO frame of mind}.  The following
is a preliminary concept that helps to set the stage.

\begindefinition [Ahlfors regularity]
\label{ahlfors regular}
Fix $n$ and $d$, with $n$ a positive integer and $0 < d \le n$.  A set
$E$ contained in ${\bf R}^n$ is said to be \emph{(Ahlfors) regular of
dimension $d$} if it is closed, and if there is a positive Borel
measure $\mu$ supported on $E$ and a constant $C > 0$ such that
\begin{equation}
	C^{-1} \, r^d \le \mu(B(x,r)) \le C \, r^d
\end{equation}
for all $x \in E$ and $0 < r \le \diam E$.  Here $B(x,r)$ denotes
the (open) ball with center $x$ and radius $r$.
\end{definition}

	Roughly speaking, this definition asks that $E$ behave like
ordinary Euclidean space in terms of the distribution of its mass.
Notice that $d$-planes satisfy this condition automatically, with
$\mu$ equal to the ordinary $d$-dimensional volume.  The same is true
for compact smooth manifolds, and finite polyhedra which are given as
unions of $d$-dimensional simplices (i.e., with no lower-dimensional
pieces sticking off in an isolated manner).  There are also plenty of
``fractal'' examples, like self-similar Cantor sets and snowflake
curves.  In particular, the dimension $d$ can be any (positive) real
number.

	A basic fact is that if $E$ is regular and $\mu$ is as in
Definition \ref{ahlfors regular}, then $\mu$ is practically the same
as $d$-dimensional Hausdorff measure $H^d$ restricted to $E$.
Specifically, $\mu$ and $H^d$ are each bounded by constant multiples
of the other when applied to subsets of $E$.  This is not hard to
prove, and it shows that $\mu$ is essentially unique.  Definition
\ref{ahlfors regular} could have been formulated directly in terms of
Hausdorff measure, but the version above is a bit more elementary.

	Let us recall the definition of a bilipschitz mapping.  Let
$A$ be a set in ${\bf R}^n$, and let $f$ be a mapping from $A$ to some
other set in ${\bf R}^n$.  We say that $f$ is \emph{$k$-bilipschitz},
where $k$ is a positive number, if
\begin{equation}
	k^{-1} \, |x-y| \le |f(x) - f(y)| \le k \, |x-y|
\end{equation}
for all $x, y \in A$.

\begindefinition [Uniform rectifiability]
\label{unif rec}
Let $E$ be a subset of ${\bf R}^n$ which is Ahlfors regular of
dimension $d$, where $d$ is a positive integer, $d < n$, and let $\mu$
be a positive measure on $E$ as in Definition \ref{ahlfors regular}.
Then $E$ is \emph{uniformly rectifiable} if there exists a positive
constant $k$ so that for each $x \in E$ and each $r > 0$ with $r \le
\diam E$ there is a closed subset $A$ of $E \cap \overline{B}(x,r)$
such that
\begin{equation}
\label{lower bound for mu(A)}
	\mu(A) \ge  \frac{9}{10} \cdot 
				\mu(E \cap \overline{B}(x,r))
\end{equation}
and
\begin{equation}
\label{there is a bilip map}
	\hbox{there is a $k$-bilipschitz mapping } f			
	 	\hbox{ from } A \hbox{ into } {\bf R}^d.	
\end{equation}
\end{definition}

	In other words, inside of each ``snapshot'' $E \cap
\overline{B}(x,r)$ of $E$ there should be a large subset, with at
least $90 \% $ of the points, which is bilipschitz equivalent to a
subset of ${\bf R}^d$, and with a uniform bound on the bilipschitz
constant.  This is like asking for a controlled parameterization,
except that we allow for holes and singularities.

	Definition \ref{unif rec} should be compared with the
classical notion of (countable) rectifiability, in which one asks that
$E$ be covered, except for a set of measure $0$, by a countable union
of sets, each of which is bilipschitz equivalent to a subset of ${\bf
R}^d$.  Uniform rectifiability implies this condition, but it is
stronger, because it provides quantitative information at definite
scales, while the classical notion really only gives asymptotic
information as one zooms in at almost any point.  See \cite{Fa, Fe,
Mat} for more information about classical rectifiability.

	Normally one would be much happier to simply have bilipschitz
coordinates outright, without having to allow for bad sets of small
measure where this does not work.  In practice bilipschitz coordinates
simply do not exist in many situations where one might otherwise hope
to have them.  This is illustrated by the double-suspension spheres of
Edwards and Cannon \cite{C1, C2, Dm2, E}, and the observations about
them in \cite{SS}.  Further examples are given in \cite{Seqs, Sebil}.

	The use of arbitrary scales and locations is an important part
of the story here, and is very similar to the concept of BMO.  At the
level of a single snapshot, a fixed ball $\overline{B}(x,r)$ centered
on $E$, the bad set may seem pretty wild, as nothing is said about
what goes on there in (\ref{lower bound for mu(A)}) or (\ref{there is
a bilip map}).  However, uniform rectifiability, like BMO, applies to
all snapshots equally, and in particular to balls in which the bad set
is concentrated.  Thus, inside the bad set, there are in fact further
controls.  We shall see other manifestations of this later, and the
same basic principle is used in the John--Nirenberg theorem for BMO
functions (discussed in Section \ref{BMO frame of mind}).

	Uniform rectifiability provides a substitute for (complete)
bilipschitz coordinates in much the same way that BMO provides a
substitute for $L^\infty$ bounds, as in Section \ref{BMO frame of
mind}.  Note that $L^\infty$ bounds and bilipschitz coordinates
automatically entail uniform control over all scales and locations.
This is true just because of the way they are defined, i.e., a bounded
function is bounded in all snapshots, and with a uniform majorant.
With BMO and uniform rectifiability the scale-invariance is imposed by
hand.

	It may be a little surprising that one can get anything new
through concepts like BMO and uniform rectifiability.  For instance,
suppose that $f$ is a locally-integrable function on ${\bf R}^k$, and
that the averages
\begin{equation}
	\frac{1}{\omega_k \, t^k} \int_{B(z,t)} |f(w)| \, dw
\end{equation}
are uniformly bounded, independently of $z$ and $t$.  Here $\omega_k$
denotes the volume of the unit ball in ${\bf R}^k$, so that $\omega_k
\, t^k$ is the volume of $B(z,t)$.  This implies that $f$ must itself
be bounded by the same amount almost everywhere on ${\bf R}^k$, since
\begin{equation}
	f(u) = 
	\lim_{t \to 0} \frac{1}{\omega_k \, t^k} \int_{B(z,t)} f(w) \, dw
\end{equation}
almost everywhere on ${\bf R}^k$.  Thus a uniform bound for the size
of the snapshots does imply a uniform bound outright.  For BMO the
situation is different because one asks only for a uniform bound on
the mean \emph{oscillation} in every ball.  In other words, one also
has the freedom to make renormalizations by additive constants when
moving from place to place, and this gives enough room for some
unbounded functions, like $\log |x|$.  Uniform rectifiability is like
this as well, although with different kinds of ``renormalizations''
available.

	These remarks might explain why \emph{some} condition like
uniform rectifiability could be useful or natural, but why the
specific version above in particular?  Part of the answer to this is
that nearly all definitions of this nature are equivalent to the
formulation given above.  For instance, the $9/10$ in (\ref{lower
bound for mu(A)}) can be replaced by any number strictly between $0$
and $1$.  See \cite{DS3, DS5} for more information.

	Another answer lies in a theme often articulated by Coifman,
about the way that operator theory can provide a good guide for
geometry.  One of the original motivations for uniform rectifiability
came from the ``Calder\'on program'' \cite{Ca2}, concerning the
$L^p$-boundedness of certain singular operators on curves and surfaces
of minimal smoothness.  David \cite{Da1, Da2, Da4} showed that uniform
rectifiability of a set $E$ implies $L^p$-boundedness of wide classes
of singular operators on $E$.  (See \cite{Ca1, Ca2, CDM, CMM} and the
references therein for related work connected to the Calder\'on
program.)  In \cite{DS3}, a converse was established, so that uniform
rectifiability of an Ahlfors-regular set $E$ is actually equivalent to
the boundedness of a suitable class of singular integral operators
(inherited from the ambient Euclidean space ${\bf R}^n$).  See also
\cite{DS-pubmat, DS5, MMV, MP}.

	Here is a concrete statement about uniform rectifiability in
situations where well-behaved parameterizations would be natural but
may not exist.

\begintheorem 
\label{unif rec theorem}
Let $E$ be a subset of ${\bf R}^n$ which is regular of dimension $d$.
If $E$ is also a $d$-dimensional topological manifold and satisfies
the local linear contractability condition (Definition \ref{local lin
con}), then $E$ is uniformly rectifiable.
\end{theorem}

	Note that Ahlfors-regularity automatically implies the doubling
condition (Definition \ref{doubling con}).
	
	Theorem \ref{unif rec theorem} has been proved by G. David and
myself.  Now-a-days we have better technology, which allows for
versions of this which are localized to individual ``snapshots'',
rather than using all scales and locations at once.  See
\cite{DS-q-min} (with some of the remarks in Section 12.3 of
\cite{DS-q-min} helping to provide a bridge to the present
formulation).  We shall say a bit more about this, near the end of
Subsection \ref{A class of variational problems}.

	The requirement that $E$ be a topological manifold is
convenient, but weaker conditions could be used.  For that matter,
there are natural variations of local linear contractability too.

	One can think of Theorem \ref{unif rec theorem} and related
results in the following terms.  Given a compact set $K$, upper bounds
for the $d$-dimensional Hausdorff measure of $K$ together with lower
bounds for the $d$-dimensional topology of $K$ should lead to strong
information about the geometric behavior of $K$.  See \cite{DS8,
DS-q-min, Secon} for more on this.

	To understand better what Theorem \ref{unif rec theorem}
means, let us begin by observing that the hypotheses of Theorem
\ref{unif rec theorem} would hold automatically if $E$ were
bilipschitz equivalent to ${\bf R}^d$, or if $E$ were compact and
admitted bilipschitz local coordinates from ${\bf R}^d$.  Under these
conditions, a test of the hypotheses of Theorem \ref{unif rec theorem}
on $E$ can be converted into a similar test on ${\bf R}^d$, where it
can then be resolved in a straightforward manner.

	A similar argument shows that the hypotheses of Theorem
\ref{unif rec theorem} are ``bilipschitz invariant''.  More
precisely, if $F$ is another subset of ${\bf R}^n$ which is
bilipschitz equivalent to $E$, and if the hypotheses of Theorem
\ref{unif rec theorem} holds for one of $E$ and $F$, then it
automatically holds for the other.

	Since the existence of bilipschitz coordinates implies the
hypotheses of Theorem \ref{unif rec theorem}, we cannot ask for more
than that in the conclusions.  In other words, bilipschitz coordinates
are at the high end of what one can hope for in the context of Theorem
\ref{unif rec theorem}.  The hypotheses of Theorem \ref{unif rec
theorem} do in fact rule out a lot of basic obstructions to the
existence of bilipschitz coordinates, like cusps, fractal behavior,
self-intersections and approximate self-intersections, and bubbles
with very small necks.  (Compare with Section \ref{Quantitative
topology, and calculus on singular spaces}, especially Theorem
\ref{theorem} and the discussion of its proof and consequences.)
Nonetheless, it can easily happen that a set $E$ satisfies the
hypotheses of Theorem \ref{unif rec theorem} but does not admit
bilipschitz local coordinates.  Double-suspension spheres provide
spectacular counterexamples for this (using the observations of
\cite{SS}).  Additional counterexamples are given in \cite{Seqs,
Sebil}.

	We should perhaps emphasize that the assumption of being a
topological manifold in Theorem \ref{unif rec theorem} does not
involve bounds.  By contrast, uniform rectifiability does involve
bounds, which is part of the point.  In the context of Theorem
\ref{unif rec theorem}, the proof shows that the uniform
rectifiability constants for the conclusion are controlled in terms of
the constants that are implicit in the hypotheses, i.e., in
Ahlfors-regularity, the linear contractability condition, and the
dimension.

	If bilipschitz coordinates are at the high end of what one
could hope for, what happens if one asks for less?  What if one asks
for homeomorphic local coordinates with some control, but not as much?
For instance, instead of bounding the ``rate'' of continuity through
Lipschitz conditions like
\begin{equation}
\label{Lipschitz}
	|f(x) - f(y)| \le C \, |x-y|
\end{equation}
(for some $C$ and all $x$, $y$ in the domain of $f$), one could work
with H\"older continuity conditions, which have the form
\begin{equation}
\label{holder cont.}
	|f(x) - f(y)| \le C' \, |x-y|^\gamma.
\end{equation}
Here $\gamma$ is a positive number, sometimes called the H\"older
``exponent''.  As usual, (\ref{holder cont.}) is supposed to hold
simultaneously for all $x$ and $y$ in the domain of $f$, and with a
fixed constant $C'$.  When $x$ and $y$ are close to each other and
$\gamma$ is less than $1$, this type of condition is strictly weaker
than that of being Lipschitz.  Just as $f(x) = |x|$ is a standard
example of a Lipschitz function that is not differentiable at the
origin, $g(x) = |x|^\gamma$ is a basic example of a function that is
H\"older continuous of order $\gamma$, $\gamma \le 1$, but not of any
order larger than $\gamma$, in any neighborhood of the origin.

	Instead of local coordinates which are bilipschitz, one could
consider ones that are ``bi-H\"older'', i.e., H\"older continuous and
with H\"older continuous inverse.  It turns out that double-suspension
spheres do not admit bi-H\"older local coordinates when the H\"older
exponent $\gamma$ lies above an explicit threshold.  Specifically, if
$P$ is an $n$-dimensional polyhedron which is the double-suspension of
an $(n-2)$-dimensional homology sphere that is not simply connected,
then there are points in $P$ (along the ``suspension circle'') for
which bi-H\"older local coordinates of exponent $\gamma > 1/(n-2)$ do
not exist.  This comes from the same argument as in \cite{SS}.  More
precisely, around these points in $P$, there do not exist homeomorphic
local coordinates from subsets of ${\bf R}^n$ for which the inverse
mapping is H\"older continuous of order $\gamma > 1/(n-2)$ (without
requiring a H\"older condition for the mapping itself).  

	Given any positive number $a$, there are examples in
\cite{Sebil} so that local coordinates (at some points) cannot have
their inverses be H\"older continuous of order $a$.  These examples do
admit bi-H\"older local coordinates (with a smaller exponent), and
even ``quasisymmetric'' \cite{TV} coordinates, and they satisfy the
hypotheses of Theorem \ref{unif rec theorem}.  In \cite{Seqs} there
are examples which satisfy the hypotheses of Theorem \ref{unif rec
theorem}, but for which no uniform modulus of continuity for local
coordinate mappings and their inverses is possible (over all scales
and locations).

	Related topics will be discussed in Appendix \ref{more on
existence and behavior of homeomorphisms}.

\subsection{Smoothness of Lipschitz and bilipschitz mappings}
\label{first subsection}

	Another aspect of uniform rectifiability is that it provides
the same amount of ``smoothness'' as when there is a global
bilipschitz parameterization.  To make this precise, let us first look
at the smoothness of Lipschitz and bilipschitz mappings.

	A mapping $f : {\bf R}^d \to {\bf R}^n$ is Lipschitz if there
is a constant $C$ so that (\ref{Lipschitz}) holds for all $x$, $y$ in
${\bf R}^d$.  The space of Lipschitz mappings is a bit simpler than
the space of bilipschitz mappings, because the former is a vector
space (and even a Banach space) while the latter is not.  For the
purposes of ``smoothness'' properties, though, there is not really any
difference between the two.  Bilipschitz mappings are always
Lipschitz, and anything that can happen with Lipschitz mappings can
also happen with bilipschitz mappings (by adding new components, or
considering $x + h(x)$ when $h(x)$ has Lipschitz norm less than $1$ to
get a bilipschitz mapping).

	One should also not worry too much about the difference
between Lipschitz mappings which are defined on all of ${\bf R}^d$,
and ones that are only defined on a subset.  Lipschitz mappings into
${\bf R}^n$ that are defined on a subset of ${\bf R}^d$ can always be
extended to Lipschitz mappings on all of ${\bf R}^d$.  This is a
standard fact.  There are also extension results for bilipschitz
mappings, if one permits oneself to replace the image ${\bf R}^n$ with
a Euclidean space of larger dimension (which is not too serious in the
present context).

	For considerations of ``smoothness'' we might as well restrict
our attention to functions which are real-valued, since the ${\bf
R}^n$-valued case can always be reduced to that.

	Two basic facts about Lipschitz functions on ${\bf R}^d$ are
that they are differentiable almost everywhere (with respect to
Lebesgue measure), and that for each $\eta > 0$ they can be modified
on sets of Lebesgue measure less than $\eta$ (depending on the
function) in such a way as to become \emph{continuously}
differentiable everywhere.  See \cite{Fe}.

	These are well-known results, but they do not tell the whole
story.  They are not quantitative; they say a lot about the asymptotic
behavior (on average) of Lipschitz mappings at very small scales, but
they do not say anything about what happens at scales of definite
size.

	To make this precise, let a Lipschitz mapping $f : {\bf R}^d
\to {\bf R}$ be given, and fix a point $x \in {\bf R}^d$ and a radius
$t > 0$.  We want to measure how well $f$ is approximated by an affine
function on the ball $B(x,t)$.  To do this we define the quantity
$\alpha(x,t)$ by
\begin{equation}
\label{def of alpha(x,t)}
	\alpha(x,t) = \inf_{A \in \mathcal{A}} \, \sup_{y \in B(x,t)} 
						t^{-1} \, |f(y) - A(y)|.
\end{equation}
Here $\mathcal{A}$ denotes the (vector space) of affine functions on
${\bf R}^d$.  The part on the right side of (\ref{def of alpha(x,t)})
with just the supremum (and not the infimum) measures how well the
particular affine function $A$ approximates $f$ inside $B(x,t)$, and
then the infimum gives us the best approximation by any affine
function for a particular choice of $x$ and $t$.  The factor of
$t^{-1}$ makes $\alpha(x,t)$ scale properly, and be dimensionless.  In
particular, $\alpha(x,t)$ is uniformly bounded in $x$ and $t$ when $f$
is Lipschitz, because we can take $A(y)$ to be the constant function
equal to the value of $f$ at $x$.

	The smallness of $\alpha(x,t)$ provides a manifestation of
the smoothness of $f$.  For functions which are twice-continuously
differentiable one can get estimates like
\begin{equation}
	\alpha(x,t) = O(t),
\end{equation}
using Taylor's theorem.  If $0 < \delta < 1$, then estimates like
\begin{equation}
\label{C^{1,delta}}
	\alpha(x,t) = O(t^\delta)
\end{equation}
(locally uniformly in $x$) correspond to H\"older continuity of the
gradient of $f$ of order $\delta$.  Differentiability almost
everywhere of $f$ implies that
\begin{equation}
\label{a.e. statement}
	\lim_{t \to 0} \alpha(x,t) = 0 \qquad \hbox{for almost every } x.
\end{equation}
This does not say anything about any particular $t$, because one does not
know how long one might have to wait before the limiting behavior kicks in.

	Here is a simple example.  Let us take $d=1$, and consider the function
\begin{equation}
\label{g_rho}
	g_\rho(x) = \rho \cdot \sin (x/\rho)
\end{equation}
on ${\bf R}$.  Here $\rho$ is any positive number.  Now, $g_\rho(x)$
is Lipschitz with norm $1$ no matter how $\rho$ is chosen.  This is
not hard to check; for instance, one can take the derivative to get
that
\begin{equation}
	g_\rho'(x) = \cos (x/\rho)
\end{equation}
so that $|g_\rho'(x)| \le 1$ everywhere.  This implies that
\begin{equation}
	|g_\rho(u) - g_\rho(v)| \le |u-v|
\end{equation}
for all $u$ and $v$ (and all $\rho$), because of the mean-value
theorem, or the fundamental theorem of calculus.  (One also has that
$|g_\rho'(x)|=1$ at some points, so that the Lipschitz norm is always
equal to $1$.)

	If $\rho$ is very small, then one has to wait a long time
before the limit in (\ref{a.e. statement}) takes full effect, because
$\alpha(x,t)$ will \emph{not} be small when $t = \rho$.  In fact,
there is then a positive lower bound for $\alpha(x,t)$ that does not
depend on $x$ or $\rho$ (assuming that $t$ is taken to be equal to
$\rho$).  This is not hard to verify directly.  One does not really
have to worry about $\rho$ here, because one can use scaling arguments
to reduce the lower bound to the case where $\rho = 1$.

	Thus a bound on the Lipschitz norm is not enough to say
anything about when the limit in (\ref{a.e. statement}) will take
effect.  These examples work uniformly in $x$, so that one cannot
avoid the problem by removing a set of small measure or anything like
that.

	However, there is something else that one can observe about
these examples.  Fix a $\rho$, no matter how large or small.  The
corresponding quantities $\alpha(x,t)$ will not be too small when $t$
is equal to $\rho$, as mentioned above, but they will be small when
$t$ is either much smaller than $\rho$, or much larger than $\rho$.
At scales much smaller than $\rho$, $g_\rho(x)$ is approximately
affine, because the smoothness of the sine function has a chance to
kick in, while at larger scales $g_\rho(x)$ is simply small outright
compared to $t$ (and one can take $A = 0$ as the approximating affine
function).

	In other words, for the functions $g_\rho(x)$ there is always
a bad scale where the $\alpha(x,t)$'s may not be small, and that bad
scale can be arbitrarily large or small, but the bad behavior is
confined to approximately just one scale.

	It turns out that something similar happens for arbitrary
Lipschitz functions.  The bad behavior cannot always be confined to a
single scale --- one might have sums of functions like the $g_\rho$'s,
but with very different choices of $\rho$ --- but, on average, the bad
behavior is limited to a bounded number of scales.

	Let us be more precise, and define a family of functions which
try to count the number of ``bad'' scales associated to a given point
$x$.  Fix a radius $r > 0$, and also a small number $\epsilon$, which
will provide our threshold for what is considered ``small''.  We
assume that a lipschitz function $f$ on ${\bf R}^d$ has been fixed, as
before.  Given $x \in {\bf R}^d$, define $N_r(x)$ to be the number
of nonnegative integers $j$ such that
\begin{equation}
\label{bad j's}
	\alpha(x, 2^{-j} \, r) \ge \epsilon.
\end{equation}
These $j$'s represent the ``bad'' scales for the point $x$, and below the
radius $r$.

	It is easy to see that (\ref{a.e. statement}) implies that
$N_r(x) < \infty$ for almost all $x$.  There is a more quantitative
statement which is true, namely that the average of $N_r(x)$ over any
ball $B$ in ${\bf R}^d$ of radius $r$ is finite and uniformly bounded,
independently of the ball $B$ and the choice of $r$.  That is,
\begin{equation}
\label{bound on average}
	r^{-d} \int_B N_r(x) \, dx \le C(n,\epsilon^{-1} \|f\|_{Lip}),
\end{equation}
where $C(n, s)$ is a constant that depends only on $n$ and
$s$, and $\|f\|_{Lip}$ denotes the Lipschitz norm of $f$.
This is a kind of ``Carleson measure condition''.

	Before we get to the reason for this bound, let us consider
some examples.  For the functions $g_\rho(x)$ in (\ref{g_rho}), the
functions $N_r(x)$ are simply uniformly bounded, independently of $x$,
$r$, and $\rho$.  This is not hard to check.  Notice that the bound
does depend on $\epsilon$, i.e., it blows up as $\epsilon \to 0$.  As
another example, consider the function $f$ defined by
\begin{equation}
	f(x) = |x|.
\end{equation}
For this function we have that $N_r(0) = \infty$ as soon as $\epsilon$
is small enough.  This is because $\alpha(0,t)$ is positive and
independent of $t$, so that (\ref{bad j's}) holds for all $j$ when
$\epsilon$ is sufficiently small.  Thus $N_r(x)$ is not uniformly
bounded in this case.  In fact it has a logarithmic singularity
near $0$, with $N_r(x)$ behaving roughly like $\log (r/|x|)$, and this
is compatible with (\ref{bound on average}) for $B$ centered at $0$.
If one is far enough away from the origin (compared to $r$), then
$N_r(x)$ simply vanishes, and there is nothing to do.

	In general one can have mixtures of the two types of phenomena.
Another interesting class of examples to consider are functions of the
form
\begin{equation}
	f(x) = \dist(x,F),
\end{equation}
where $F$ is some nonempty subset of ${\bf R}^d$ which is not all of
${\bf R}^d$, and $\dist(x,F)$ is defined (as usual) by
\begin{equation}
\label{def of dist}
	\dist(x,F) = \inf \{|x-z| : z \in F \}.
\end{equation}
It is a standard exercise that such a function $f$ is always Lipschitz
with norm at most $1$.  Depending on the behavior of the set $F$, this
function can have plenty of sharp corners, like $|x|$ has at the
origin, and plenty of oscillations roughly like the ones in the
functions $g_\rho$.  In particular, the oscillations can occur at lots
of different scales as one moves from point to point.  However, one
does not really have oscillations at different scales overlapping each
other.  Whenever the elements of $F$ become dense enough to make a lot
of oscillations, the values of $f$ become small in compensation.  (One
can consider situations where $F$ has points at regularly-spaced
intervals, for instance.)

	How might one prove an estimate like (\ref{bound on average})?
This is part of a larger story in harmonic analysis, called
Littlewood--Paley theory, some of whose classical manifestations are
described in \cite{St1}.  The present discussion is closer in spirit
to \cite{Dr} for the measurements of oscillation used, and indeed
(\ref{bound on average}) can be derived from the results in \cite{Dr}.

	There are stronger estimates available than (\ref{bound on
average}).  Instead of simply counting how often the $\alpha(x,t)$'s
are larger than some threshold, as in the definition of $N_r(x)$
above, one can work with sums of the form
\begin{equation}
\label{q-sum}
	\bigg(\sum_{j=0}^\infty \alpha(x,2^{-j} \, r)^q \bigg)^\frac{1}{q}.
\end{equation}
The ``right'' choice of $q$ is $2$, but to get this one should modify
the definition of $\alpha(x,t)$ (in most dimensions) so that the
measurement of approximation of $f$ by an affine function uses a
suitable $L^p$ norm, rather than the supremum.  (That the choice of
$q=2$ is the ``right'' one reflects some underlying orthogonality, and
is a basic point of Littlewood--Paley theory.  At a more practical
level, $q=2$ is best because it works for the estimates for the
$\alpha(x,t)$'s and allows reverse estimates for the size of the
gradient of $f$ in terms of the sizes of the $\alpha(x,t)$'s.)

	In short, harmonic analysis provides a fairly thorough
understanding of the sizes of the $\alpha(x,t)$'s and related
quantities, and with quantitative estimates.  This works for Lipschitz
functions, and more generally for functions in Sobolev spaces.

	There is more to the matter of smoothness of Lipschitz
functions than this, however.  The $\alpha(x,t)$'s measure how well a
given function $f$ can be approximated by an affine function on a ball
$B(x,t)$, but they do not say too much about how these approximating
affine functions might change with $x$ and $t$.  In fact, there are
classical examples of functions for which the $\alpha(x,t)$'s tend to
$0$ uniformly as $t \to 0$, and yet the derivative fails to exist at
almost every point.  Roughly speaking, the affine approximations keep
spinning around as $t \to 0$, without settling down on a particular
affine function, as would happen when the derivative exists.  (A
faster rate of decay for the $\alpha(x,t)$'s, as in
(\ref{C^{1,delta}}), would prevent this from happening.)

	For Lipschitz functions, the existence of the differential
almost everywhere implies that for almost every $x$ the gradients of
the approximating affine functions on $B(x,t)$ do not have to spin
around by more than a finite amount as $t$ goes between some fixed
number $r$ and $0$.  In fact, quantitative estimates are possible, in
much the same manner as before.  Again one fixes a threshold
$\epsilon$, and one can measure how many oscillations of size at least
$\epsilon$ there are in the gradients of the affine approximations as
$t$ ranges between $0$ and $r$.  For this there are uniform bounds on
the averages of these numbers, just as in (\ref{bound on average}).

	This type of quantitative control on the oscillations of the
gradients of the affine approximations of $f$ comes from Carleson's
Corona construction, as in \cite{Ga}.  This construction was initially
applied to the behavior of bounded holomorphic functions in the unit
disk of the complex plane, but in fact it is a very robust
real-variable method.  For example, the type of bound just mentioned
in the previous paragraph (on the average number of oscillations of
the gradients of the affine approximations as $t$ goes from $r$ to
$0$) is completely analogous to one for the boundary behavior of
harmonic functions given in Corollary 6.2 on p348 of \cite{Ga}.

	A more detailed discussion of the Corona construction in
the context of Lipschitz functions can be found in Chapter IV.2
of \cite{DS5}. 

	The Corona construction and the known estimates for affine
approximations as discussed above provide a fairly complete picture of
the ``smoothness'' of Lipschitz functions.  They also provide an
interesting way to look at ``complexity'' of Lipschitz functions, and
one that is quite different from what is suggested more naively by the
definition (\ref{Lipschitz}).

\subsection{Smoothness and uniform rectifiability}
\label{Smoothness and uniform rectifiability}

	The preceding discussion of smoothness for Lipschitz and
bilipschitz mappings has natural extensions to the geometry of sets in
Euclidean spaces.  Instead of approximations of functions by affine
functions, one can consider approximations of sets by affine planes.
Differentials of mappings correspond to tangent planes for sets.

	One can think of ``embedding'' the discussion for functions
into one for sets by taking a function and replacing it with its
graph.  This is consistent with the correspondence between affine
functions and $d$-planes, and between differentials and tangent planes.

	How might one generalize the $\alpha(x,t)$'s (\ref{def of
alpha(x,t)}) to the context of sets?  Fix a set $E$ in ${\bf R}^n$ and
a dimension $d < n$, and let $x \in E$ and $t > 0$ be given.  In
analogy with (\ref{def of alpha(x,t)}), consider the quantity
$\beta(x,t)$ defined by
\begin{equation}
\label{def of beta}
	\beta(x,t) = \inf_{P \in \mathcal{P}_d} \, 
			\sup \{ t^{-1} \dist(y,P) : y \in E \cap B(x,t) \}.
\end{equation}
Here $\mathcal{P}_d$ denotes the set of $d$-dimensional affine planes in
${\bf R}^n$, and $\dist(y,P)$ is defined as in (\ref{def of dist}).
In other words, we take a ``snapshot'' of $E$ inside the ball
$B(x,t)$, and we look at the optimal degree of approximation of $E$ by
$d$-planes in $B(x,t)$.  The factor of $t^{-1}$ in (\ref{def of beta})
makes $\beta(x,t)$ a scale-invariant, dimensionless quantity.  Notice
that $\beta(x,t)$ is always less than or equal to $1$, no
matter the behavior of $E$, as one can see by taking $P$ to be any
$d$-plane that goes through $x$.  The smoothness of $E$ is reflected
in how small $\beta(x,t)$ is.

	If $E$ is the image of a bilipschitz mapping $\phi : {\bf R}^d
\to {\bf R}^n$, then there is a simple correspondence between the
$\beta(x,t)$'s on $E$ and the $\alpha(z,s)$'s for $\phi$ on ${\bf R}^d$.
This permits one to transfer the estimates for the $\alpha$'s on
${\bf R}^d$ to estimates for the $\beta$'s on $E$, and one could
also go backwards.

	It turns out that the type of estimates that one gets for the
$\beta$'s in this way when $E$ is bilipschitz equivalent to ${\bf
R}^d$ also work when $E$ is uniformly rectifiable.  Roughly speaking,
this because the estimates for the $\alpha$'s and $\beta$'s are not
uniform ones, but involve some kind of integration, and in a way which
is compatible with the measure-theoretic aspects of the Definition
\ref{unif rec}.  This is very much analogous to results in the context
of BMO functions, especially a theorem of Str\"omberg.  (See Chapter
IV.1 in \cite{DS5} for some general statements of this nature.)

	To my knowledge, the first person to look at estimates like
these for sets was P. Jones \cite{Jo1}.  In particular, he used the
sharp quadratic estimates that correspond to Littlewood--Paley theory
to give a new approach to the $L^2$ boundedness of the Cauchy integral
operator on nonsmooth curves.  Here ``quadratic'' means $q=2$ in the
context of (\ref{q-sum}).

	In \cite{Jo3}, Jones showed how quadratic estimates on the
$\beta$'s could actually be used to \emph{characterize} subsets of
rectifiable curves.  The quadratic nature of the estimates, which
come naturally from orthogonality considerations in Littlewood--Paley
theory and harmonic analysis, can, in this context, be more directly
linked to the ordinary Pythagorean theorem, as in \cite{Jo3}.  A
completion of Jones' results for $1$-dimensional sets in Euclidean
spaces of higher dimension was given in \cite{O}.

	Analogues of Jones' results for (Ahlfors-regular) sets of
higher dimension are given in \cite{DS3}.  More precisely, if $E$ is
a $d$-dimensional Ahlfors-regular set in ${\bf R}^n$, then the uniform
rectifiability of $E$ is equivalent to certain quadratic Carleson
measure conditions for quantities like $\beta(x,t)$ in (\ref{def of
beta}).  One cannot use $\beta(x,t)$ itself in general, with the
supremum on the right side of (\ref{def of beta}) (there are
counterexamples due to Fang and Jones), but instead one can replace
the supremum with a suitable $L^p$ norm for a range of $p$'s that
depends on the dimension (and is connected to Sobolev embeddings).
This corresponds to the situation for sharp estimates of quantities
like $\alpha(x,t)$ in the context of Lipschitz functions, as in
\cite{Dr}.

	The problem of building parameterizations is quite different
when $d > 1$ than in the $1$-dimensional case.  This is a basic fact,
and a recurring theme of classical topology.  Making parameterizations
for $1$-dimensional sets is largely a matter of ordering, i.e., lining
up the points in a good way.  For rectifiable curves there is a
canonical way to regulate the ``speed'' of a parameterization, using
arclength.  In higher dimensions none of these things are true,
although conformal coordinates sometimes provide a partial substitute
when $d=2$.  See \cite{DS-trans, MuS, HK1, Se-casscII, Seqs}.  In
\cite{deturck-yang} a different kind of ``normalized coordinates'' are
discussed for $d=3$, but the underlying partial differential equation
is unfortunately not elliptic.  Part of the point of uniform
rectifiability was exactly to try to come to grips with the issue of
parameterizations in higher dimensions.  (See also Appendix \ref{more
on existence and behavior of homeomorphisms} in connection with these
topics.)

	Although this definition (\ref{def of beta}) of $\beta(x,t)$
provides a natural version of the $\alpha(x,t)$'s from (\ref{def of
alpha(x,t)}), it is not the only choice to consider.  There is a
``bilateral'' version, in which one measures both the distance from
points in $E$ to the approximating $d$-plane (as in (\ref{def of
beta})) as well as distances from points in the $d$-plane to $E$.
Specifically, given a set $E$ in ${\bf R}^n$, a point $x \in E$,
a radius $t > 0$, and a $d$-plane $P$ in ${\bf R}^n$, set
\begin{eqnarray}
	{\rm Approx}(E, P, x, t) 
     & = & \sup \{ t^{-1} \dist(y,P) : y \in E \cap B(x,t) \}		\\
	&& + \sup \{ t^{-1} \dist(z,E) : z \in P \cap B(x,t) \} \nonumber
\end{eqnarray}
and then define the bilateral version of $\beta(x,t)$ by
\begin{equation}
	b\beta(x,t) = \inf_{P \in \mathcal{P}_d} {\rm Approx}(E, P, x, t).
\end{equation}
This takes ``holes'' in $E$ into account, which the definition of
$\beta(x,t)$ does not.  For instance, $\beta(x,t) = 0$ if and only if
there is a $d$-plane $P_0$ such that every point in $E \cap B(x,t)$
lies in $P_0$, while for $b\beta(x,t)$ to be $0$ it should also be
true that every point in $P_0 \cap B(x,t)$ lies in $E$ (assuming that
$E$ is closed, as in the definition of Ahlfors regularity).

	It turns out that the $b\beta(x,t)$'s behave a bit differently
from the $\beta(x,t)$'s, in the following sense.  Imagine that we do
not look for something like sharp quadratic estimates, as we did
before, but settle for cruder ``thresholding'' conditions, as
discussed in Subsection \ref{first subsection}.  In other words, one
might fix an $\epsilon > 0$, and define a function $N_r'(x)$ which
counts the number of times that $b\beta(x,2^{-j} \, r)$ is greater
than or equal to $\epsilon$, with $j \in {\bf Z}_+$ (as in the
discussion around (\ref{bad j's})).  For uniformly rectifiable sets
one has bounds on the averages of $N_r'(x)$ exactly as in (\ref{bound
on average}), but now integrating over $E$ instead of ${\bf R}^d$.  A
slightly surprising fact is that the converse is also true, i.e.,
estimates like these for the $b\beta$'s are sufficient to imply the
uniform rectifiability of the set $E$, at least if $E$ is
Ahlfors-regular of dimension $d$.  This was proved in
\cite{DS5}.

	In the context of functions, this type of thresholding
condition is too weak, in that one can have the $\alpha(x,t)$'s going
to $0$ uniformly as $t \to 0$ for functions which are differentiable
almost nowhere, as mentioned in Subsection \ref{first subsection}.
Similarly, there are Ahlfors-regular sets which are ``totally
unrectifiable'' (in the sense of \cite{Fa, Fe, Mat}) and have the
$\beta(x,t)$'s tending to $0$ uniformly as $t \to 0$.  (See
\cite{DS3}.)  For the $b\beta$'s the story is simply different.  On
the other hand, the fact that suitable thresholding conditions on the
$b\beta$'s are sufficient to imply uniform rectifiability relies
heavily on the assumption that $E$ be Ahlfors-regular, while mass
bounds are part of the conclusion (rather than the hypothesis) in
Jones' results, and no counterpart to the mass bounds are included in
the above-mentioned examples for functions.  One does have mass bounds
for the examples in \cite{DS3} (of totally-unrectifiable
Ahlfors-regular sets for which the $\beta(x,t)$'s tend to $0$
uniformly as $t \to 0$), and there the issue is more in the size of
the holes in the set.  The $b\beta$'s, by definition, control the
sizes of holes.  Note that this result for the $b\beta$'s does have
antecedents for the classical notion of (countable) rectifiability, as
in \cite{Mat}.

	There are a number of variants of the $b\beta$'s, in which
one makes comparisons with other collections of sets besides $d$-planes,
like unions of $d$-planes, for instance.  See \cite{DS5}.

	Perhaps the strongest formulation of smoothness for uniformly
rectifiable sets is the existence of a ``Corona decomposition''.  This
is a geometric version of the information that one can get about a
Lipschitz function from the methods of Carleson's Corona construction
(as mentioned in Subsection \ref{first subsection}).  Roughly
speaking, in this condition one controls not only how often $E$ is
well-approximated by a $d$-plane, but how fast the $d$-planes turn as
well.  This can also be formulated in terms of good approximations of
$E$ by flat Lipschitz graphs.

	Although a bit technical, the existence of a Corona
decomposition is perhaps the most useful way of managing the
complexity of a uniformly rectifiable set.  Once one has a Corona
decomposition, it is generally pretty easy to derive whatever else one
would like to know.  Conversely, in practice the existence of a Corona
decomposition can be a good place to start if one wants to prove that
a set is uniformly rectifiable.

	In fact, there is a general procedure for finding a Corona
decomposition when it exists, and one which is fairly simple (and very
similar to Carleson's Corona construction).  The difficult part is to
show that this procedure works in the right way, with the correct
estimates.  Specifically, it is a stopping-time argument, and one does
not want to have to stop too often.  This is a nice point, because in
general it is not so easy to build something like a good
parameterization of a set, even if one knows a priori that it exists.
In this context, there are in principle methods for doing this.

	See \cite{DS-pubmat, DS3, DS5, DS-compg, Se-ind} for more
information about Corona decompositions of uniformly rectifiable sets
and the way that they can be used.  The paper \cite{DS-compg} is
written in such a way as to try to convey some of the basic concepts
and constructions without worrying about why the theorems are true
(which is much more complicated).  In particular, the basic procedure
for finding Corona decompositions when they exist is discussed.  See
\cite{GaJ, Jo1, Jo3} for some other situations in which Carleson's
Corona construction is used geometrically.

\subsection{A class of variational problems}
\label{A class of variational problems}

	Uniform rectifiability is a pretty robust condition.  If one
has a set which looks roughly as though it ought to be uniformly
rectifiable, then there is a good chance that it is.  This as opposed
to sets which look roughly as though they should admit a well-behaved
(homeomorphic) parameterization, and do not (as discussed before).

	In this subsection we would like to briefly mention a result
of this type, concerning a minimal surface problem with nonsmooth
coefficients.  Let $g(x)$ be a Borel measurable function on ${\bf
R}^n$, and assume that $g$ is positive, bounded, and bounded away from
$0$, so that
\begin{equation}
	0 < m \le g(x) \le M
\end{equation}
for some constants $m$, $M$ and all $x \in{\bf R}^n$.  Let $Q_0$, $Q_1$
be a pair of (closed) cubes in ${\bf R}^n$, with sides parallel to the
axes, and assume that $Q_0$ is contained in the interior of $Q_1$.

	Let $U$ be an open subset of $Q_1$ which contains the interior
of $Q_0$.  Consider an integral like
\begin{equation}
\label{functional}
	\int_{\partial U} g(x) \, d\nu_U(x),
\end{equation}
where $d\nu_U$ denotes the measure that describes the
$(n-1)$-dimensional volume of subsets of $\partial U$.  This would be
defined as in calculus when $\partial U$ is at least a little bit
smooth (like $C^1$), but in general one has to be more careful.  One
can simply take for $d\nu_U$ the restriction of $(n-1)$-dimensional
Hausdorff measure to $\partial U$, but for technical reasons it is
often better to define $d\nu_U$ using distributional derivatives of
the characteristic function of $U$, as in \cite{Giusti}.  For this one
would work with sets $U$ which have ``finite perimeter'', which means
exactly that the distributional first derivatives of the
characteristic function of $U$ are measures of finite mass. 

	Here is one way in which this kind of functional, and the
minimization of this kind of functional, can come up.  Let $F$ be a
closed subset of ${\bf R}^n$.  Imagine that one is particularly
interested in domains $U$ which have their boundary contained in $F$,
or very nearly so.  On the other hand, one might also wish to limit
irregularities in the behavior of the boundary of $U$.  For this type
of situation one could choose $g$ so that it is much smaller on $F$
than on the complement of $F$, and then look for minimizers of
(\ref{functional}) to find domains with a good balance between the
behavior of $\partial U$ and the desire to have it be contained (as
much as possible) in $F$.  (See \cite{DS8} for an example of this.)

	When do minimizers of (\ref{functional}) exist, and how do
they behave?  If one works with sets of finite perimeter, and if the
function $g$ is lower semi-continuous, then one can obtain the
existence of minimizers through standard techniques (as in
\cite{Giusti}).  That is, one takes limits of minimizing sequences for
(\ref{functional}) using weak compactness, and one uses the lower
semi-continuity of $g$ to get lower semi-continuity of
(\ref{functional}) with respect to suitable convergence of the $U$'s.
The latter ensures that the limit of the minimizing sequence is
actually a minimum.  Note that the ``obstacle'' conditions that $U$
contain the interior of $Q_0$ and be contained in $Q_1$ prevents the
minimization from collapsing into something trivial.

	As to the behavior of minimizers of (\ref{functional}), one
cannot expect much in the way of smoothness in general.  For instance,
if the boundary of $U$ can be represented locally as the graph of a
Lipschitz function, then $U$ in fact minimizes (\ref{functional}) for
a suitable choice of $g$.  Specifically, one can take $g$ to be a
sufficiently small positive constant on $\partial U$, and to be equal
to $1$ everywhere else.  That such a choice of $g$ works is not very
hard to establish, and more precise results are given in \cite{DS8}.

	Conversely, minimizers of (\ref{functional}) are always
Ahlfors-regular sets of dimension $n-1$, and uniformly rectifiable.
This is shown in \cite{DS8}, along with some additional
geometric information which is sufficient to characterize the class of
sets $U$ which occur as minimizers for functionals of the form
(\ref{functional}) (with $g$ bounded and bounded away from $0$).  If a
set $U$ arises as the minimizer for some $g$, it is also a minimizer
with $g$ chosen as above, i.e., a small positive constant on $\partial
U$ and equal to $1$ everywhere else.

	The same regularity results work for a suitable class of
``quasiminimizers'' of the usual area functional, and one that
includes minimizers for (\ref{functional}) as a special case.

	Uniform rectifiability provides a natural level of structure
for situations like this, where stronger forms of smoothness cannot be
expected, but quantitative bounds are reasonable to seek.  Note that
properties of ordinary rectifiability always hold for boundaries of
sets of finite perimeter, regardless of any minimizing or
quasiminimizing properties.  See \cite{Giusti}.

	Analogous results about regularity work for sets of higher
codimension as well, although this case is more complicated
technically.  See \cite{DS-q-min} for more information.  One can use
this framework of minimization (with respect to nonsmooth coefficient
functions $g$) as a tool for studying the structure of sets in ${\bf
R}^n$ with upper bounds on their $d$-dimensional Hausdorff measure and
lower bounds for their $d$-dimensional topology.  This brings one back
to Theorem \ref{unif rec theorem} and related questions, and in
particular more ``localized'' versions of it.

	To put it another way, minimization of functionals like these
can provide useful means for obtaining ``existence results'' for
approximate parameterizations with good behavior, through uniform
rectifiability.  See \cite{DS8, DS-q-min}.  Part of the motivation for
this came from an earlier argument of Morel and Solimini
\cite{morel-solimini}.  Their argument concerned the existence of
curves containing a given set, with good properties in terms of the
distribution of the arc-length measure of these curves, under more
localized conditions on the given set (at all locations and scales).
See Lemma 16.27 on p207 of \cite{morel-solimini}.

\appendix
\addtocontents{toc}{{\protect\bigskip\protect\hbox{\protect\bf\protect\large
Appendices}}}

\section{Fourier transform calculations}
\label{Some Fourier transform calculations}
\setcounter{equation}{0}

\renewcommand{\thetheorem}{A.\arabic{equation}}
\renewcommand{\theequation}{A.\arabic{equation}}

	If $\phi(x)$ is an integrable function on ${\bf R}^n$, then
its \emph{Fourier transform} $\widehat{\phi}(\xi)$ is defined (for
$\xi \in {\bf R}^n$) by
\begin{equation}
	\widehat{\phi}(\xi) 
		= \int_{{\bf R}^n} e^{i \, \langle x , \xi \rangle} 
							\, \phi(x) \, dx.
\end{equation}
Here $\langle x , \xi \rangle$ denotes the usual inner product for $x,
\xi \in {\bf R}^n$, and $i = \sqrt{-1}$.  Often one makes slightly
different conventions for this definition --- with some extra factors
of $\pi$ around, for instance --- but we shall not bother with this.

	A key feature of the Fourier transform is that it diagonalizes
differential operators.  Specifically, if $\partial_k$ denotes the
operator $\partial/\partial x_k$ on ${\bf R}^n$, then
\begin{equation}
\label{diagonalization identity}
  (\partial_k \phi)^{\widehat{}} \, (\xi) = i \, \xi_k \, \widehat{\phi}(\xi),
\end{equation}
i.e., differentiation is converted into mere multiplication.  For this
one should either make some differentiability assumptions on $\phi$,
so that the left side can be defined in particular, or one should
interpret this equation in the sense of tempered distributions on
${\bf R}^n$.  The Fourier transform also carries out this
diagonalization in a controlled manner.  That is, there is an explicit
inversion formula (which looks a lot like the Fourier transform
itself), and the Fourier transform preserves the $L^2$ norm of the
function $\phi$, except for a multiplicative constant, by the
Plancherel theorem.  See \cite{SW} for these and other basic facts
about the Fourier transform.

	Using Plancherel's theorem, it is very easy to give another
proof of the $L^2$ estimate (\ref{L^2 estimate}) from Section \ref{BMO
frame of mind}, and to derive many other inequalities of a similar
nature.  One can also use the Fourier transform to give a precise
definition of the operator $R = \partial_j \partial_k / \Delta$, where
$\Delta$ is the Laplace operator $\sum_{\ell = 1}^n \partial_\ell^2$.
Specifically, one can define it through the equation
\begin{equation}
\label{specific m(xi)}
	(R\phi)^{\widehat{}} \, (\xi) 	
		= {\xi_j \xi_k \over |\xi|^2} \, \widehat{\phi}(\xi).
\end{equation}

	If $m(\xi)$ is any bounded function on ${\bf R}^n$, then
\begin{equation}
\label{operator associated to m}
	(T\phi)^{\widehat{}} \, (\xi) 	
		= m(\xi) \, \widehat{\phi}(\xi)
\end{equation}
defines a bounded operator on $L^2({\bf R}^n)$.  In general these
operators are not bounded on $L^p$ for any other value of $p$, but
this is true for many of the operators that arise naturally in
analysis.  For instance, suppose that $m(\xi)$ is homogeneous of
degree $0$, so that
\begin{equation}
	m(t \xi) = m(\xi) \qquad \hbox{when } t > 0,
\end{equation}
and that $m(\xi)$ is smooth away from the origin.  Then the associated
operator $T$ is bounded on $L^p$ for all $p$ with $1 < p < \infty$.
See \cite{St1, SW}.  Note that this criterion applies to the specific
choice of $m(\xi)$ in (\ref{specific m(xi)}) above.

	For a multiplier operator as in (\ref{operator associated to
m}) to be bounded on $L^1$ or $L^\infty$ is even more exceptional than
for $L^p$ boundedness when $1 < p < \infty$.  (See \cite{SW}.)  For
instance, if $m$ is homogeneous, as above, and not constant, then the
corresponding operator cannot be bounded on $L^1$ or $L^\infty$.
However, if $m$ is homogeneous and smooth away from the origin, then
the operator $T$ in (\ref{operator associated to m}) does determine a
bounded operator from $L^\infty$ into BMO.  See \cite{GCRF, Ga,
Journe, St2}.  In fact, $T$ determines a bounded operator from BMO to
itself.

	Here is another example.  Let $\phi$ now be a mapping from
${\bf R}^2$ to itself, with components $\phi_1$, $\phi_2$.  Consider
the differential $d \phi$ of $\phi$ as a matrix-valued function, namely,
\begin{equation}
	{\partial_1 \phi_1 \enspace \partial_1 \phi_2 \choose
			\partial_2 \phi_1 \enspace \partial_2 \phi_2}.
\end{equation}
(Let us assume that $\phi$ is smooth enough that the differential is
at least some kind of function when taken in the sense of
distributions, although one can perfectly well think of $d\phi$ as a
matrix-valued distribution.)  Let $A$ and $S$ denote the antisymmetric
and symmetric parts of $d\phi$, respectively, so that
\begin{equation}
	A = {d\phi - d\phi^t \over 2}, \qquad S = {d\phi + d\phi^t \over 2},
\end{equation}
where $d\phi^t$ denotes the transpose of $d\phi$.  

	In this case of $2 \times 2$ matrices, the antisymmetric part $A$
really contains only one piece of information, namely
\begin{equation}
	\partial_1 \phi_2 - \partial_2 \phi_1.
\end{equation}
It is not hard to check that this function can be reconstructed from
the entries of $S$ through operators of the form (\ref{operator
associated to m}), using functions $m$ which are homogeneous of degree
$0$ and smooth away from the origin.  For this one should add some
mild conditions on $\phi$, like compact support, to avoid the possibility
that $S$ vanishes identically but $A$ does not.  

	Under these conditions, we conclude that the $L^p$ norm of $A$
is always bounded by a constant multiple of the $L^p$ norm of $S$, $1
< p < \infty$, and that the BMO norm of $A$ is controlled by the
$L^\infty$ (or BMO) norm of $S$.  (For the case of BMO norms, the
possibility that $S$ vanishes but $A$ does not causes no trouble,
because $A$ will be constant in that case.)

	This example is really a ``linearized'' version of the problem
discussed in Section \ref{Mappings and distortion}.  Specifically, let
us think of $f : {\bf R}^2 \to {\bf R}^2$ as being of the form
\begin{equation}
	f(x) = x + \epsilon \, \phi(x),
\end{equation}
where $\epsilon$ is a small parameter.  The extent to which $f$ distorts
distances is governed by the matrix-valued function $df^t \, df$,  which
we can write out as
\begin{equation}
	df^t \, df = I + 4 \, \epsilon \, S + \epsilon^2 \, d\phi^t d\phi.
\end{equation}
Thus the linear term in $\epsilon$ is governed by $S$, while $A$
controls the leading behavior in $\epsilon$ of the ``rotational'' part
of $df$.

\section{Mappings with branching}
\label{branching}
\setcounter{equation}{0}

\renewcommand{\thetheorem}{B.\arabic{equation}}
\renewcommand{\theequation}{B.\arabic{equation}}

	In general, there can be a lot of trouble with existence and
complexity of homeomorphisms (with particular properties, like
specified domain and range).  If one allows mappings with
\emph{branching}, then the story can be very different.  

	As a basic example of this, there is a classical result
originating with Alexander to the effect that any oriented
pseudomanifold of dimension $n$ admits an orientation-preserving
branched covering over the $n$-sphere.  Let us state this more
carefully, and then see how it is proved.

	Let $M$ be a finite polyhedron.  We assume that $M$ is given
as a finite union of $n$-dimensional simplices that meet only in their
boundary faces (so that $M$ is really a simplicial complex).  To be a
pseudomanifold means that every $(n-1)$-dimensional face in $M$ arises
as the boundary face of exactly \emph{two} $n$-dimensional simplices.
In effect this says that $M$ looks like a manifold away from its
codimension-$2$ skeleton (the corresponding statement for the
codimension-$1$ skeleton being automatic).  For the present purposes
it would be enough to ask that every $(n-1)$-dimensional face in $M$
arise as the boundary face of \emph{at most} two $n$-dimensional
simplices, which would be like a ``pseudomanifold with boundary''.

	An \emph{orientation} for an $n$-dimensional pseudomanifold
$M$ means a choice of orientation (in the usual sense) for each of the
constituent $n$-dimensional simplices in $M$, with compatibility of
orientations of adjacent $n$-dimensional simplices along the common
$(n-1)$-dimensional face.  In terms of algebraic topology, this means
that the sum of the $n$-dimensional simplices in $M$, with their
orientations, defines an $n$-dimensional \emph{cycle} on $M$. 

	For the purposes of the Alexander-type result, it will be
convenient to think of the $n$-sphere as consisting of two standard
simplices $S_1$ and $S_2$ glued together along the boundary.  This is
not quite a polyhedron in the usual (affine) sense, but one could
easily repair this by subdividing $S_1$ or $S_2$.  We also assume that
$S_1$ and $S_2$ have been oriented, and have opposite orientations
relative to their common boundary.

	To define a mapping from $M$ to the $n$-sphere one would like
to simply identify each of the constituent $n$-dimensional simplices
in $M$ with $S_1$ or $S_2$ in a suitable manner.  Unfortunately, this
does not work, even when $n=2$, but the problem can be fixed using a
\emph{barycentric subdivision} of $M$.  Recall that the
\emph{barycenter} of a simplex (embedded in some vector space) is the
point in the interior of the simplex which is the average of the
vertices of the simplex.  The \emph{set of barycenters} for $M$ means
the set of barycenters of all of the constituent simplices in $M$
(viewed as a simplicial complex), of all dimensions, including $0$.
In particular, the set of barycenters for $M$ includes the vertices of
$M$ (which are themselves $0$-dimensional simplices, and their own
barycenters).  The barycentric subdivision of $M$ is a refinement of
$M$ as a simplicial complex whose vertices are exactly the set of
barycenters of $M$.  In other words, the set $M$ as a whole does not
change, just its decomposition into simplices, which is replaced by a
finer decomposition.

	Here is a precise description of the simplices in the
barycentric subdivision of $M$.  Let $s_0, s_1, \ldots, s_k$ be a
finite sequence of simplices in $M$, with each $s_i$ an
$i$-dimensional simplex which is a face of $s_{i+1}$ (when $i < k$).
Let $b(s_i)$ denote the barycenter of $s_i$.  Then $b(s_0), b(s_1),
\ldots, b(s_k)$ are affinely independent, and hence determine a
$k$-dimensional simplex.  The simplices that arise in this manner are
precisely the ones used for the barycentric subdivision of $M$.  (See
p123 of \cite{Sp} for more details.)

	Let $\widetilde{V}$ denote the set of all vertices in the
barycentric subdivision of $M$.  This is the same as the set of points
which arise as barycenters of simplices in the original version of
$M$, and in particular we have a natural mapping from $\widetilde{V}$
to the integers $\{0, 1, \ldots, n\}$, defined by associating to each
point $b$ in $\widetilde{V}$ the dimension of the simplex from which
it was derived.

	If $T$ is a $k$-dimensional simplex in the barycentric
subdivision of $M$, then the mapping from $\widetilde{V}$ to $\{0, 1,
\ldots, n\}$ just described induces a one-to-one correspondence
between the $k+1$ vertices of $T$ and the set $\{0, 1, \ldots, k\}$.
This follows easily from the definitions.

	We are now ready to define our mapping from $M$ to the
$n$-sphere.  There are exactly $n+1$ vertices in our realization of
the $n$-sphere as the gluing of $S_1$ and $S_2$.  Let us identify
these vertices with the integers from $0$ to $n$.  Thus our mapping
from $\widetilde{V}$ to $\{0, 1, \ldots, n\}$ can now be interpreted
as a mapping from the vertices of the barycentric subdivision of $M$
to the vertices of the $n$-sphere.

	This mapping between vertices admits a canonical linear
extension to each $k$-dimensional simplex, $k < n$, in the barycentric
subdivision of $M$.  For the $n$-dimensional simplices the extension
is uniquely determined once one chooses $S_1$ or $S_2$ for the image
of the simplex.  Because of the orientations, there is only one
natural choice of $S_1$ or $S_2$ for each $n$-dimensional simplex $T$,
namely the one so that the linear mapping from $T$ onto $S_j$ is
orientation-preserving.

	In the end we get a mapping from the barycentric subdivision
of $M$ to the $n$-sphere which preserves orientations and which
defines an affine isomorphism from each $n$-dimensional simplex $T$ in
the domain onto one of $S_1$ and $S_2$.  This uses the fact that our
initial mapping between vertices was always one-to-one on the set of
vertices in any given simplex in the domain, by construction.

	This completes the proof.  We should emphasize that the
singularities of the mapping from $M$ to the $n$-sphere --- i.e., the
places where it fails to be a local homeomorphism --- are confined to
the codimension-$2$ skeleton of the barycentric subdivision of $M$.
This is because of the orientation and pseudomanifold conditions,
which ensure that if a point in $M$ lies in the interior of an
$(n-1)$-dimensional simplex in the barycentric subdivision of $M$,
then the (two) adjacent $n$-simplices at that point are not sent to
the same $S_j$ in the image.

	The idea of branching also makes sense for mappings that are
not piecewise-linear, and there are well-developed notions of
``controlled geometry'' in this case, as with the classes of
quasiregular mappings and mappings of bounded length distortion.  See
\cite{HKM, MaRV1, MaRV2, MaRV, MaV, Res, Ri1, V-icm, Vuorinen}, for
instance.  In \cite{HR1, HR2} there are examples where branching maps
of controlled geometry can be constructed but suitable homeomorphisms
either do not exist or must distort distances more severely.

	Sullivan \cite{Su2, Su3} has proposed some mechanisms by which
the existence of local (controlled) branching maps can be deduced, and
some ideas for studying obstructions to controlled homeomorphic
coordinates.

	  See \cite{GMRV, HKi, MRyVu} for some recent results about
branching and regularity conditions under which it does not occur.  A
broader and more detailed discussion of mappings with branching can be
found in \cite{HR2}.  For some real-variable considerations of
mappings which may branch but enjoy substantial geometric properties,
see \cite{Da3, Jo2, DS-trans}.

\section{More on existence and behavior of homeomorphisms}
\label{more on existence and behavior of homeomorphisms}
\setcounter{equation}{0}

\renewcommand{\thetheorem}{C.\arabic{equation}}
\renewcommand{\theequation}{C.\arabic{equation}}

\subsection{Wildness and tameness phenomena}
\label{Wildness and tameness phenomena}

	Consider the following question.  Let $n$ be a positive integer,
and let $K$ be a compact subset of ${\bf R}^n$.  If $K$ is homeomorphic
to the unit interval $[0,1]$, is there 
\begin{eqnarray}
\label{tameness condition}
	&& \hbox{a global homeomorphism from ${\bf R}^n$ onto itself}	
                                                                        \\
	&& \hbox{which maps $K$ to a straight line segment?}        
                                                                 \nonumber
\end{eqnarray}

	If $n=1$, then $K$ itself is a closed line segment, and the
answer is ``yes''.  When $n=2$, the answer is also ``yes'', but this
is more complicated, and is more in the spirit of the Sch\"onflies theorem
in the plane.  See \cite{Moise}, especially Chapter 10.

	When $n \ge 3$, the answer to the question above can be
``no''.  An arc $K$ is said to be ``tame'' (or flat) when a
homeomorphism does exist as in (\ref{tameness condition}), and
``wild'' when it does not exist.  See \cite{Moise} for some examples
of wild arcs in ${\bf R}^3$.

	Smooth arcs are always tame, as are polygonal arcs, i.e., arcs
made up of \emph{finitely} many straight line segments.  For these one
can take the corresponding homeomorphism to be smooth or
piecewise-linear as well.  (Compare with Theorem 1 on p134 of
\cite{Moise}, for instance.)  In order for an arc to be wild, some
amount of infinite processes are needed.

	A simple closed curve in ${\bf R}^3$ might be smooth or
polygonal and still \emph{knotted}, so that there does not exist a
homeomorphism of ${\bf R}^3$ onto itself which maps the curve onto a
standard circle (inside a standard $2$-dimensional plane in ${\bf
R}^3$).  There are many well-known examples of this, like the trefoil
knot.  Thus, for a closed curve, one defines ``wildness'' in a
slightly differently way, in terms of the existence of local
flattenings, for instance.  This turns out to be compatible with the
case of arcs (for which there is no issue of knottedness), and there
are some other natural variants of this.

	Here is another basic example, for sets of higher dimension.
Suppose that $\gamma$ is a simple closed curve in ${\bf R}^3$, which
is a polygonal curve, and which represents the trefoil knot.  Consider
the cone over $\gamma$, which gives a $2$-dimensional polyhedron in
${\bf R}^4$, and which is in fact piecewise-linearly equivalent to a
standard $2$-dimensional cell.  One can show that this embedding of
the $2$-cell is not locally flat at the cone point, i.e., it cannot be
straightened out to agree with a standard (geometrically flat)
embedding by a homeomorphism defined on a neighborhood in ${\bf R}^4$
of the cone point.  Similar phenomena occur for codimension-$2$
embeddings in ${\bf R}^n$ for all $n \ge 4$, as in Example 2.3.2 on
p59-60 of \cite{Ru-book}.

	This phenomenon is special to codimension $2$, however.
Specifically, a piecewise-linear embedding of a $k$-dimensional
piecewise-linear manifold into ${\bf R}^n$ is locally topologically
flat if $n - k \ne 2$ (or if $k = 1$ and $n = 3$, as before).  See
Theorem 1.7.2 on p34 of \cite{Ru-book}.

	In the context of piecewise-linear embeddings, one can also
look for local flattenings which are piecewise-linear.  A similar
remark applies to other categories of mappings.  We shall not pursue
this here.

	Wild embeddings of cells and spheres (and other manifolds)
exist in ${\bf R}^n$ for all $n \ge 3$, and for all dimensions of the
cells and spheres (from $1$ to $n-1$).  This includes embeddings of
cells and spheres which are not equivalent to piecewise-linear
embeddings in codimension $2$.  We shall mostly consider here issues
of existence of topological flattenings or local flattenings, and
embeddings which are not normally given as piecewise-linear.

	See \cite{Bing-book, Bu, BuC, C1, Dm1, Dm2, E-demension,
Moise, Ru-book, Ru} for more information, and for further references.
Let us also mention that embeddings, although wild, may still enjoy
substantial good behavior.  For instance, they may be bilipschitz, as
in (\ref{k-bilipschitz}), or quasisymmetric, in the sense of
\cite{TV}.  (Roughly speaking, an embedding is \emph{quasisymmetric}
if relative distances are approximately preserved, rather than
distances themselves, as for a bilipschitz mapping.)  See
\cite{Ge-icm1, LuV, V1} for some basic results about this.

	As another version of wildness for embeddings, imagine that
one has a compact set $C$ in some ${\bf R}^n$, and that $C$ is
homeomorphic to the usual middle-thirds Cantor set.  Can one move $C$
to a subset of a straight line in ${\bf R}^n$, through a homeomorphism
from ${\bf R}^n$ onto itself?

	When $n=1$ this is automatically true.  It is also true when
$n=2$; see \cite{Moise}, especially Chapter 13.  In higher dimensions
it is not true in general, as is shown by a famous construction of
Antoine (``Antoine's necklaces'').  See Chapter 18 of \cite{Moise} and
\cite{Blankinship}.

	How can one tell when a set is embedded wildly or not?  As a
simple case, let us consider Cantor sets.  If $C$ is a compact subset
of ${\bf R}^n$ which lies in a line and is homeomorphic to the Cantor
set, and if $n$ is at least $3$, then the complement of $C$ in ${\bf
R}^n$ is simply-connected.  This is not hard to see.  Basically, if
one takes a loop in the complement of $C$ and fills it with a disk in
${\bf R}^n$, and if that disk happens to run into $C$, then one can
make small perturbations of the disk to avoid intersecting $C$.

	The complement of $C$ is also simply-connected if there is a
global homeomorphism from ${\bf R}^n$ onto itself which maps $C$ into
a line.  This is merely because the homeomorphism itself permits one
to reduce to the previous case.

	However, Antoine's necklaces have the property that their
complements are \emph{not} simply-connected.  See \cite{Moise,
Blankinship}.  Note that the \emph{homology} of the complement of a
compact set in ${\bf R}^n$ is controlled through the intrinsic
topology of the set itself, as in Alexander duality \cite{Sp}.  In
particular, while the complement of an Antoine's necklace may not be
simply-connected, its $1$-dimensional homology does vanish.

	Versions of the fundamental group play an important role for
wildness and taming in general, and not just for Cantor sets.  For
this one may not take (or want to take) the fundamental group of the
whole complementary set, but look at more localized forms of the
fundamental group.  A specific and basic version of this is the
following.  Suppose that $F$ is a closed set inside of some ${\bf
R}^n$.  Given a point $p \in F$, and a loop $\gamma$ in ${\bf R}^n
\backslash F$ which lies close to $p$, one would like to know whether
it is possible to contract $\gamma$ to a point in ${\bf R}^n
\backslash F$ while staying in a small neighborhood of $p$.  This
second neighborhood of $p$ might not be quite as small the first one;
a precise statement would say that for every $\epsilon > 0$ there is a
$\delta > 0$ so that if $\gamma$ lies in $B(p,\delta) \cap ({\bf R}^n
\backslash F)$, then $\gamma$ can be contracted to a point in
$B(p,\epsilon) \cap ({\bf R}^n \backslash F)$.  

	This type of condition is satisfied by standard embeddings of
sets into ${\bf R}^n$, like Cantor sets, cells, and spheres, at least
when the dimension of the set is different from $n-2$.  For a point in
${\bf R}^2$, or a line segment in ${\bf R}^3$, etc., one would get
${\bf Z}$ for the corresponding localized fundamental group of the
complement of the set.  (In the case of a line segment in ${\bf R}^3$,
one should restrict one's attention to points $p$ in the interior of
the segment for this.)  

	Conversely, there are results which permit one to go
backwards, and say that localized fundamental group conditions for the
complement like these (localized simple-connectedness conditions in
particular) lead to tameness of a given set, or other kind of
``standard'' (non-wild) behavior.  See \cite{Bing5, Bing-book,
Bing-coll, Bu, BuC, C1, C-bulletin, Dm1, Dm2, E-demension, Moise,
Quinn-bull, QuinnI, Ru-book, Ru} for more information about localized
fundamental groups and their role in wildness phenomena and taming
theorems (and for related matters and further references).

	Fundamental groups and localized versions of them have a basic
role in geometric topology in general.  Some aspects of this came up
before in Section \ref{parameterization problems}, and we shall
encounter some more in this appendix.  

	One might wonder why $\pi_1$ and localized versions of it play
such an important role.  Some basic points behind this are as follows.
For homology (or cohomology), one often has good information from data
in the given situation through standard results in algebraic topology,
like duality theorems.  This definitely does not work for $\pi_1$, as
shown by many examples, including the ones mentioned above, and others
later in this appendix.  In circumstances with suitable
simple-connectivity, one can pass from information about homology to
information about homotopy (in general dimensions), as in the Hurewicz
and Whitehead theorems.  See \cite{Bredon, Sp}.  For many constructions,
homotopy is closer to what one really needs.

	Another basic point concerns the effect of
\emph{stabilization}.  A wild embedding of a set into some ${\bf R}^n$
can become tame when viewed as an embedding into an ${\bf R}^m$ with
$m > n$ ($m = n+1$ in particular).  The same is true for knotting.
For instance, a smooth loop may be knotted in ${\bf R}^3$, but when
viewed as a subset of ${\bf R}^4$, it is always unknotted.  This is
easy to see in explicit examples (like a trefoil knot).

	For some simple and general results about wild embeddings in
${\bf R}^n$ becoming tame in a larger ${\bf R}^m$, see Proposition 4
on p84 of \cite{Dm2}, and the corollaries on p85 of \cite{Dm2}.  These
involve a famous device of Klee.  In concrete examples, one can often
see the taming in a larger-dimensional space directly, and explicitly.
Examples of wild sets are often made with the help of various
linkings, or something like that, and in a higher-dimensional space
one can disentangle the linked parts.  This can be accomplished by
taking individual pieces and pulling them into a new dimension, moving
them around freely there, and then putting them back into the original
${\bf R}^n$ in a different way.

	This simplifying effect of stabilization also fits with the
role of localized versions of fundamental groups indicated before.
Let $F$ be a closed set inside of some ${\bf R}^n$, and imagine that
one has a loop $\gamma$ in ${\bf R}^n \backslash F$ which lies in a
small ball centered at a point in $F$.  In the condition that was
discussed earlier, one would like to contract $\gamma$ to a point in
the complement of $F$, while remaining in a small ball.  If one thinks
of $\gamma$ and $F$ as being also inside ${\bf R}^{n+1}$, then it is
easy to contract $\gamma$ to a point in the complement of $F$ in ${\bf
R}^{n+1}$, while remaining in a small ball.  Specifically, one can
first translate $\gamma$ into a parallel copy of ${\bf R}^n$ inside of
${\bf R}^{n+1}$, i.e., into ${\bf R}^n \times \{a\}$ for some $a \ne
0$ rather than ${\bf R}^n \times \{0\}$ in ${\bf R}^{n+1}$ (using the
obvious identifications).  This parallel copy is then disjoint from
$F$, and one can contract the loop in a standard way.  Note that this
argument works independently of the behavior of $F$.

	Using taming theorems based on localized fundamental group
conditions, and considerations like those in the previous paragraph,
one can get stronger results on tameness that occurs from
stabilization than the ones on p84-85 in \cite{Dm2} mentioned above.
More precisely, instead of needing $k$ extra dimensions in some cases,
it is enough to go from ${\bf R}^n$ to ${\bf R}^{n+1}$.  Compare with
the bottom of p390 and the top of p391 in \cite{Dm1}, and the
references indicated there.  (Compare also with the remarks on p452
of \cite{C1}.)

	At any rate, this type of phenomenon, of objects becoming more
``tame'' or simple after stabilization, is a very basic one in
geometric topology (as well as other areas, for that matter).  We
shall encounter a number of other instances of this in this appendix.
As in the present setting, the effect of stabilization is also
frequently related to conditions concerning localized fundamental
groups.  I.e., such conditions often become true, and in a simple way,
after stabilization.

\subsection{Contractable open sets}
\label{Contractable open sets}

	Fix a positive integer $n$.
\begin{eqnarray}
\label{question about whether a contractable open set is a topological ball}
    &&  \hbox{If $U$ is a nonempty contractable open subset of ${\bf R}^n$,}
                                                                          \\
    &&  \hbox{is $U$ necessarily homeomorphic to the open unit ball in }
                                                               {\bf R}^n?
                                                                   \nonumber
\end{eqnarray}
For the record, to say that $U$ is contractable means that the
identity mapping on $U$ is homotopic to a constant, through
(continuous) mappings from $U$ into itself.  In particular, the
homotopy and homology groups of $U$ (of positive dimension) would then
vanish, just as for an $n$-dimensional ball.

	When $n=1$, the answer to the question in (\ref{question about
whether a contractable open set is a topological ball}) is ``yes''.
In this case, $U$ is either the whole real line, an open segment in
the real line, or an open ray.  Each of these is easily seen to be
homeomorphic to the interval $(-1,1)$, which is the unit ball in this
case.

	If $n=2$, then the answer to the question in (\ref{question
about whether a contractable open set is a topological ball}) is
``yes'' again.  This is a well-known fact, and we shall return to it
later, in Subsection \ref{Geometric and analytic results about the
existence of good coordinates}.

	Starting in dimension $3$, the answer to the question in
(\ref{question about whether a contractable open set is a topological
ball}) is ``no''.  We shall say something about examples for this in a
moment, but let us first ask ourselves the following: how might one be
able to tell that a given contractable open set in ${\bf R}^n$ is
\emph{not} homeomorphic to an $n$-dimensional ball?

	Here again a localized version of the fundamental group is
important.  If $n \ge 3$, then a necessary condition for a set $U$ to
be homeomorphic to an $n$-dimensional ball is that $U$ be ``simply
connected at infinity''.  Roughly speaking, this means that if one
takes a closed loop $\gamma$ out near infinity in $U$, then it should
be possible to contract $\gamma$ to a point, while staying out near
infinity too (although perhaps not as much as $\gamma$ itself is).

	Here is a more formal definition.  For this we also include
``connectedness at infinity'' as a first part.

\begindefinition
\label{def of simple connectivity at infinity}
Let $U$ be an open set in ${\bf R}^n$, or a topological space more
generally.  (Normally one might at least ask that $U$ be locally
compact.)

	$U$ is \emph{connected at infinity} if for each compact set
$K_0 \subseteq U$ there is a larger compact set $L_0 \subseteq U$ such
that every pair of points in $U \backslash L_0$ is contained in a
connected set which is itself contained in $U \backslash K_0$.  (One
can define ``arcwise connectedness at infinity'' in a similar manner.
The two notions are equivalent under assumptions of local arcwise
connectedness, and for topological manifolds in particular.)

	$U$ is \emph{simply-connected at infinity} if it is connected
at infinity, and if for every compact set $K_1 \subseteq U$ there is a
larger compact set $L_1 \subseteq U$ so that if $\gamma$ is an
arbitrary closed loop in $U \backslash L_1$ (i.e., an arbitrary
continuous mapping from the unit circle ${\bf S}^1$ into $U \backslash
L_1$), then $\gamma$ is homotopic to a constant through continuous
mappings from the circle into $U \backslash K_1$.
\end{definition}

	If $U$ is the unit ball in ${\bf R}^n$, $n \ge 3$, then $U$ is
simply-connected at infinity.  Indeed, let $B(0,r)$ denote the open
ball in ${\bf R}^n$ with center $0$ and radius $r$, and let
$\overline{B}(0,r)$ denote the corresponding closed ball.  Then every
compact subset of $B(0,1)$ is contained in $\overline{B}(0,r)$ for
some $r < 1$.  For each $r < 1$, $\overline{B}(0,r)$ is a compact
subset of $B(0,1)$, and $B(0,1) \backslash \overline{B}(0,r)$ is
connected when $n \ge 2$, and simply-connected when $n \ge 3$.  This
is because $B(0,1) \backslash \overline{B}(0,r)$ is homeomorphic to
${\bf S}^{n-1} \times (r,1)$, and ${\bf S}^{j}$ is connected when
$j \ge 1$, and simply-connected when $j \ge 2$.

	The property of being simply-connected at infinity is clearly
preserved by homeomorphisms.  Thus to get a contractable open set $U$
in ${\bf R}^n$ which is \emph{not} homeomorphic to an $n$-dimensional
ball, it suffices to choose $U$ so that it is not simply-connected at
infinity.

	If $U$ is contractable, then it is simply-connected itself in
particular.  If $U$ is simply-connected, connected at infinity, and
not simply-connected at infinity, then it means that there is a
compact set $K \subseteq U$ and loops $\gamma$ in $U$ which lie as far
towards infinity as one would like (i.e., in the complement of any given
compact subset of $U$) such that (a) $\gamma$ can be contracted to a
point in $U$, and (b) $\gamma$ cannot be contracted to a point in $U
\backslash K$.  To put it another way, these loops $\gamma$ can be
contracted to points in $U$, but in doing this one always has to pass
through at least one element of the compact set $K$.

	A mechanism for having this happen for a set $U$ contained in
${\bf R}^3$ is given by the construction of ``the Whitehead
continuum'' \cite{Whitehead}.  (See also \cite{Dm2, K}.)  Here is an
outline of the procedure.

	Start with a standard smooth ``round'' solid torus $T$ in
${\bf R}^3$.  Here $T$ should be a compact set, i.e., it should
contain its boundary.

	Next one chooses another smooth solid torus $T_1$ inside $T$.
More precisely, $T_1$ should lie in the interior of $T$.  One chooses
$T_1$ in a particular way, which can be imagined as follows.
(Pictures can be found on p68 of \cite{Dm2} and p82 of \cite{K}.)
First take a ``small'' solid torus in $T$, small enough to be
contained in a topological ball in $T$.  One can think of grabbing
hold of this small solid torus at two ends, and then stretching them
around the ``hole'' in the larger torus $T$.  One stretches them
around the two different sides of the hole in $T$.  To get $T_1$,
these two ends should hook around each other on the other side of the
hole.  In other words, one might imagine having the two ends of the
small solid torus from before, stretched around opposite sides of $T$,
and then passing one across the other, until they do not touch any
more, but are clasped together, like two hooks, or two links in a
chain.  The configuration looks \emph{locally} like two hooks or links
clasped together, but in fact one has two ends of the single solid
torus $T_1$, wrapped around the hole in $T$.

	If $T_1$ is chosen in this way, then it has the following two
basic properties.  The first is that it is \emph{homotopically}
trivial in $T$.  That is, the identity mapping on $T_1$ is homotopic
to a constant mapping through (continuous) mappings from $T_1$ into
$T$.  This follows exactly the description above; in making the
homotopy, one is allowed to stretch or move $T_1$ around as much as
one like, and one is allowed to have \emph{different parts of (images
of) $T_1$ cross each other in $T$}.  To put it a bit differently, the
mappings being deformed are \emph{not} required to be injective.

	The second property is that $T_1$ is not ``isotopically
trivial''.  This means that one cannot continuously deform $T_1$
through an isotopy of $T$ into a set which lies in a ball contained in
$T$.  In effect, this means that one cannot continuously deform $T_1$
inside $T$ in such a way that $T_1$ ends up in a ball in $T$, and so
that the deformations do \emph{not} ever cross each other (unlike the
homotopy in the previous paragraph).  If one could get $T_1$ inside a
ball in $T$, then one could continue the deformation to get an isotopy
into an arbitrarily small ball.  One would not ask for shrinking $T_1$
to a point here, because this is automatically prevented by
injectivity (independent of clasping or not).

	This explains how $T$ and $T_1$ should be chosen.  Since $T_1$
is a $3$-dimensional smooth solid torus in its own right, one can
repeat the process to get another smooth solid torus $T_2$ contained
in it, and in fact contained in the interior of $T_1$.  In other
words, since $T$ and $T_1$ are both smooth solid tori, they are
diffeomorphic to each other in particular, and this can be used to
make precise the idea of ``repeating the process''.  Specifically, if
$\phi : T \to T_1$ is such a diffeomorphism, then one can take $T_2$
to be $\phi(T_1)$.

	One then repeats the process indefinitely, getting smooth
solid tori $T_j$ for $j = 1, 2, \ldots$ such that $T_{j+1}$ is
contained in the interior of $T_j$ for each $j$, and so that $T_{j+1}$
is arranged in $T_j$ in the same way as $T_1$ is arranged in $T$.

	Now let $W$ be the intersection of all these solid tori $T_j$.
This gives a nonempty compact set in ${\bf R}^3$.  We can think of $W$
as lying inside of ${\bf S}^3$, and then take $U = {\bf S}^3
\backslash W$.  One can also rotate this around so that $U$ actually
lies in ${\bf R}^3$.  One can show that $U$ is contractable, but not
simply-connected at infinity.  See \cite{Dm2, K, Whitehead} for more
information.  (For the purposes of looking at $U$, the complement of
$W$, it can be convenient to use a modestly different description of
the construction, in which one builds $U$ up from smaller pieces in an
``increasing'' manner, analogous to the ``decreasing'' construction
for $W$ above.)

	Although $U$ is not homeomorphic to a $3$-dimensional ball in
this case, the Cartesian product of $U$ with a nonempty open interval
is homeomorphic to a $4$-dimensional ball.  This is attributed to
Arnold Shapiro in \cite{Bing3}; see also Section 10 of \cite{Bing4}
and \cite{K}.  This is analogous to the effect of stabilization
before, in Subsection \ref{Wildness and tameness phenomena}.  In
particular, one can check directly that taking the Cartesian product
with the interval gets rid of the problem that $U$ itself has with
simple-connectivity at infinity.  This is a general phenomenon, which
is relevant as well for other situations mentioned in this appendix.
A similar point came up Subsection \ref{Wildness and tameness
phenomena}.

	Beginning in dimension $4$, there are contractable open sets
in ${\bf R}^n$ which are not topological $n$-balls, and which have the
additional feature that their closures are compact manifolds with
boundary.  This last does not work in dimension $3$, and, for that
matter, the complement of the Whitehead continuum in ${\bf S}^3$
cannot be realized as the interior of a compact manifold with
boundary, whether or not this compact manifold should occur as the
closure of the set in ${\bf S}^3$.  The reason is that if such a
compact manifold did exist, its boundary would be a $2$-dimensional
surface with the homology of the $2$-sphere.  We shall say more about
this in a moment.  In this case the boundary would have to be
homeomorphic to the $2$-sphere.  This would contradict the failure of
simple-connectivity at infinity for the original space, since ${\bf
S}^2$ is simply-connected.

	The difference with $n \ge 4$ is that the boundary can be a
homology $(n-1)$-sphere (i.e., a manifold with the same homology as
${\bf S}^{n-1}$) which is not simply-connected.  The interior then
fails to be simply-connected at infinity again, and is not
homeomorphic to an $n$-ball in particular.

	For some related information and references concerning these
examples in dimensions greater than or equal to $4$, see \cite{Dm2},
including the top of p94, and the discussion on p103-104.

\subsubsection {Some positive results} 
\label{Some positive results}

	For dimensions $n \ge 4$, it is known that every contractable
topological manifold $M$ which is simply-connected at infinity is
homeomorphic to ${\bf R}^n$.  See \cite{Stallings} for $n \ge 5$, and
Corollary 1.2 on p366 of \cite{Fr} for $n = 4$.  A related reference
is \cite{McMillan-Zeeman}.  Actually, \cite{Stallings} is stated for
the piecewise-linear category; one can go from there to the
topological category via \cite{KS}.  The four-dimensional result does
not work in the smooth or piecewise-linear categories (which are
equivalent in dimension $4$), because of the existence of ``fake ${\bf
R}^4$'s (smooth manifolds homeomorphic to ${\bf R}^4$, but not
diffeomorphic to it).  Concerning the latter, see \cite{FQ} (p122 in
particular) and \cite{K} (Chapter XIV).

	These topics are also related to ``McMillan's cellularity
criterion'', in \cite{McMillan}.  A $4$-dimensional version of this is
given in \cite{Fr}, in Theorem 1.11 on p373.  We shall discuss
cellularity and this criterion further in Subsubsection
\ref{Cellularity, and the cellularity criterion}.

	Now let us look more closely at the case of compact manifolds
with boundary.  Suppose that $N$ is an $n$-dimensional compact
topological manifold with boundary.  Consider the following question:
\begin{eqnarray}
\label{contractable manifolds with spherical boundary}
    && 
\hbox{If $N$ is contractable and $\partial N$ is a topological $(n-1)$-sphere,}
                                                                     \\
    && \hbox{is $N$ homeomorphic to the closed unit ball in ${\bf R}^n$?}
                                                                  \nonumber
\end{eqnarray}
This question is actually equivalent to the Poincar\'e conjecture (in
dimension $n$, and in the topological category).  This is a well-known
fact.  The main points are the following.  If one is given a compact
$n$-dimensional topological manifold without boundary which is a
homotopy $n$-sphere, then one can get an $n$-dimensional manifold $N$
as in (\ref{contractable manifolds with spherical boundary}) from it
by cutting out a topological ball (with tame boundary).  If this
manifold $N$ is homeomorphic to the closed unit ball in ${\bf R}^n$,
then one can obtain that the original space was homeomorphic to the
standard $n$-sphere, by gluing the ball which was removed back in.
(We shall say more about this in Remark \ref{remark about gluing balls
together}.)  Conversely, given a manifold $N$ as in (\ref{contractable
manifolds with spherical boundary}), one can get a homotopy $n$-sphere
from it by gluing in a ball along the boundary of $N$.  To go from
this and the Poincar\'e conjecture to the conclusion that $N$ is
homeomorphic to a closed ball, one can use the ``generalized
Sch\"onflies theorem'', discussed later in this subsection (after
Remark \ref{remark about gluing balls together}).

	In particular, the answer to (\ref{contractable manifolds with
spherical boundary}) is known to be ``yes'' when $n \ne 3$, and the
problem is open for $n=3$.

	There are some analogous relationships between the Poincar\'e
conjecture and contractable open manifolds.  Namely, if one starts
with a compact $n$-dimensional topological manifold without boundary
$P$ which is a homotopy $n$-sphere, then one can get an
$n$-dimensional contractable open manifold by removing a point $x$
from $P$.  If $n \ge 3$, then $P \backslash \{x\}$ will also be
simply-connected at infinity, as one can check using the manifold
structure of $P$ around $x$.  If one knows that $P \backslash \{x\}$
is homeomorphic to ${\bf R}^n$, then one can deduce that $P$, which is
topologically the same as the one-point compactification of $P
\backslash \{x\}$, is homeomorphic to ${\bf S}^n$.

	However, it is not as easy to go in the other direction, from
contractable open manifolds which are simply-connected at infinity to
compact manifolds which are homotopy-equivalent to a sphere, as it is
in (\ref{contractable manifolds with spherical boundary}).  One can
take the one-point compactification of the open manifold to get a
compact space, but it is not immediately clear that this space is a
manifold.  The simple-connectivity at infinity for the open manifold
is a necessary condition for this (when $n \ge 3$), but the converse
is more complicated.  There are broader issues concerning the behavior
of open manifolds at infinity, and we shall mention some aspects of
this in Subsubsection \ref{Ends of manifolds} and Subsection
\ref{Interlude: looking at infinity, or looking near a point}.

	For the first part, about going from $P$ to an open manifold,
suppose that one is in a situation where there is a general result to
the effect that an $n$-dimensional contractable open manifold which is
simply-connected at infinity is homeomorphic to ${\bf R}^n$ for some
fixed $n$.  As above, one can use this to show that a compact
$n$-dimensional manifold $P$ (without boundary) which is
homotopy-equivalent to ${\bf S}^n$ is homeomorphic to ${\bf S}^n$.  A
complication with this type of argument is that one does not
necessarily say too much about the behavior of the homeomorphism at
the point $x$ which was removed and added back again (in the notation
before), even if one knows more about $P$ and the homeomorphism
between $P \backslash \{x\}$ and ${\bf R}^n$.  In this respect,
arguments that go through compact manifolds with boundary, as in
(\ref{contractable manifolds with spherical boundary}), can work
better; there are also some tricky aspects in this case, though, and
we shall say more about this next.

\beginremark 
\label{remark about gluing balls together}
{\rm There are some subtleties about gluing in balls in
the context of (\ref{contractable manifolds with spherical boundary})
and its correspondence with the Poincar\'e conjecture in the
\emph{smooth} category.  If one takes two copies of the closed unit
ball in ${\bf R}^n$, and glues them together using a homeomorphism
between their boundaries, then the resulting space is homeomorphic to
a standard $n$-dimensional sphere.  This is a standard observation
(which can be proved using the fact mentioned in the next paragraph),
and it works for any gluing homeomorphism.  If the gluing map is a
diffeomorphism, then the resulting space is a smooth manifold in a
natural way, but it may not be diffeomorphic to a standard sphere.
Exotic spheres can be viewed in this manner, as gluings of standard
closed balls through (tricky) diffeomorphisms along their boundaries.

	In the topological case, one can use the following fact.  Let
$B_n$ denote the closed unit ball in ${\bf R}^n$.  If $h$ is a
homeomorphism from $\partial B_n$ onto itself, then $h$ can be
extended to a homeomorphism from $B_n$ onto itself.  One can do this
by a straightforward ``radial extension''.  This method also works for
the analogous statement in the piecewise-linear category.  However, in
the smooth category, a radial extension like this is not smooth in
general at the origin in $B_n$.  An extension to a diffeomorphism may
simply not exist (radial or not).

	In any of the three categories, once one has an extension like
this, one can use it to get an equivalence between the space obtained
by gluing together the two copies of $B_n$, and the standard
$n$-dimensional sphere.  The extension unwinds the effect of the
gluing map, if the gluing map is not the standard one.  In the smooth
case, this may not be possible, and this occurs with exotic spheres.

}
\end{remark}

	Let us look some more at (\ref{contractable manifolds with
spherical boundary}), in the topological category.  If $N$ happens to
be given as a subset of ${\bf R}^n$, in addition to the conditions in
(\ref{contractable manifolds with spherical boundary}), then $N$ is
homeomorphic to the closed unit ball in ${\bf R}^n$.  This can be
derived from the ``generalized Sch\"onflies theorem'' \cite{Brown1,
Mazur, Morse}.  This result says that if one has an embedding $f$ of
${\bf S}^{n-1} \times [-1,1]$ into ${\bf R}^n$, then $f({\bf S}^{n-1}
\times \{0\})$ can be realized as the image of ${\bf S}^{n-1}$ under a
homeomorphism mapping all of ${\bf R}^n$ onto itself.  See also
Theorem 6 on p38 of \cite{Dm2}.

	Let us be a bit more precise about the way that the
generalized Sch\"onflies theorem is used here.  The first point is
that the boundary of $N$ is ``collared'' in $N$.  This means that
there is a neighborhood of $\partial N$ in $N$ which is homeomorphic
to $\partial N \times [0,1)$, and where the homeomorphism maps each
point $z \in N$ to $(z,0) \in \partial N \times \{0\}$.  The
assumption that $N$ be a topological manifold with boundary gives a
local version of this at each point in $\partial N$, and one can
derive the existence of a global collaring from a result of Brown.
See \cite{Brown2, Connelly} and Theorem 8 on p40 of \cite{Dm2}.

	On the other hand, to apply the generalized Sch\"onflies
theorem, one needs a topological $(n-1)$-sphere in ${\bf R}^n$ which
is ``bicollared'', i.e., occurs as $f({\bf S}^{n-1} \times \{0\})$ for
some embedding $f : {\bf S}^{n-1} \times [-1,1] \to {\bf R}^n$.
The boundary of $N$ is collared inside of $N$, but may not be
bi-collared inside ${\bf R}^n$.  To deal with this, one can use a
parallel copy of $\partial N$ in the interior of $N$, provided by the
collaring of $\partial N$ inside $N$.  This parallel copy is now
bi-collared in the interior of $N$, because of the collaring that we
have for $N$.  

	If $N$ lies inside ${\bf R}^n$, then this parallel copy is
also bi-collared inside ${\bf R}^n$.  One can apply the generalized
Sch\"onflies theorem, to get that the region in ${\bf R}^n$ bounded by
this parallel copy of $\partial N$, together with this copy of
$\partial N$ itself, is homeomorphic to the closed unit ball in ${\bf
R}^n$.  To get back to $N$ in its entirety, one uses the original
collaring of $\partial N$ inside $N$, to know that the missing part of
$N$ is homeomorphic to the product of $\partial N \cong {\bf S}^{n-1}$
with an interval, and to glue this to the other piece without causing
trouble.

	Thus one can get a positive answer to (\ref{contractable
manifolds with spherical boundary}) when $N$ lies inside ${\bf R}^n$,
using the generalized Sch\"onflies theorem.  This is simpler than the
solutions of the Poincar\'e conjecture, and it does not require any
restrictions on the dimension $n$.  The assumption of contractability
of $N$ is not needed for this either.  (For arbitrary manifolds, not
necessarily embedded in ${\bf R}^n$, this assumption would be
crucial.)

	In general, if $N$ is an $n$-dimensional compact topological
manifold with boundary which is contractable, then the boundary
$\partial N$ is always a homology sphere (has the same homology groups
as an $(n-1)$-dimensional sphere).  This is a well-known fact.  One
could use Theorem 9.2 on p357 of \cite{Bredon}, for instance.
Conversely, any compact $(n-1)$-dimensional topological manifold
without boundary which is a homology sphere can be realized as the
boundary of an $n$-dimensional compact topological manifold with
boundary which is contractable.  This is elementary for $n \le 3$,
where the homology spheres are all ordinary spheres, and can be filled
with balls.  For $n \ge 5$, this is given in \cite{Kervaire1}, in the
piecewise-linear category (for both the homology sphere and its
filling by a contractable manifold).  For $n \ge 6$, one can convert
this into a statement about topological manifolds, through the
Kirby--Siebenmann theory \cite{KS}.  See also the bottom of p184 of
\cite{Dm2}.  For $n=5$ in the topological category, see the corollary
on p197 of \cite{FQ}.  See also Corollary 2B on p287 of \cite{Dm2} for
$n \ge 5$ and the topological category.  (For the smooth category in
high dimensions, there are complications which come from the existence
of exotic spheres, as in the discovery of Milnor.)

	For $n=4$, see \cite{Fr, FQ}.  In particular, see Theorem 1.4'
on p367 of \cite{Fr}, and Corollary 9.3C on p146 of \cite{FQ}.  In
this case it can happen that the filling by a contractable manifold
cannot be given as a piecewise-linear manifold.  The boundary
$\partial N$ would always admit a unique piecewise-linear structure,
by well-known results about $3$-dimensional manifolds (as in
\cite{Moise}).  Concerning the possible lack of piecewise-linear
filling for a homology $3$-sphere by a contractable $4$-manifold, see
\cite{Fr, FQ, K}.

	A famous example of a homology $3$-sphere which is not
simply-connected is given by the ``Poincar\'e homology sphere''.  This
is a quotient of the standard ${\bf S}^3$ by the (finite) icosahedral
group.  See Theorem 8.10 on p353 of \cite{Bredon}.  This is a
particular example where a contractable filling exists among
topological $4$-manifolds, but not among piecewise-linear manifolds.
See \cite{Fr, FQ, K}.

	If $H$ is a $k$-dimensional compact manifold (without
boundary) which is $k$-dimensional homology sphere, and if $H$ is also
simply-connected, then $H$ is homotopy-equivalent to the standard
$k$-dimensional sphere.  This is a standard fact from topology, which
was also mentioned in Section \ref{parameterization problems}.  In
this case, the Poincar\'e conjecture in dimension $k$ would seek to
say that $H$ should be homeomorphic to ${\bf S}^k$.

	Note that if $N$ is a compact manifold with boundary, then
$\partial N$ is simply-connected if and only if the interior of $N$ is
simply-connected at infinity in the sense of Definition \ref{def of
simple connectivity at infinity}.

\subsubsection{Ends of manifolds}
\label{Ends of manifolds}

	Suppose that $M$ is an $n$-dimensional manifold without boundary
which is ``open'', i.e., not compact.  What can one say about the ``ends''
of $M$?

	In particular, when can $M$ be realized as the interior of a
compact manifold with boundary?  This would be a nice way of
``taming'' the end.

	This type of issue is clearly related to the questions
considered throughout this subsection.  It also makes sense in
general, whether or not $M$ is contractable, or one expects it to be
homeomorphic to a ball, or one expects the end to be spherical.

	Some sufficient conditions for realizing an open manifold as
the interior of a compact manifold with boundary in high dimensions
are given in \cite{BLL}.  A characterization for this is given in
\cite{Sie-thesis}.  See also \cite{Kervaire2} concerning the latter.

	For dimension $5$ (with $4$-dimensional boundaries), see
\cite{Quinn} and Section 11.9 of \cite{FS}.  For dimension $4$ (with
$3$-dimensional boundaries), see Theorem 1.12 on p373 of \cite{Fr},
and Section 11.9 in \cite{FQ}.  Concerning dimension $3$, see p216 of
\cite{FQ}.

	In all of these, the fundamental group at infinity plays an
important role.

\subsection{Interlude: looking at infinity, or looking near a point}
\label{Interlude: looking at infinity, or looking near a point}

	Let $M$ be a topological manifold of dimension $n$, and without
boundary.  Assume that $M$ is \emph{open}, i.e., not compact.

	Define $\widehat{M}$ to be the one-point compactification of
$M$, through the usual recipe.  That is, one adds to $M$ a special
point $q$, the point at infinity, and the neighborhoods of $q$ in
$\widehat{M}$ are given by sets of the form $\widehat{M} \backslash
K$, where $K$ is a compact subset of $M$.  

	Consider the following question:
\begin{equation}
\label{question about when widehat{M} is a topological manifold}
	\hbox{Under what conditions is $\widehat{M}$ a topological manifold?}
\end{equation}

	One might look at this as a kind of \emph{local} question,
about the behavior of a space at a given point, or as a question about
\emph{large-scale} behavior of $M$.  It is not hard to see that
$\widehat{M}$ will be a topological manifold exactly when $M$ looks
like (is homeomorphic to) ${\bf S}^{n-1} \times [0,1)$ outside of a
set with compact closure.  Equivalently, $\widehat{M}$ is a manifold
exactly when $M$ can be realized as the interior of a compact manifold
with boundary, where the boundary is homeomorphic to ${\bf S}^{n-1}$.
(This uses Brown's theorem about the existence of collars for
boundaries of manifolds with boundary, as in \cite{Brown2, Connelly}
and Theorem 8 on p40 of \cite{Dm2}.)

	The large-scale perspective of (\ref{question about when
widehat{M} is a topological manifold}) is somewhat close in outlook to
Subsection \ref{Contractable open sets}, especially Subsubsection
\ref{Some positive results}, while the local view is perhaps more like
the perspective in Subsection \ref{Wildness and tameness phenomena}.
Concerning the latter, one might think of local taming properties of
embedded sets in terms of existence of ``normal bundles'' for the
embedded sets.  Similarly, one can think of (\ref{question about when
widehat{M} is a topological manifold}) as asking about the existence
of a normal bundle for $\widehat{M}$ at the point $q$.  In this
regard, one might compare with the discussion in Section 9.3 in
\cite{FQ}, especially Theorem 9.3A and Corollary 9.3B.  

	For the record, let us note that a necessary condition for
$\widehat{M}$ to be a manifold is that
\begin{equation}
\label{M is simply-connected at infinity}
	\hbox{$M$ is simply-connected at infinity,}
\end{equation}
at least if $n \ge 3$.  In ``local'' language, we can reformulate this
condition as follows: for every neighborhood $U$ of $q$ in
$\widehat{M}$, there is a neighborhood $V$ of $q$ such that $V
\subseteq U$,
\begin{equation}
\label{local connectedness of widehat{M} backslash {q}}
	\qquad \hbox{every pair of points } x, y \in V \backslash \{q\}
         \hbox{ lies in a connected set in } U \backslash \{q\},
\end{equation}
and
\begin{equation}
\label{local vanishing of pi_1 in widehat{M} backslash {q}}
	\qquad \hbox{every loop $\gamma$ in } V \backslash \{q\}
          \hbox{ can be contracted to a point in } U \backslash \{q\}.
\end{equation}
This is similar to the localized fundamental group conditions
mentioned in Subsection \ref{Wildness and tameness phenomena}.

	Let us now think of $\widehat{M}$ as being any topological
space, and not necessarily the one-point compactification of an open
manifold.  For a given point $q \in \widehat{M}$, one can still ask
whether $\widehat{M}$ is an $n$-dimensional manifold at $q$, i.e., if
there is a neighborhood of $q$ in $\widehat{M}$ which is homeomorphic
to an open ball in ${\bf R}^n$.  The necessary condition in the
preceding paragraph still applies (when $n \ge 3$), concerning local
simple-connectivity of $\widehat{M} \backslash \{q\}$ near $q$ (as in
(\ref{local connectedness of widehat{M} backslash {q}}) and
(\ref{local vanishing of pi_1 in widehat{M} backslash {q}})).

	For the rest of this subsection, let us assume that $n \ge 3$.
Note that there are special results for detecting manifold behavior in
a space of dimension $1$ or $2$.  This is reviewed in the introduction
of \cite{Fe4}.

	As a special case, imagine now that $\widehat{M}$ is a finite
polyhedron of dimension $n$.  Let $L$ denote the codimension-$1$ link
of $q$ in $\widehat{M}$, as discussed in Section \ref{parameterization
problems}.  Thus $L$ is an $(n-1)$-dimensional finite polyhedron, and
$\widehat{M}$ looks locally at $q$ like a cone over $L$, as in Section
\ref{parameterization problems}.  (As in Section \ref{parameterization
problems}, $L$ is determined up to piecewise-linear equivalence, but
not as a polyhedron.)

	In order for $\widehat{M}$ to be an $n$-dimensional
topological manifold in a neighborhood of $q$, the link $L$ should be
fairly close to a standard $(n-1)$-sphere.  In particular, it is not
hard to see that $L$ should be homotopy-equivalent to ${\bf S}^{n-1}$.
This implies that $L$ should be connected and simply-connected, under
our assumption that $n$ is at least $3$.

	In fact, in the case where $\widehat{M}$ is a finite
polyhedron, the connectedness and simple-connectedness of the link $L$
around $q$ are equivalent to the local connectivity and
simple-connectivity conditions for $\widehat{M} \backslash \{q\}$ near
$q$ indicated above, with (\ref{local connectedness of widehat{M}
backslash {q}}) and (\ref{local vanishing of pi_1 in widehat{M}
backslash {q}}).  This is not hard to see, and it is also rather nice.
To put it a bit differently, imagine that one starts with the class of
finite polyhedra, and then tries to go to more general contexts of
topological spaces.  The local connectivity and simple-connectivity
conditions for $\widehat{M} \backslash \{q\}$ at $q$ as described
above provide a way to capture the information in the connectedness
and simple-connectedness of the codimension-$1$ link at $q$ in the
case where $\widehat{M}$ is a polyhedron, in a manner that makes sense
for arbitrary topological spaces, without special structure as one has
for finite polyhedra.

	These local connectedness and simple-connectedness conditions
for $\widehat{M} \backslash \{q\}$ at $q$ should be compared with
local connectedness and simple-connectedness conditions for
$\widehat{M}$ itself.  If $\widehat{M}$ is a polyhedron, then any
point $q$ in $\widehat{M}$ automatically has the feature that there
are arbitrarily small neighborhoods $U$ of $q$ in $\widehat{M}$ which
can be contracted to $q$ while staying near $q$, and, in fact, while
staying in $U$ itself.  This is because $\widehat{M}$ looks locally
like a cone at $q$.

	In the next subsection we shall look at another case of this kind
of ``local manifold'' question.

\subsection{Decomposition spaces, 1}
\label{Decomposition spaces, 1}

	Let $n$ be a positive integer, and let $K$ be a nonempty
compact subset of ${\bf R}^n$.  One could also consider general
manifolds instead of ${\bf R}^n$ here, but we shall generally stick to
Euclidean spaces for simplicity.  The main ideas come up in this case
anyway.

	Imagine shrinking $K$ to a single point, while leaving the
rest of ${\bf R}^n$ alone, and looking at the topological space that
results.  This can be defined more formally as follows.  Let us write
${\bf R}^n / K$ for the set which consists of the points in ${\bf
R}^n$ which do not lie in $K$, together with a single point which
corresponds to $K$ itself.  In other words, this is where we shrink
$K$ to a single point.  This set can be given a topology in a standard
way, so that a subset $U$ of ${\bf R}^n / K$ is open if and only if
its inverse image back in ${\bf R}^n$ is open.  Here ``inverse image''
uses the automatic quotient mapping ${\bf R}^n$ to ${\bf R}^n / K$.
(In concrete terms, the inverse image of $U$ in ${\bf R}^n$ means the
set of points in ${\bf R}^n$ which correspond to elements of $U$,
where one includes all points in $K$ if the element of ${\bf R}^n / K$
associated to $K$ lies in $U$.)

	This type of quotient ${\bf R}^n / K$ is a special case of a
``decomposition space''.  We shall discuss the general situation
further in Subsection \ref{Decomposition spaces, 2}, but this special
case already includes a lot of interesting examples and phenomena.

	Now let us consider the following question:
\begin{equation}
\label{when is {bf R}^n / K a topological manifold?}
  \hbox{Given $K$ as above, when is ${\bf R}^n / K$ a topological manifold?}
\end{equation}

	This is really a special case of the situation in Subsection
\ref{Interlude: looking at infinity, or looking near a point}.  For
this it is better to use ${\bf S}^n$ instead of ${\bf R}^n$, so that
${\bf S}^n / K$ --- defined in the same manner as above --- is
equivalent to the one-point compactification of ${\bf S}^n \backslash
K$.

	Let us consider some basic examples.  If $K$ consists of only
a single point, then ${\bf R}^n / K$ is automatically the same as
${\bf R}^n$ itself, and there is nothing to do.  If $K$ is a finite
set with more than one element, then it is easy to see that ${\bf R}^n
/ K$ is not a manifold.  If we let $q$ denote the point in ${\bf R}^n
/ K$ which corresponds to $K$, then $({\bf R}^n / K) \backslash \{q\}$
does not enjoy the local connectedness property that it should if
${\bf R}^n / K$ were a manifold at $q$, as in (\ref{local
connectedness of widehat{M} backslash {q}}) in Subsection
\ref{Interlude: looking at infinity, or looking near a point}.  More
precisely, this local connectedness property for the complement of
$\{q\}$ would be necessary only when $n \ge 2$.  When $n = 1$, one
does not have to have this local connectedness condition, but then
$({\bf R}^n / K) \backslash \{q\}$ would have too many local
components near $q$ for ${\bf R}^n / K$ to be a manifold at $q$.
(That is, there would be more than $2$ such local components.)

	Now suppose that $K$ is a straight line segment in ${\bf
R}^n$.  In this event, ${\bf R}^n \backslash K$ is homeomorphic to
${\bf R}^n$ again.  This is not hard to check.  This would also work
if $K$ were a standard rectangular cell of higher dimension in ${\bf
R}^n$.

	More generally, this works if $K$ is a \emph{tame cell} in
${\bf R}^n$, meaning the image of a standard rectangular cell under a
homeomorphism of ${\bf R}^n$ onto itself.  This follows automatically
from the case of standard rectangular cells.

	However, if one merely assumes that $K$ is homeomorphic to a
standard rectangular cell, then it is not necessarily true that ${\bf
R}^n / K$ is a manifold!  This is another aspect of wild embeddings,
from Subsection \ref{Wildness and tameness phenomena}.  We shall say
more about this as the subsection goes on.  A concrete example is
given by taking $K$ to be a copy of the Fox--Artin wild arc in ${\bf
R}^3$.  (Compare with \cite{Fe4}.)  

	Note that we are not saying that ${\bf R}^n / K$ is
\emph{always} not a manifold when $K$ is wildly embedded.  The
converse is true, that $K$ must be wildly embedded when ${\bf R}^n /
K$ is not a manifold (and $K$ is a topological cell).  This is just a
rephrasal of the remark above, that ${\bf R}^n / K$ is a manifold when
$K$ is a tamely embedded cell.

	Here is a slightly more foolish example, which one might view
as a generalization of the earlier comments about the case where $K$
is a finite set with more than a single point.  Imagine now that $K$
is a copy of the $j$-dimensional sphere ${\bf S}^j$, $1 \le j \le
n-1$.  For this let us use a standard, smooth, round sphere; it is not
a matter of wildness that we want to consider.

	In this case ${\bf R}^n / K$ is \emph{never} a topological
manifold.  If $j = n-1$, then ${\bf R}^n / K$ is homeomorphic to the
union of ${\bf R}^n$ and an $n$-sphere, with the two meeting at a
single point.  This point is the one that corresponds to $K$ in ${\bf
R}^n / K$.  Let us denote this point by $q$ again, as above.  In this
case $({\bf R}^n / K) \backslash \{q\}$ does not have the right
local-connectedness property at $q$ in order for ${\bf R}^n / K$ to be
a manifold, as in (\ref{local connectedness of widehat{M} backslash
{q}}) in Subsection \ref{Interlude: looking at infinity, or looking
near a point}.

	If $j = n-2$, then one runs into trouble with local
simple-connectivity of $({\bf R}^n / K) \backslash \{q\}$ at $q$, as
in (\ref{local vanishing of pi_1 in widehat{M} backslash {q}}) in
Subsection \ref{Interlude: looking at infinity, or looking near a
point}.  For this one might think about the special case where $n =3$,
so that $K$ is a standard circle in ${\bf R}^3$.  It is easy to take
small loops in ${\bf R}^3 \backslash K$, lying close to $K$, which are
nonetheless linked with $K$.  These loops then project down into
$({\bf R}^3 / K) \backslash \{q\}$, where they can be as close to the
point $q$ as one likes, but they are never contractable in $({\bf R}^3
/ K) \backslash \{q\}$ at all, let alone in small neighborhoods of $q$
(as in (\ref{local vanishing of pi_1 in widehat{M} backslash {q}}) in
Subsection \ref{Interlude: looking at infinity, or looking near a
point}).  This is the same as saying that these loops are not
contractable inside of ${\bf R}^n \backslash K$, which is equivalent
to $({\bf R}^n / K) \backslash \{q\}$.

	When $j < n-2$, then one has similar obstructions to ${\bf
R}^n / K$ being a manifold, but in terms of the failure of
higher-dimensional forms of local connectedness of $({\bf R}^n / K)
\backslash \{q\}$ (using homology or homotopy).  This is analogous to
the cases already described, when $j = n-1$ or $n-2$.  We shall say
more about this soon, but for the moment let us go on to some other
matters.

	For this example, where $K$ is taken to be a standard
$j$-dimensional sphere, note that ${\bf R}^n / K$ itself is locally
contractable at $q$.  This is as opposed to connectedness properties
of $({\bf R}^n / K) \backslash \{q\}$, and it is analogous to what
happens in the case of finite polyhedra.  Specifically, for finite
polyhedra one always has local contractability, but the behavior near
a given point of \emph{punctured} neighborhoods around that point are
another matter.  The latter is connected to the behavior of the
codimension-$1$ link of the polyhedron around the given point, as in
Subsection \ref{Interlude: looking at infinity, or looking near a
point}.

	In the present case, where we have ${\bf R}^n / K$ with $K$ a
standard round $j$-dimensional sphere, one can see the local
contractability of ${\bf R}^n / K$ at the point $q$ (corresponding to
$K$) as follows.  In ${\bf R}^n$, one can take a tubular neighborhood
of $K$, which is homeomorphic to the Cartesian product of the
$j$-sphere $K$ and an $(n-j)$-dimensional ball.  This neighborhood can
be contracted onto $K$ in a simple way, and this leads to the local
contractability of ${\bf R}^n / K$ at $q$.

	Now let us consider the case of the \emph{Whitehead
continuum}, from Subsection \ref{Contractable open sets}.  We should
not really say \emph{the} Whitehead continuum here, as there is some
flexibility in the construction, which can lead to the resulting set
$W$ not being pinned down completely.  This ambiguity will not really
cause trouble for us here, and we can work with any compact set $W$ in
${\bf R}^3$ which is obtained as in the procedure described in
Subsection \ref{Contractable open sets}.

	The set $W$ has the feature of being \emph{cell-like}, as in the
following definition.

\begindefinition [Cell-like sets]
\label{cell-likeness}
A compact set $K$ in ${\bf R}^n$ is said to be \emph{cell-like} if $K$
can be contracted to a point inside of any neighborhood $U$ of itself
in ${\bf R}^n$.
\end{definition}

	Compare with \cite{Dm2}, especially p120.  That the Whitehead
continuum $W$ is cell-like is not hard to see from the construction of
$W$, as the intersection of a decreasing sequence of solid tori with
certain properties.  Specifically, for this the key point is that the
$\ell$th solid torus can be contracted to a point inside the previous
one.

	If $K$ is a topological cell, then $K$ is contractable to a
point inside of itself, without using the extra bit of room provided
by a small neighborhood of itself.  This is also independent of the
way that $K$ might be embedded into some ${\bf R}^n$, i.e., wildly or
tamely.  In this respect, $W$ is like a topological cell (and hence the
name ``cell-like'' for the property in Definition \ref{cell-likeness}).

	For $W$, it is \emph{not} true that ${\bf R}^3 / W$ is a
topological manifold.  If we let $q$ denote the point in ${\bf R}^3 /
W$ corresponding to $W$, then $({\bf R}^3 / W) \backslash \{q\}$ is
\emph{not} locally simply-connected at $q$ (in the sense of the
condition in Subsection \ref{Interlude: looking at infinity, or
looking near a point}, around (\ref{local vanishing of pi_1 in
widehat{M} backslash {q}})).  In concrete terms, this means that
there are loops in ${\bf R}^3 \backslash W$ (which is the same
as $({\bf R}^3 / W) \backslash \{q\}$) which lie as close to $W$
as one likes (in their entirety), but which cannot be contracted
to a point in ${\bf R}^3 \backslash W$ while remaining reasonably
close to $W$.

	These loops can be described concretely, as meridians in the
solid tori whose intersection gives $W$.  The loop from the solid
torus $T_j$ can be filled with a disk inside $T_j$, but not without
crossing the smaller torus $T_{j+1}$, or any of its successors.  This
comes back to the way that each $T_{\ell+1}$ is ``clasped'' inside of
$T_{\ell}$.  See \cite{Dm2, K} for more information (including
Proposition 9 on p76 of \cite{Dm2}).

	In any event, the failure of the local simple-connectivity of
$({\bf R}^3 / W) \backslash \{q\}$ at $q$ is equivalent to ${\bf S}^3
\backslash W$ not being simply-connected at infinity, as in Subsection
\ref{Contractable open sets}.  This also follows the discussion in
Subsection \ref{Interlude: looking at infinity, or looking near a
point}, and the comment just after (\ref{when is {bf R}^n / K a
topological manifold?}).

	This case is quite different from the one of embedding round
spheres in ${\bf R}^n$, as discussed before.  More precisely, let us
compare the situation with $W$ and the example before where $K$ is a
standard circle inside of ${\bf R}^3$.  For the latter, there are
loops in ${\bf R}^3 \backslash K$ which lie as close to $K$ as one
wants, and which are not contractable to a point in ${\bf R}^3
\backslash K$ at all, let alone in a neighborhood of $K$.  For $W$,
one has that ${\bf S}^3 \backslash W$ is contractable (as mentioned in
Subsection \ref{Contractable open sets}), and this implies that ${\bf
R}^3 \backslash W$ is simply-connected.  (This is a straightforward
exercise.)  Thus these loops near $W$ can be contracted to a point in
${\bf R}^3 \backslash W$, if one allows oneself to go away from $W$
for the contraction.

	Here is another aspect of this.  Although one has these loops
in ${\bf R}^3 \backslash W$ which lie near $W$ but cannot be
contracted to a point in ${\bf R}^3 \backslash W$ while staying near
$W$, these loops can be made \emph{homologically} trivial in ${\bf
R}^3 \backslash W$ while staying near $W$.  That is, one can fill the
loops with surfaces inside ${\bf R}^3 \backslash W$ while staying
close to $W$, if one allows the surfaces to have handles (rather than
simply being a disk, as in the case of homotopic triviality).  This is
something that one can easily see from the pictures (as in \cite{Dm2,
K}).  The basic idea is that one can fill the loops with disks, where
the disks stay close to $W$, but also pass through $W$ (and so are not
in ${\bf R}^3 \backslash W$).  However, one can avoid the intersection
with $W$ by cutting out a couple of small holes in the disk, and
attaching a handle to them which goes along the boundary of the solid
torus in the next generation of the construction.  Then $W$ will stay
inside this next solid torus, throughout the rest of the construction,
and this surface gives a way of filling the loop without intersecting
$W$ (or being forced to go far away from it).

	This kind of filling by surfaces does \emph{not} work in the
case where we take $K$ to be a standard circle in ${\bf R}^3$.  In
this situation, we have loops in ${\bf R}^3 \backslash K$ which lie
close to $K$, and which are linked \emph{homologically} with the
circle $K$.  In other words, the \emph{linking number} of the loop
with $K$ is nonzero, and this linking number is a \emph{homological}
invariant which would vanish if the loop could be filled with a
surface without intersecting $K$.  (For more about ``linking
numbers'', see \cite{BT, Bredon, Fl, Sp}.)

	In this respect, the case of Whitehead continua is much more
tame than that of embedded circles and spheres of other dimensions,
even if it is still singular.  We shall encounter other versions of
this, now and further in this appendix.

	Here is another feature of $W$, which distinguishes it from
ordinary circles in ${\bf R}^3$ (or spheres in ${\bf R}^n$ more
generally).  Let us think of $W$ now as lying in ${\bf R}^4$ rather
than ${\bf R}^3$, through the inclusion of ${\bf R}^3$ in ${\bf R}^4$
by taking the fourth coordinate to be $0$.  

	For ${\bf R}^4$, we have that ${\bf R}^4 / W$ \emph{is} a
topological manifold (homeomorphic to ${\bf R}^4$).  The basic point
behind this is the following.  In the realization of $W$ as the
intersection of a decreasing sequence of solid tori in ${\bf R}^3$,
the $\ell$th solid torus was always ``clasped'' in the previous one
(as in Subsection \ref{Contractable open sets}, and \cite{Dm2, K}).
In ${\bf R}^4$, the extra dimension provides a lot of extra room, in
such a way that this ``clasping'' is not really present any more.  If
$T'$ is a solid torus which is embedded and clasped inside of another
solid torus $T$ in ${\bf R}^3$, one can ``unclasp'' $T'$ in ${\bf
R}^4$ by lifting one end up, bringing it around the hole in $T$, and
leaving the other end alone.  This is a standard observation, and it
is analogous to the way that knots in ${\bf R}^3$ become unknotted in
${\bf R}^4$.

	In other words, this procedure gives a way to make a
deformation of ${\bf R}^4$, in which the solid torus $T'$ is mapped to
a set of small diameter, while not moving points some distance away at
all.  By contrast, back in ${\bf R}^3$, it is not possible to make an
isotopy which shrinks $T'$ to a set of small diameter, while leaving
the points in the complement of the larger solid torus fixed.  This is
exactly because of the way that $T'$ is ``clasped'' in $T$, so that it
cannot be ``unclasped'' by an isotopy in $T$.  When one has the extra
dimension in ${\bf R}^4$, one can ``undo'' the clasping, by lifting
one end up and moving it around, as indicated above.

	Once one has this kind of ``shrinking'' in ${\bf R}^4$, one
can use this to show that ${\bf R}^4 / W$ is homeomorphic to ${\bf
R}^4$.  One can do this directly, using shrinking homeomorphisms like
this, and combinations of them, to make a mapping from ${\bf R}^4$ to
itself which shrinks $W$ to a point while remaining injective (and
continuous) everywhere else.  One puts homeomorphisms like this on top
of each other, and deeper and deeper in the construction of $W$, until
$W$ itself is shrunk all the way to a point.  The various solid tori
$T_j$ in the construction, of which $W$ is the intersection, are made
smaller and smaller in this process.  The trick is to do this without
shrinking everything, so that the mapping that results remains a
homeomorphism on the complement of $W$.

	This idea of shrinking can be given a general form, and is
discussed in detail in \cite{Dm2}.  See also \cite{E, K}.

	By contrast, let us consider the case of a circle $K$ in ${\bf
R}^3$.  If one views $K$ as a subset of ${\bf R}^4$ in the same way,
then ${\bf R}^4 / K$ is still \emph{not} a topological manifold.  This
follows from our earlier discussion about circles and spheres of
higher dimensions inside of ${\bf R}^n$ in general.  One also does not
get a manifold by replacing ${\bf R}^4$ with ${\bf R}^m$ for larger
$m$'s.

	Notice, however, that there is a kind of ``improvement'' that
occurs in adding dimensions in this way.  If $K$ is a circle in ${\bf
R}^3$, and if $q$ denotes the point in ${\bf R}^3 / K$ which
corresponds to $K$, then $({\bf R}^3 / K) \backslash \{q\}$ is not
locally simply-connected at $q$.  For that matter, $({\bf R}^3 / K)
\backslash \{q\} \cong {\bf R}^3 \backslash K$ is not simply-connected
at all.  When one considers $K$ as a subset of ${\bf R}^4$, and asks
analogous questions for ${\bf R}^4 / K$ (or ${\bf R}^4 \backslash K$),
then there is no longer any trouble with \emph{simple}-connectivity.
The basic underlying problem continues, though, in the form of
\emph{$2$-dimensional} connectivity.  This is not hard to see.

	Similarly, if one views $K$ as a subset of ${\bf R}^n$ for
larger $n$, then the trouble with connectivity in lower dimensions
goes away, but $(n-2)$-dimensional connectivity still does not work.

	With the Whitehead continuum we are more fortunate.  The
problem with local simple-connectivity goes away when we proceed from
${\bf R}^3$ to ${\bf R}^4$, but difficulties with higher-dimensional
connectivity do not then arise in their place.  One should not be too
surprised about this, since the Whitehead continuum is cell-like,
while circles or spheres of higher dimension are not at all cell-like.
In other words, with circles or spheres (and their complements in
${\bf R}^n$), there is some clear and simple nontrivial topology
around, while the Whitehead continuum is much closer to something like
a standard cell, which causes less trouble.  (One can look at this
more precisely, but we shall not pursue this here.)

\subsubsection{Cellularity, and the cellularity criterion}
\label{Cellularity, and the cellularity criterion}

	Now let us look at some general notions and results,
concerning the possibility that ${\bf R}^n / K$ be a topological
manifold (and, in fact, homeomorphic to ${\bf R}^n$).

\begindefinition [Cellularity]
\label{definition of cellularity}
A compact set $K$ in ${\bf R}^n$ (or, more generally, an
$n$-dimensional topological manifold) is said to be \emph{cellular} if
it can be realized as the intersection of a countable family of sets
$B_i$, where each $B_i$ is a topological $n$-cell (or, equivalently,
homeomorphic to the closed unit ball in ${\bf R}^n$), and if each
$B_{i+1}$ is contained in the interior of the preceding $B_i$.
\end{definition}

	Compare with \cite{Dm2}, especially p35, \cite{E}, and p44 of
\cite{Ru-book}.  Alternatively, a compact set $K$ is cellular if and
only if any neighborhood of $K$ contains an open set which contains
$K$ and is homeomorphic to the standard $n$-dimensional ball.

\begintheorem
\label{R^n/K a topological manifold if and only if K is cellular}
Let $K$ be a compact subset of ${\bf R}^n$.  Then ${\bf R}^n / K$ is
a topological manifold if and only if $K$ is cellular in ${\bf R}^n$.
In this case, ${\bf R}^n / K$ is homeomorphic to ${\bf R}^n$.
\end{theorem}

	See Exercise 7 on p41 of \cite{Dm2} for the first assertion,
and Proposition 2 on p36 of \cite{Dm2} for the second one.
(Concerning the latter, see Section 5 in \cite{Dm2} too.  Note that
some of the notation in Exercise 7 on p41 in \cite{Dm2} is explained
in the statement of Proposition 2 on p36 of \cite{Dm2}.)  See also
\cite{E}, especially the theorem on p114, and p44ff of \cite{Ru-book}.

	For the record, let us mention the following.

\beginproposition
\label{cellularity and cell-likeness}
Let $K$ be a compact subset of ${\bf R}^n$.  If $K$ is cellular, then
$K$ is cell-like.  Conversely, if $n$ is equal to $1$ or $2$, then $K$
is cellular if it is cell-like.
\end{proposition}

	The fact that cellularity implies cell-likeness follows easily
from the definitions.  When $n=1$, the converse is very simple, since
connectedness implies that a set is an interval, and hence cellular.
In ${\bf R}^2$, the argument uses special features of plane topology.
See Corollary 4C on p122 of \cite{Dm2}.

	In higher dimensions, cell-like sets need not be cellular.
Examples are given by Whitehead continua, and some wild embeddings of
cells.  However, there is an exact characterization of cellular sets
among cell-like sets, which is the following.  Basically, the point is
to include the same kind of localized simple-connectivity of ${\bf
R}^n \backslash K$ around $K$ as discussed before.

\begintheorem
\label{cellularity equivalence with cell-likeness and local simple-conn.}
Let $K$ be a compact set in ${\bf R}^n$, with $n \ge 3$.  Then $K$ is
cellular inside of ${\bf R}^n$ if and only if {\rm (a)} it is
cell-like, and {\rm (b)} for every open neighborhood $U$ of $K$ in
${\bf R}^n$ there is another open neighborhood $V$ of $K$ so that
every continuous mapping from ${\bf S}^1$ into $V \backslash K$ can be
contracted to a point inside of $U \backslash K$.
\end{theorem}

	This characterization of cellularity is stated in Theorem 5 on
p145 of \cite{Dm2}.  This uses also the definition of the
\emph{cellularity criterion} given on p143 of \cite{Dm2}.  When $n \ge
4$, this result works for subsets of general $n$-dimensional
topological manifolds, and not just ${\bf R}^n$.  When $n=3$, there is
trouble with the general case of manifolds, related to the
$3$-dimensional Poincar\'e conjecture being unsettled; if the
cellularity criterion holds for general manifolds, then the
$3$-dimensional Poincar\'e conjecture would follow, as discussed on
p145 of \cite{Dm2}.  See Theorem 1.11 on p373 of \cite{Fr} concerning
the $4$-dimensional case, and \cite{McMillan} and Section 4.8 of
\cite{Ru-book} for dimensions $5$ and higher.

\begincorollary
\label{cell-like sets are cellular in one higher dimension}
Let $K$ be a compact subset of ${\bf R}^n$, $n \ge 3$.  If $K$ is
cell-like in ${\bf R}^n$, then $K \times \{0\}$ is cellular
in ${\bf R}^{n+1}$.
\end{corollary}

	See Corollary 5A on p145 of \cite{Dm2}.  

	Corollary \ref{cell-like sets are cellular in one higher
dimension} is analogous to the fact that wild embeddings into ${\bf
R}^n$ can become tame when one passes from ${\bf R}^n$ as the ambient
space to an ${\bf R}^m$ with $m > n$.  This was discussed briefly in
Subsection \ref{Wildness and tameness phenomena}, towards the end.
Similarly, Theorem \ref{cellularity equivalence with cell-likeness and
local simple-conn.} is analogous to taming theorems for embeddings
mentioned in Subsection \ref{Wildness and tameness phenomena}.

	Theorem \ref{cellularity equivalence with cell-likeness and
local simple-conn.} and Corollary \ref{cell-like sets are cellular in
one higher dimension} are also close to some of the matters described
in Subsection \ref{Contractable open sets}.

	The main point behind the derivation of Corollary
\ref{cell-like sets are cellular in one higher dimension} from Theorem
\ref{cellularity equivalence with cell-likeness and local
simple-conn.} is that by passing to a Euclidean space of one higher
dimension, potential trouble with local simple-connectedness of the
complement of $K$ goes away.  This fits with basic examples, and the
Whitehead continuum in particular. 

	We saw before that when $K$ is a sphere, the passage from
${\bf R}^n$ to ${\bf R}^{n+1}$ can get rid of the trouble with local
simple-connectivity of the complement of $K$, but that problems remain
with higher-dimensional connectivity of the complement.  For cell-like
sets, unlike general sets, it is only the localized $1$-dimensional
connectivity of the complement which is needed to get cellularity.
This is shown by Theorem \ref{cellularity equivalence with
cell-likeness and local simple-conn.}.

	In this regard, let us also notice the following simple
converse to Corollary \ref{cell-like sets are cellular in one higher
dimension}.

\beginlemma
\label{K times {0} cell-like in {bf R}^{n+1} implies K cell-like in {bf R}^n}
Suppose that $K$ is a compact subset of ${\bf R}^n$.  If $K \times
\{0\}$ is cellular in ${\bf R}^{n+1}$, then $K$ is cell-like in ${\bf
R}^n$.
\end{lemma}

	Indeed, if $K \times \{0\}$ is cellular in ${\bf R}^{n+1}$,
then it is also cell-like in ${\bf R}^{n+1}$, as in Proposition
\ref{cellularity and cell-likeness}.  It is easy to check that
cell-likeness for $K \times \{0\}$ in ${\bf R}^{n+1}$ implies
cell-likeness for $K$ inside ${\bf R}^n$, just by the definitions.
(Thus cell-likeness, unlike cellularity, is not made more feasible by
the extra room of extra dimensions.)  This implies Lemma \ref{K times
{0} cell-like in {bf R}^{n+1} implies K cell-like in {bf R}^n}.

	For concrete examples of cell-like sets, often the cellularity
in higher-dimensional spaces, as in Corollary \ref{cell-like sets are
cellular in one higher dimension}, can be seen in fairly direct and
simple terms.  The room from the extra dimensions makes it easy to
move pieces of the set apart, without the claspings, knottings, etc.,
which occurred originally.  Some aspects of this came up earlier,
concerning Whitehead continua.

	Before leaving this subsection, let us observe that the
localized simple-connectivity conditions that are used here are a bit
different from those employed in the context of taming theorems, as in
Subsection \ref{Wildness and tameness phenomena}.  To make this
precise, let $K$ be a compact subset of some ${\bf R}^n$.  The
conditions that come up in the present subsection involve the behavior
of ${\bf R}^n \backslash K$, localized around $K$ (the whole of $K$).
That is, one looks at the behavior of ${\bf R}^n \backslash K$ within
arbitrarily-small neighborhoods of $K$ in ${\bf R}^n$.  In the context
of Subsection \ref{Wildness and tameness phenomena}, one would look at
the behavior of ${\bf R}^n \backslash K$ near individual points in
$K$.

	To put it another way, here one seeks to contract loops in
${\bf R}^n \backslash K$ that are close to $K$ to points, while
staying close to $K$.  In the context of Subsection \ref{Wildness and
tameness phenomena}, one looks at \emph{small} loops in ${\bf R}^n
\backslash K$ near $K$, and tries to contract them to points in the
complement while staying in \emph{small balls}, and not just staying
near $K$.

\subsection{Manifold factors}
\label{Manifold factors}

	Let $W$ be a Whitehead continuum, constructed through a
decreasing sequence of solid tori in ${\bf R}^3$, as in Subsection
\ref{Contractable open sets}.

\begintheorem
\label{(R^3/W) times R is homeomorphic to R^4}
If ${\bf R}^3 / W$ is defined as in Subsection \ref{Decomposition
spaces, 1}, then $({\bf R}^3 / W) \times {\bf R}$ is homeomorphic to
${\bf R}^4$.
\end{theorem}

	In particular, $({\bf R}^3 / W) \times {\bf R}$ is a
topological manifold, even though ${\bf R}^3 / W$ itself is not.  Thus
$({\bf R}^3 / W) \times {\bf R}$ is a \emph{manifold factor}.

	The fact that $({\bf R}^3 / W) \times {\bf R}$ is homeomorphic
to ${\bf R}^4$ is given as Corollary 3B on p84 of \cite{Dm2}.  See
also \cite{Andrews-Rubin, K}.

	Note that the existence of a homeomorphism from $({\bf R}^3 /
W) \times {\bf R}$ onto ${\bf R}^4$ is \emph{not} the same as the
observation mentioned in Subsection \ref{Decomposition spaces, 1},
that ${\bf R}^4 / (W \times \{0\})$ is homeomorphic to ${\bf R}^4$.
In considering $({\bf R}^3 / W) \times {\bf R}$, one is in effect
taking ${\bf R}^4$, and then shrinking each copy $W \times \{u\}$ of
$W$ to a point, where $u$ runs through all real numbers.  For ${\bf
R}^4 / (W \times \{0\})$, one shrinks only a single copy of $W$ to a
point.

	Although the construction is more complicated for $({\bf R}^3
/ W) \times {\bf R}$ than for ${\bf R}^4 / (W \times \{0\})$, there
are some common aspects.  As before, one of the main points is that
the solid tori in ${\bf R}^3$ which are ``clasped'' (inside of other
solid tori) become unclasped in ${\bf R}^4$.  With the extra dimension
in ${\bf R}^4$, one can pick up one end of one of these tori, bring it
around, and then lay it down again, so that the clasping is undone.
For the present situation with $({\bf R}^3 / W) \times {\bf R}$, one
performs this kind of action for all of the copies $W \times \{u\}$ of
$W$ at once, $u \in {\bf R}$, rather than just a single copy.
(Compare also with Subsection \ref{Decomposition spaces, 2}, and the
general notion of decomposition spaces mentioned there.)

	In Subsection \ref{Contractable open sets}, it was mentioned
that ${\bf S}^3 \backslash W$ is a contractable open set which is
\emph{not} homeomorphic to a $3$-ball (because it is not
simply-connected at infinity), and that $({\bf S}^3 \backslash W)
\times {\bf R}$ \emph{is} homeomorphic to a $4$-dimensional open ball.
(See \cite{Bing3, Bing4, K}.)  This result is similar in some ways to
Theorem \ref{(R^3/W) times R is homeomorphic to R^4}, but the
conclusions are not quite the same either.

	In this vein, let us make the following observation.  As
usual, denote by $q$ the (singular) point in ${\bf R}^3 / W$ that
corresponds to $W$.  Let us write $L$ for the subset of $({\bf R}^3 /
W) \times {\bf R}$ given by $q \times {\bf R}$.  Thus $L$ is
homeomorphic to a line.

	Using a homeomorphism from $({\bf R}^3 / W) \times {\bf R}$ to
${\bf R}^4$, one gets an embedding of $L$ into ${\bf R}^4$.  It is not
hard to see that any such embedding of $L$ into ${\bf R}^4$ has to be
wild.  Just as ${\bf R}^3 \backslash W$ is not locally
simply-connected near $W$, if $\mathcal{L}$ denotes the image of $L$
in ${\bf R}^4$ by an embedding as above, then ${\bf R}^4 \backslash
\mathcal{L}$ is not locally simply-connected near $\mathcal{L}$.
(Note that ${\bf R}^4 \backslash \mathcal{L}$ is homeomorphic to
$({\bf R}^3 \backslash W) \times {\bf R}$, by construction.)  This
ensures that $\mathcal{L}$ is wild in ${\bf R}^4$, no matter what
homeomorphism from $({\bf R}^3 / W) \times {\bf R}$ onto ${\bf R}^4$
one might use, since ordinary straight lines in ${\bf R}^4$ do not
behave in this way.

	One can also make local versions of this argument, to show
that ${\bf L}$ is locally wild in the same manner.

	If one were to want to pass from a homeomorphism from $({\bf
R}^3 / W) \times {\bf R}$ onto ${\bf R}^4$ in Theorem \ref{(R^3/W)
times R is homeomorphic to R^4} to a homeomorphism from $({\bf S}^3
\backslash W) \times {\bf R}$ onto ${\bf R}^4$, then in particular one
could be lead to try to figure out something about what happens when
one deletes $\mathcal{L}$ from ${\bf R}^4$.  Conversely, if one wanted
to go in the other direction, one might have to figure out something
about how to put the topological line back in.  These endeavors should
be at least somewhat complicated, because of the wildness of
$\mathcal{L}$ inside of ${\bf R}^4$.

	The next fact helps to give an idea of how wild $\mathcal{L}$ can
have to be.

\begintheorem
\label{hausdorff dimension and local fundamental groups}
Let $U$ be an open set in some ${\bf R}^n$, and let $F$ be a closed set
in ${\bf R}^n$.  If the Hausdorff dimension of $F$ is less than
$n-2$, then any open loop in $U \backslash F$ that can be contracted
to a point in $U$ can also be contracted to a point in $U \backslash F$.
In particular, ${\bf R}^n \backslash F$ is simply-connected.
\end{theorem}

	In other words, if $F$ is closed and has Hausdorff dimension
$< n-2$, then $F$ is practically invisible for considerations of
fundamental groups, even local ones.  From this one can check that
$\mathcal{L}$ as above necessarily has Hausdorff dimension at least $2$,
no matter the homeomorphism from $({\bf R}^3 / W) \times {\bf R}$ onto
${\bf R}^4$ which produced it.  This is also true locally, i.e., each
nontrivial arc of $\mathcal{L}$ has to have Hausdorff dimension at least
$2$, and for the same reasons.  This is somewhat remarkable, since
$\mathcal{L}$ is homeomorphic to a line (and straight lines have
Hausdorff dimension $1$).  (In general, Hausdorff dimension need not
be preserved by homeomorphisms, though, and this is an instance of
that.)

	Theorem \ref{hausdorff dimension and local fundamental groups}
is given (in a slightly different form) in \cite{MaRV}, in Lemma 3.3
on p9.  See also \cite{Ge-link, LuV, SS, V1} for related results.  In
particular, \cite{SS} uses Theorem \ref{hausdorff dimension and local
fundamental groups} in a manner very similar to this, in the context
of double-suspension spheres and homeomorphic parameterizations of
them.  

	Note that instead of requiring that $F$ have Hausdorff
dimension less than $n-2$ in Theorem \ref{hausdorff dimension and
local fundamental groups}, it is enough to ask that the
$(n-2)$-dimensional Hausdorff measure of $F$ be zero.  

	By now, there are many examples known of spaces which are
manifold factors (and not manifolds).  The original discovery was in
\cite{Bing3, Bing4}, using a different space.  This example, Bing's
``dogbone'' space, will come up again in Subsection \ref{Decomposition
spaces, 2}.

	As mentioned near the bottom of p93 of \cite{Dm2}, $3$ is the
smallest dimension in which this type of phenomenon can occur (where a
non-manifold becomes a manifold after taking the Cartesian product
with ${\bf R}$), because of results about recognizing manifolds in
dimensions $1$ and $2$.  It does occur in all dimensions greater than
or equal to $3$.

	As a basic class of examples, if $K$ is a compact set in ${\bf
R}^n$ which is a cell, then ${\bf R}^n / K$ may not be a manifold if
$K$ is wild, but $({\bf R}^n / K ) \times {\bf R}$ is homeomorphic to
${\bf R}^{n+1}$.  See \cite{Andrews-Curtis, Bryant1, Bryant2, Dm2,
Ru-book}.

	More generally, if $K$ is any compact set in ${\bf R}^n$ which
is cell-like, then $({\bf R}^n / K) \times {\bf R}$ is homeomorphic to
${\bf R}^{n+1}$.  When $n=1$ or $2$, ${\bf R}^n / K$ is itself
homeomorphic to ${\bf R}^n$.  (Compare with Proposition
\ref{cellularity and cell-likeness} and Theorem \ref{R^n/K a
topological manifold if and only if K is cellular} in Subsubsection
\ref{Cellularity, and the cellularity criterion}.)  If $n \ge 3$, then
${\bf R}^n / K$ may not be homeomorphic to ${\bf R}^n$, but it is true
that $({\bf R}^n / K) \times {\bf R}$ is homeomorphic to ${\bf
R}^{n+1}$ in this situation.  See Proposition 2 on p206 of \cite{Dm2}
for $n=3$, and Theorem 9 on p196 or Theorem 13 on p200 of \cite{Dm2}
for $n \ge 4$.  (For Theorem 9 on p196 of \cite{Dm2}, note that the
definition of a ``$k$-dimensional decomposition'' of a manifold is
given near the top of p152 of \cite{Dm2}.)

	Another class of examples (which is in fact closely related to
the previous ones) comes from the celebrated double-suspension
theorems of Edwards and Cannon \cite{C1, C-bulletin, C2, Dm2, E}
(mentioned in Section \ref{parameterization problems}).  From these
one has the remarkable fact that there are \emph{finite polyhedra} $P$
which are not topological manifolds, but for which $P \times {\bf R}$
is a manifold.

	There are compact sets $K$ in ${\bf R}^n$ such that $K$
is \emph{not} cell-like, and not cellular in particular, but $({\bf
R}^n / K) \times {\bf R}$ is homeomorphic to ${\bf R}^{n+1}$.  This
happens for every $n \ge 4$.  See Corollary 3E on p185 of \cite{Dm2}.
(The proof uses the double-suspension theorems.)

	There are also non-manifold spaces which become manifolds
after taking the product with other non-manifold spaces.  See Section
29 of \cite{Dm2}, beginning on p223.

\subsection{Decomposition spaces, 2}
\label{Decomposition spaces, 2}

	The construction of the quotient ${\bf R}^n / K$, given a
compact subset $K$ of ${\bf R}^n$, as in Subsection \ref{Decomposition
spaces, 1}, is an example of a \emph{decomposition space}.  More
generally, one can allow many subsets of ${\bf R}^n$ (or some other
space) to be contracted to points at the same time, rather than just a
single set.

	In general, a \emph{decomposition} of ${\bf R}^n$, or some
other space, means a partition of it, i.e., a collection of subsets
which are pairwise disjoint, and whose union is the whole space.  One
can then form the corresponding quotient space, first as a set --- by
collapsing the sets in the partition to individual points --- and then
as a topological space.  The topology on the quotient is the richest
one (the one with the most open sets) so that the canonical mapping
from the space to the quotient is continuous.

	All of this makes sense in general, but in order to have some
good properties (like the Hausdorff condition for the quotient space),
some assumptions about the decomposition should be made.  As a start,
typical assumptions would be that the individual sets that make up the
decomposition be closed, and that the decomposition satisfy a certain
upper semi-continuity property.  See \cite{Dm2} for details, including
the definition on p13 of \cite{Dm2}.  For the present discussion, we
shall always assume that some conditions like these hold, even if we
do not say so explicitly.  

	If $G$ is a decomposition of ${\bf R}^n$, then one writes
${\bf R}^n / G$ for the corresponding quotient space.  

	Given a set $K$ in ${\bf R}^n$, one can always consider the
decomposition of ${\bf R}^n$ consisting of $K$ and sets with only
single elements, with the latter running through all points in ${\bf
R}^n \backslash K$.  This decomposition is sometimes denoted $G_K$,
and the quotient ${\bf R}^n / G_K$ in the general sense of
decompositions is the same as the space ${\bf R}^n / K$ from
Subsection \ref{Decomposition spaces, 1}.

	As another basic situation, product spaces of the form $({\bf
R}^n / K) \times {\bf R}$ can be viewed as decomposition spaces.
Specifically, one gets a decomposition of ${\bf R}^{n+1} \cong {\bf
R}^n \times {\bf R}$ using sets of the form $K \times \{u\}$ in ${\bf
R}^{n+1}$ for each $u \in {\bf R}$, together with single-element sets
for all of the points in ${\bf R}^{n+1} \backslash (K \times {\bf
R})$.  The resulting decomposition space is equivalent topologically
to $({\bf R}^n / K) \times {\bf R}$.

	An important general point is that wild or interesting
embeddings can often occur in simple or useful ways through
decompositions.  For instance, in the decomposition space described in
the preceding paragraph, one has a particular ``line'', corresponding
to the copies of $K$.  See \cite{C1, Dm1, Dm2} for more information,
including p451 of \cite{C1}, the last paragraph in Section 2 on p380
of \cite{Dm1}, and the remarks near the top of p37 in \cite{Dm2}.

	The following theorem of R. L. Moore \cite{Moore} is an early
result about when decomposition spaces are manifolds, and homeomorphic
to the original space.

\begintheorem
\label{Moore's theorem}
Let $X$ be a compact Hausdorff topological space.  Suppose that $f$
is a continuous mapping from the $2$-sphere ${\bf S}^2$ to $X$,
and that for each $x \in X$, $f^{-1}(x)$ is nonempty and connected,
and ${\bf S}^2 \backslash f^{-1}(x)$ is nonempty and connected.  Then
$X$ is homeomorphic to ${\bf S}^2$.
\end{theorem}

	In this theorem, the mapping $f$ itself may not be a
homeomorphism.  As in Subsection \ref{Decomposition spaces, 1}, $f$
might have the effect of collapsing some line segments down to single
points, for instance.  It is true that $f$ can always be approximated
by homeomorphisms, however.  See \cite{Dm2} for more information and
references.

	Similar results hold in dimension $1$.  For this it is enough
to assume that the inverse images of points under the mapping be
\emph{connected} (and nonempty proper subsets of the domain), without
imposing conditions on their complements.  In this situation, the
inverse images of points will simply be arcs.  Compare with
\cite{Dm2}, including the remarks near the bottom of p17.

	What might be reasonable analogues of Theorem \ref{Moore's
theorem} in higher dimensions?  One should not keep the hypotheses
literally as they are above, where the fibers are connected and have
connected complements, because of counterexamples like the
non-manifold spaces that one can get by contracting a circle to a
point (as in Subsection \ref{Decomposition spaces, 1}).

	However, in dimension $2$, the property of a set in ${\bf
S}^2$ or ${\bf R}^2$ being connected and having connected complement
is quite strong.  For a closed subset of ${\bf S}^2$ which is not
empty nor all of ${\bf S}^2$, these conditions imply that the
complement of the set is homeomorphic to a $2$-dimensional disk, and
that the set itself is \emph{cellular} (Definition \ref{definition of
cellularity} in Subsection \ref{Decomposition spaces, 1}).

	A decomposition $G$ of ${\bf R}^n$ is said to be \emph{cellular}
if each of the subsets of ${\bf R}^n$ of which it is composed is cellular.
(Compare with \cite{Dm2}, in the statement of Corollary 2A on p36.)
As an analogy with Moore's theorem, one might hope that a quotient of
${\bf R}^n$ (or ${\bf S}^n$, or other topological manifolds) by a cellular
decomposition is a manifold, and homeomorphic to ${\bf R}^n$ again (or to
the original manifold, whatever it might be).

	This is true for decompositions which consist of a single
cellular subset of the space, together with all the remaining points
in the space as one-element sets.  In other words, this statement is
true in the context of quotients as in Subsection \ref{Decomposition
spaces, 1}.  See Theorem \ref{R^n/K a topological manifold if and only
if K is cellular} in Subsection \ref{Decomposition spaces, 1}, and
Proposition 2 on p36 of \cite{Dm2}.

	For decompositions in general, it is not true that cellularity
is sufficient to ensure that the quotient space is a manifold.  This
fails already for decompositions of ${\bf R}^3$.  The first example of
this was provided by Bing's ``dogbone'' construction in \cite{Bing2}.
See also \cite{Dm2}, especially Section 9, for this and other
examples.  For the statement that the decomposition space is not a
manifold, and not just not homeomorphic to ${\bf R}^3$, see Theorem 13
on 498 of \cite{Bing2}.  A recent paper related to this is
\cite{Armentrout}.

	The dogbone space was also used in the initial discovery of
manifold factors.  See \cite{Bing3, Bing4}.

	Although quotients by cellular decompositions do not in
general give manifolds, there are many nontrivial examples where this
does occur, and results about when it should take place.  A
particularly nice and fundamental example is given by ``Bing
doubling''.  See \cite{Bing1, Dm2} (Example 1 in Section 9 of
\cite{Dm2}).  This is a decomposition of ${\bf R}^3$ for which the
corresponding quotient space is homeomorphic to ${\bf R}^3$.  While
the quotient space is standard, the decomposition has some interesting
features, giving rise to some wild embeddings in particular.  This
decomposition has a symmetry to it, which leads to a homeomorphic
involution on ${\bf S}^3$ which is highly nonstandard.  The
fixed-point set of this homeomorphism is a wildly-embedded $2$-sphere
in ${\bf S}^3$.  This construction apparently gave the first examples
of involutions on ${\bf R}^3$ which were not topologically equivalent
to ``standard'' ones, made from rotations, reflections, and
translations.  See \cite{Bing1} for more information, especially
Section 4 of \cite{Bing1}.

	See also Section 9 of \cite{Bing4} for a wild involution on
${\bf R}^4$, whose fixed point set is homeomorphic to the dogbone
space.  This uses the fact that the product of the dogbone space with
the real line is homeomorphic to ${\bf R}^4$.

	The results mentioned in Subsection \ref{Manifold factors} ---
about quotients ${\bf R}^n / K$ being homeomorphic to a Euclidean
space after taking the Cartesian product with ${\bf R}$ --- can also be
seen as providing nontrivial examples of cellular decompositions of
${\bf R}^{n+1}$ for which the corresponding quotients are manifolds,
and are homeomorphic to ${\bf R}^{n+1}$.

\subsection{Geometric structures for decomposition spaces}
\label{Geometric structures for decomposition spaces}

\subsubsection{A basic class of constructions}
\label{A basic scenario}

	One feature of decomposition spaces is that they do not a
priori come with a canonical or especially nice \emph{geometry}, or
anything like that.  The \emph{topology} is canonical, but this is
somewhat different.  Note that there are general results about
existence of metrics which are compatible with the topology, as in
Proposition 2 on p13 of \cite{Dm2}.  Once one has one such metric,
there are many which define the same topology.  This is true just as
well for ordinary Euclidean spaces, or the spheres ${\bf S}^n$, even
if there are also special metrics (like the Euclidean metric) that one
might normally like or use.

	In some cases (for decomposition spaces), there are some
particularly nice or special geometries that one can consider.  A
number of basic examples --- like the Whitehead continuum, Bing
doubling \cite{Bing1}, and Bing's dogbone space \cite{Bing2} --- have a
natural topological ``self-similarity'' to them, which can be
converted into geometric self-similarity.

	Let us be more precise.  In these cases, the nondegenerate
elements of the decomposition are generated by repeating a simple
``rule''.  The ``rule'' can be described by a smooth domain $D$ in
${\bf R}^n$ together with some copies of $D$ embedded in the interior
of $D$, in a pairwise-disjoint manner.  To generate the nontrivial
elements of the decomposition, one starts with $D$, and then passes to
the copies of $D$ inside of itself.  For each of these copies of $D$
inside of $D$, one can get a new collection of smaller copies of $D$,
by applying the ``rule'' to the copies of $D$ with which we started.
One then repeats this indefinitely.

	Thus, if the original ``rule'' involves $m$ copies of $D$
inside of itself, then the $j$th step of this process gives rise to
$m^{j-1}$ copies of $D$, with the first step corresponding to $D$
alone.

	The limiting sets which arise from this procedure are pairwise
disjoint, and are used to define the decomposition.  To do this
carefully, one can think of the $j$th step of the process as producing
a compact set $C_j$, which is the union of the $m^{j-1}$ individual
copies of $D$ indicated above.  The construction gives $C_j \subseteq
C_{j-1}$ for all $j \ge 2$.  To pass to the limit, one can take the
set $C = \bigcap_{j=1}^\infty C_j$, and then use the components of $C$
as subsets of ${\bf R}^n$ to be employed in the decomposition $G$ of
${\bf R}^n$.  For each point $x \in {\bf R}^n \backslash C$, one also
includes the one-point set $\{x\}$ in the decomposition $G$.

	The Whitehead continuum (discussed in Subsection
\ref{Contractable open sets}) is an example of this.  There the
``rule'' is particularly simple, in that it is based on an embedding
of a single solid torus $T_1$ inside a larger one $T$.  At each stage
of the process there is only one domain, and only one nondegenerate
set being produced in the end (i.e., the Whitehead continuum).  In
particular, one can have $m=1$ in the general set-up described above,
and with the result being nontrivial.  For Bing doubling and Bing's
dogbone space, one has $m > 1$, i.e., there are more than one
embedding being used at each step, and more than one copy of the basic
domain.  (For Bing doubling the basic domain is again a solid torus,
while for the dogbone space it is a solid $2$-handled torus.  In Bing
doubling one has $m=2$, as the name suggests.  For the dogbone space,
$m=4$.)  In these cases the number of components grows exponentially
in the process, and is infinite in the end (after taking the limit).

	At any rate, there are a number of basic examples of
decompositions generated in this manner, with individual copies of a
domain being replaced systematically by some embedded subdomains,
copies of itself, following repetitions of a single basic ``rule''.
See \cite{Dm2}, especially Section 9.

	In typical situations, the basic ``rule'' involves nontrivial
distortion of the standard Euclidean geometry at each step.  As a
first point along these lines, when one embeds a copy of some
(bounded) domain $D$ into itself, then some change in the geometry
(along the lines of shrinkage) is unavoidable, at least if the
embedded copy is a proper subset of the original domain.

	Given that there needs to be some shrinkage in the embedding,
the next simplest possibility would be that the embeddings are made up
out of something like dilations, translations, rotations, and
reflections.  In other words, except for a uniform scale factor, one
might hope that the geometry does not have to change.

	Normally this will not be the case.  Some amount of bending or
twisting, etc., will (in general) be involved, and needed, to
accommodate the kind of topological behavior that is present.  This
includes linking, clasping, or things like that.

	For the purpose of \emph{choosing} a geometry that might fit
with a given decomposition space, however, one can modify the usual
Euclidean metric so that the embeddings involved in the basic ``rule''
do have the kind of behavior indicated above, i.e., a constant scale
factor together with an isometry.  The scale factors should be less
that $1$, to reflect the shrinking that is supposed to take place for
the decomposition spaces (even at a purely topological level).

	It is not hard to see that one can make deformations of
geometry like this.  One can do this in a kind of direct and
``intrinsic'' way, defining metrics on ${\bf R}^n$ with suitable
properties.  One can also do this more concretely, ``physically'',
through embeddings of the decomposition spaces into higher-dimensional
Euclidean spaces.  In these higher-dimensional Euclidean spaces, the
self-similarity that one wants, in typical situations, can be realized
in terms of standard \emph{linear} self-similarity, through dilations
and translations.

	More precisely, in these circumstances, the quotient of ${\bf
R}^n$ by the decomposition can be realized topologically as an
$n$-dimensional subset $X$ of some ${\bf R}^N$ (with $N = n+1$, for
instance), in such a way that $X$ is a smooth submanifold away from
the natural singularities, and $X$ is self-similar around these
singularities.  

	To build such a set $X$, one can start with the complement of
the original domain $D$ in ${\bf R}^n$.  One would view ${\bf R}^n
\backslash D$ as an $n$-dimensional submanifold of ${\bf R}^N$.  In
place of the iteration of the basic rule for the decomposition from
before, one now stacks some ``basic building blocks'' in ${\bf R}^N$
on top of ${\bf R}^n \backslash D$ (along the boundary of $D$), and
then on top of the other building blocks, over and over again.

 	These basic building blocks are given by $n$-dimensional
smooth manifolds in ${\bf R}^N$ (with boundary).  They are
diffeomorphic to a single ``model'' in ${\bf R}^n$, which is the
original domain $D$ in ${\bf R}^n$, minus the interiors of the $m$
copies of $D$ embedded inside $D$, as given by the basic ``rule'' that
generates the decomposition.  The building blocks are all
diffeomorphic to each other, since they are all diffeomorphic to this
same model, but they are also constructed in such a way as to be
``similar'' to each other.  That is, they can all be given by
translations and dilations of each other.  This is a key difference
between this construction and the original decomposition in ${\bf
R}^n$.

	Further, the building blocks are constructed in such a way
that their ends are all similar to each other (i.e., even different
ends on the same building block).  Specifically, the building blocks
are chosen so that when one goes to stack them on top of each other,
their ``ends'' fit together properly, with smoothness across the
interfaces.

	These things are not difficult to arrange.  Roughly speaking,
one uses the extra dimensions in ${\bf R}^N$ to straighten the
``ends'' in this way, so that the different building blocks can be
stacked properly.  Typically, this would involve something like the
following.  One starts with the basic model in ${\bf R}^n$, given by
$D$ minus the interiors of the $m$ embedded copies of $D$ in $D$.  One
then makes some translations of the $m$ embedded subdomains in $D$, up
into the extra dimension or dimensions in ${\bf R}^N$.  Up there,
these subdomains can be moved or bent around, until they are similar
to $D$ itself (i.e., being the same modulo translations and
dilations).  This can be done one at a time, and without changing
anything near the boundary of the original domain $D$.  In this
manner, the original model region in ${\bf R}^n$ becomes realized as
an $n$-dimensional compact smooth submanifold (with boundary) in ${\bf
R}^N$, with the ends matching up properly under similarities.

	To put it another way, the main ``trade-off'' here is that one
gives up the ``flatness'' of the original model, as a region in ${\bf
R}^n$, to get basic building blocks in ${\bf R}^N$ that are
$n$-dimensional curved submanifolds whose ends are similar to each
other.  The curving of the interiors of these building blocks
compensates for the straightening of their ends.

	As above, one then stacks these building blocks on top of each
other, one after another, to get a realization of the decomposition
space by an $n$-dimensional subset $X$ of ${\bf R}^N$.  (One also puts
in some limiting points, at the ends of the towers of the building
blocks that arise.  In other words, this makes $X$ be a \emph{closed}
subset of ${\bf R}^N$.  These extra points are the singularities of
$X$.)  This subset is smooth away from the singularities, and
self-similar at the singularities, because of the corresponding
properties of the basic building blocks.

	By choosing the scale-factors associated to the ends of the
basic building blocks to be less than $1$, the diameters of the ends
tend to $0$ (and in a good way) as one stacks the building blocks on
top of each other many times.  This corresponds to the fact that the
sets in the decomposition are supposed to be shrunk to single points
in the quotient space.  This is also part of the story of the
``limiting points'' in the previous paragraph.  The limiting points
are exactly the ones associated (in the end) to the nondegenerate sets
in the original decomposition in ${\bf R}^n$, which are being shrunk
to single points.

	The actual homeomorphic equivalence between the set $X$ in
${\bf R}^N$ produced through this method and the decomposition space
${\bf R}^n / G$ with which one starts is obtained using the
diffeomorphic equivalence between the building blocks in ${\bf R}^N$
and the original model in ${\bf R}^n$ ($D$ minus the interiors of the
$m$ embedded copies of itself, as above).  In rough terms, at the
level of the topology, the same kind of construction is occurring in
both places, $X$ and the decomposition space, and one can match them
up, by matching up the individual building blocks.  This is not hard
to track.

	At any rate, more details for these various matters are given
in \cite{Seqs}.

	Instead of stacking basic building blocks on top of each other
infinitely many times, one can stop after finitely many steps of the
construction.  This gives a set which is still smooth, and
diffeomorphic to ${\bf R}^n$, but which approximates the non-smooth
version that represents the decomposition space.  To put it another
way, although this kind of approximation is standard
\emph{topologically}, its \emph{geometry} does reflect the basic rule
from which the decomposition space is obtained.  In particular,
properties of the decomposition space can be reflected in questions of
quantitative bounds for the approximations which may or may not hold.

	Some versions of this come up in \cite{Seqs}.  In addition,
there are some slightly different but related constructions, with
copies of finite approximations repeated but getting small, at the
same time that they include more and more stages in the stacking.

	The topological structure of the decomposition space can also
be seen in terms of Gromov--Hausdorff limits of smooth approximations
like these.

	When one makes constructions like these --- either finite
approximations or infinite limits --- the self-similarity helps to
ensure that the spaces behave geometrically about as well as they
could.  For instance, this is manifested in terms of properties like
Ahlfors-regularity (Definition \ref{ahlfors regular}), and the local
linear contractability condition (Definition \ref{local lin con}), at
least under suitable assumptions or choices.  One also gets good
behavior in terms of ``calculus'', like Sobolev and Poincar\'e
inequalities.  One has uniform rectifiability as well (Definition
\ref{unif rec}), although the spaces are actually quite a bit more
regular than that.  For smooth approximations with only finitely many
levels in the construction, one can have uniform bounds, independent
of the number of levels, for conditions like these.  See \cite{Seqs}
for more information, and some slightly different versions of these
basic themes.

	To summarize a bit, although the quotient spaces that one gets
from decompositions may not be topological manifolds, in many cases
one can realize them geometrically in such a way that the behavior
that occurs is really pretty good.  With the extra structure from the
geometry, there are extra dimensions to the story as a whole, and
which are perhaps not apparent at first.  (Possibilities for doing
analysis on spaces like these, and spaces which are significantly
different from standard Euclidean spaces at that, give one form of
this.)

	We shall look at these types of geometries from some more
perspectives in the next subsubsections, and also consider some
special cases.  

	These general matters are related as well to the topics of
Appendix \ref{Working on spaces which may not have nice coordinates}.

\subsubsection{Comparisons between geometric and topological properties}
\label{A couple of comparisons between geometric and topological properties}

	Many of the geometric properties that occur when one
constructs geometries for decomposition spaces as above correspond in
natural ways to common conditions in purely topological terms that one
might consider.

	For instance, if one starts with a decomposition space ${\bf
R}^n / G$ from ${\bf R}^n$, and gives it a metric in which it becomes
Ahlfors-regular of dimension $n$, then in particular the decomposition
space has Hausdorff dimension $n$ with respect to this metric.  In the
type of examples mentioned above, ${\bf R}^n / G$ would also have
topological dimension $n$.

	The Hausdorff dimension of a metric space is always greater
than or equal to the topological dimension, as in Chapter VII of
\cite{HW}.  To have the two be equal is rather nice, and
Ahlfors-regularity in addition even nicer.  For the constructions
described in Subsubsection \ref{A basic scenario}, Ahlfors-regularity
reflects the fact that one chooses the geometries to be smooth away
from the singularities, and self-similar around the singularities.

	For cell-like decompositions in general (of compact metric
spaces, say), there is a theorem which says that the topological
dimension of the quotient is less than or equal to the topological
dimension of the space with which one starts, if the topological
dimension of the quotient is known to be finite.  See Theorem 7 on
p135 of \cite{Dm2}, and see \cite{Dm2} in general for more
information.  Note that it is possible for the topological dimension
of the quotient to be infinite, by a result of Dranishnikov
\cite{Dranishnikov} (not available at the time of the writing of
\cite{Dm2}).

	For the situations discussed in Subsubsection \ref{A basic
scenario}, Ahlfors-regularity of dimension $n$ for the geometry of the
decomposition space is a nice way to have good behavior related to the
dimension (and with the right value of the dimension, i.e., the
topological dimension).  In other words, it is a geometric property
that fits nicely with general topological considerations (as in the
preceding two paragraphs), while also being about as strong as it can
be.  It is also quantitative, and it makes sense to consider uniform
bounds for smooth approximations that use only finitely many levels of
the construction.  (Note that it is automatically preserved under
Gromov--Hausdorff limits, when there are uniform bounds.  See
\cite{DS-ff} for more information.)

	Another general result about decomposition spaces implies that
a cell-like quotient of a space is \emph{locally contractable} if the
quotient has finite dimension, under suitable conditions on the space
from which one starts (e.g., being a manifold).  See Corollary 12B on
p129 of \cite{Dm2} (as well as Corollaries 4A and 8B on p115 and p119
of \cite{Dm2} for auxiliary statements).  Here \emph{local
contractability} of a topological space $X$ means that for every point
$p \in X$, and every neighborhood $U$ of $p$ in $X$, there is a
smaller neighborhood $V \subseteq U$ of $p$ in $X$ such that $V$ can
be contracted to a point in $U$.  If the quotient space is not known
to be finite-dimensional, then there are weaker but analogous
conclusions that one can make, along the lines of local
contractability.  See Corollary 11B on p129 of \cite{Dm2}.

	These notions are purely topological ones.  The local linear
contractability condition (Definition \ref{local lin con}) is a more
quantitative version of this for metric spaces.  It is about as strong
as one can get, in terms of the sizes of the neighborhoods concerned.
That is, if one tries to contract a small ball in the space to a point
inside of another ball which is also small, but somewhat larger than
the first one, then weaker conditions besides local linear
contractability might allow somewhat larger radii for the second ball
(than a constant multiple of the initial radius).  The \emph{linear}
bound fits naturally with self-similarity.

	As indicated earlier (towards the end of Subsubsection \ref{A
basic scenario}), local linear contractability does come up in a
natural manner for the constructions in Subsubsection \ref{A basic
scenario}.  This fits nicely with the general topological results,
i.e., by giving a geometric form which reflects the topological
properties in about the best possible way.  
	
	Exactly how nice the spaces discussed in Subsubsection \ref{A
basic scenario} are depends on the features of the decompositions.
This includes not only how many singularities there are, but also how
strong they are.  For instance, cellular decompositions can lead to
nicer or more properties for the resulting spaces than cell-like
decompositions.  I.e., at the geometric level, cellularity can give
rise to more good behavior than local linear contractability.  For
example, one can look at simple-connectivity of \emph{punctured}
neighborhoods in this regard, and things like that.  This is something
that one does not have for the space obtained by taking ${\bf R}^3$
and contracting a Whitehead continuum to a single point.  Some other
versions of this are considered in \cite{Seqs}.

\subsubsection{Quotient spaces can be topologically standard, but 
geometrically tricky}
\label{Quotient spaces can be topologically standard, but geometrically tricky}

	We have seen before how decompositions of ${\bf R}^n$ might
lead to ${\bf R}^n$ again topologically in the quotient, but do so in
a manner that is still somehow nontrivial.  For instance, the
decomposition might arise from a nontrivial manifold factor, or lead
to wild embeddings in the quotient which seem very simple (like a
straight line) at the level of the decomposition.  In these
situations, one can still have highly nontrivial geometries from the
procedures described in Subsubsection \ref{A basic scenario}, even
though the underlying space is topologically equivalent to ${\bf
R}^n$.

	As a special case, the wild embeddings that can occur in the
quotient (with the topological identification with a Euclidean space)
can behave in a very special way in the geometric realization of the
decomposition space that one has here, with geometric properties which
are not possible in ${\bf R}^n$ with the standard Euclidean metric.
We shall give a concrete example of this in a moment.

	This is one way that geometric structures for decomposition
spaces (as in Subsubsection \ref{A basic scenario}) can help to add
more to the general story.  One had before the notion that a wild
embedding might have a simple realization through a decomposition
space, and now this might be made precise through metric properties of
the embedding (like Hausdorff dimension), using the type of geometry
for the decomposition space that one has here.

	Here is a concrete instance of this.  Let $W$ be a Whitehead
continuum in ${\bf R}^3$, as in Subsection \ref{Contractable open
sets}.  Consider the corresponding quotient space ${\bf R}^3 / W$, as
in Subsection \ref{Decomposition spaces, 1}.  Through the type of
construction described in Subsubsection \ref{A basic scenario}, one
can give ${\bf R}^3 / W$ a geometric structure with several nice
properties.  In short, one can realize it topologically as a subset of
${\bf R}^4$, where this subset is smooth away from the singular point,
and has a simple self-similarity at the singular point.  One could
then use the geometry induced from the usual one on the ambient ${\bf
R}^4$.

	Let us write $q$ for the singular point in ${\bf R}^3 / W$,
i.e., the point in the quotient which corresponds to $W$.  Using the
metric on ${\bf R}^3 / W$ indicated above, let us think of $({\bf R}^3
/ W) \times {\bf R}$ as being equipped with a metric.  One could also
think of $({\bf R}^3 / W) \times {\bf R}$ as being realized as a
subset of ${\bf R}^5$, namely, as the product of the one in ${\bf
R}^4$ mentioned before with ${\bf R}$.

	In this geometry, $L = \{q\} \times {\bf R}$ is a perfectly
nice line.  It is a ``straight'' line!  In particular, it has
Hausdorff dimension $1$, and locally finite length.  However, the
image of $L$ inside of ${\bf R}^4$ under a homeomorphism from $({\bf
R}^3 / W) \times {\bf R}$ onto ${\bf R}^4$ will be wild.  As in
Subsection \ref{Manifold factors}, any such embedding of $L$ in ${\bf
R}^4$ must be wild, and in fact any such embedding must have Hausdorff
dimension at least $2$, with respect to the usual Euclidean metric in
${\bf R}^4$.  This used Theorem \ref{hausdorff dimension and local
fundamental groups}.

	This shows that the geometry that we have for $({\bf R}^3 / W)
\times {\bf R}$ has to be substantially different from the usual
Euclidean geometry on ${\bf R}^4$, even though the two spaces are
topologically equivalent.  Specifically, even though there are
homeomorphisms from $({\bf R}^3 / W) \times {\bf R}$ onto ${\bf R}^4$,
no such homeomorphism can be Lipschitz (with respect to the geometries
that we have on these spaces), or even H\"older continuous of order
larger than $1/2$.

	Although $({\bf R}^3 / W) \times {\bf R}$ --- with the kind of
geometry described above --- is quite different from ${\bf R}^4$ with
the usual Euclidean metric, there is a strong and nice feature that it
has, in common with ${\bf R}^4$.  We shall call this property
``uniform local coordinates''.

	Since $({\bf R}^3 / W) \times {\bf R}$ is homeomorphic to
${\bf R}^4$, it has homeomorphic local coordinates from ${\bf R}^4$ at
every point.  ``Uniform local coordinates'' asks for a stronger
version of this, and is more quantitative.  Specifically, around each
metric ball $B$ in $({\bf R}^3 / W) \times {\bf R}$ (with respect to
the kind of geometry that we have), there are homeomorphic local
coordinates from a standard Euclidean ball $\beta$ of the same radius
in ${\bf R}^4$, such that
\begin{eqnarray}
\label{image of beta covers B}
	&& \hbox{the image of $\beta$ under the coordinate mapping}   \\
 	&& \hbox{covers the given ball $B$ in }
                    ({\bf R}^3 / W) \times {\bf R},		\nonumber
\end{eqnarray}
and
\begin{eqnarray}
\label{uniform moduli of continuity, after rescaling}
	&& \hbox{the \emph{modulus of continuity} of the coordinate mapping} \\
        && \hbox{and its inverse can be controlled, uniformly over all}
                                                                \nonumber    \\
        && \hbox{choices of metric balls $B$ in } 
                    ({\bf R}^3 / W) \times {\bf R},
                              \hbox{ and in a}                  \nonumber    \\
        && \hbox{scale-invariant manner.}                       \nonumber
\end{eqnarray}

	Here ``modulus of continuity'' means a function $\omega(r)$ so
that when two points in the domain (of a given mapping) are at
distance $\le r$, their images are at distance $\le \omega(r)$.
Also, $r$ would range through positive numbers, and $\omega(r)$ would
be nonnegative and satisfy 
\begin{equation}
	\lim_{r \to 0} \omega(r) = 0.
\end{equation}
This last captures the continuity involved, and, in fact, gives
uniform continuity.

	For a mapping from a compact metric space to another metric space,
continuity automatically implies uniform continuity, which then implies the
existence of some modulus of continuity.  This is not to say that one knows
much about the modulus of continuity, a priori.  (One can always choose it
to be monotone, for instance, but one cannot in general say how fast it
tends to $0$ as $r \to 0$.)

	Concrete examples of moduli of continuity would include
$\omega(r) = C \, r$ for some constant $C$, which corresponds to a
mapping being Lipschitz with constant $C$, or $\omega(r) = C \,
r^\alpha$, $\alpha > 0$, which corresponds to H\"older continuity of
order $\alpha$.  One can have much slower rates of vanishing, such as
$\omega(r) = (\log \log \log (1/r))^{-1}$.

	In our case, with the property of ``uniform local
coordinates'', we want to have a single modulus of continuity
$\omega(r)$ which works simultaneously for all of the local coordinate
mappings (and their inverses).  Actually, we do not look at moduli of
continuity for the mappings themselves, but renormalized versions of
them.  The renormalizations are given by dividing distances in the
domain and range by the (common) radius of $B$ and $\beta$.  In this
way, $B$ and $\beta$ are viewed as though they have radius $1$,
independently of what the radius was originally.  This gives a kind of
uniform basis for making comparisons between the behavior of the
individual local coordinate mappings.

	Let us return now to the special case of $({\bf R}^3 / W)
\times {\bf R}$, with the geometry as before.  In this case one can
get the condition of uniform local coordinates from the existence of
topological coordinates (without uniform bounds), together with the
self-similarity and smoothness properties of the set.  Here is an
outline of the argument.  (A detailed version of this, for a modestly
different situation, is given in \cite{Seqs}.  Specifically, see
Theorem 6.3 on p241 in \cite{Seqs}.  Note that the property of uniform
local coordinates is called ``Condition ($**$)'' in \cite{Seqs}, as in
Definition 1.7 on p192 of \cite{Seqs}.)

	Let $B$ be a metric ball in $({\bf R}^3 / W) \times {\bf R}$,
for which one wants to find suitable coordinates.  Assume first that
$B$ does not get too close to the singular line $\{q\} \times {\bf R}$
in $({\bf R}^3 / W) \times {\bf R}$, and in fact that the radius of
$B$ is reasonably small compared to the distance from $B$ to the
singular line.  In this case $({\bf R}^3 / W) \times {\bf R}$ is
pretty smooth and flat in $B$, by construction (through the method of
Subsubsection \ref{A basic scenario}).  This permits one to get local
coordinates around $B$ quite easily, and with suitable uniform bounds
for the moduli of continuity of the coordinate mappings and their
inverses.  The bounds that one gets are scale-invariant, because of
the self-similarity in the geometric construction (from Subsubsection
\ref{A basic scenario}).  In fact, one can have Lipschitz bounds in
this case, as well as stronger forms of smoothness.

	If the ball $B$ is reasonably close to the singular line
$\{q\} \times {\bf R}$, then one can reduce to the case where it is
actually centered on $\{q\} \times {\bf R}$.  That is, one could
replace $B$ with a ball which is centered on $\{q\} \times {\bf R}$,
and which is not too much larger.  (The radius of the new ball would
be bounded by a constant times the radius of $B$.)  This substitution
does not cause trouble for the kind of bounds that are being sought
here.

	Thus we suppose that $B$ is centered on the line $\{q\} \times
{\bf R}$.  We may as well assume that the center of $B$ is the point
$(q,0)$.  This is because $({\bf R}^3 / W) \times {\bf R}$ and the
geometry that we have on it are invariant under translations in the
${\bf R}$ direction, so that one can move the center to $(q,0)$
without difficulty, if necessary.

	Using the self-similarity of $({\bf R}^3 / W) \times {\bf R}$,
one can reduce further to the case where the radius of $B$ is
approximately $1$.  For that matter, one can reduce to the case where
it is equal to $1$, by simply increasing the radius by a bounded
factor (which again does not cause problems for the uniform bounds
that are being considered here).  (To be honest, if one takes the
geometry for ${\bf R}^3 / W$ to be flat outside a compact set, as in
Subsubsection \ref{A basic scenario}, then this reduction is not fully
covered by self-similarity.  I.e., one should handle large scales a
bit differently.  This can be done, and a similar point is discussed
in \cite{Seqs} in a slightly different situation (for examples based
on ``Bing doubling'').)

	Once one makes these reductions, one gets down to the case of
the single ball $B$ in $({\bf R}^3 / W) \times {\bf R}$, centered at
$(q,0)$ and with radius $1$.  For this single choice of scale and
location, one can use the fact that $({\bf R}^3 / W) \times {\bf R}$
is homeomorphic to ${\bf R}^4$ to get suitable local coordinates.

	For this single ball $B$, there is no issue of ``uniformity''
in the moduli of continuity for the coordinate mappings.  One simply
needs \emph{a} modulus of continuity.  

	A key point, however, is that when works backwards in the
reductions just made, to go to arbitrary balls in $({\bf R}^3 / W)
\times {\bf R}$ which are relatively close to the singular line $\{q\}
\times {\bf R}$, one does get local coordinates with uniform control
on the (normalized) moduli of continuity of the coordinate mappings
and their inverses.  This is because of the way that the reductions
cooperate with the scaling and the geometry.

	At any rate, this completes the outline of the argument for
showing that $({\bf R}^3 / W) \times {\bf R}$ has uniform local
coordinates, in the sense described before, around (\ref{image of beta
covers B}) and (\ref{uniform moduli of continuity, after rescaling}),
and with the kind of geometry for $({\bf R}^3 / W) \times {\bf R}$ as
before.  As mentioned at the beginning, this argument is nearly the
same as one given in \cite{Seqs}, in the proof of Theorem 6.3 on p241
of \cite{Seqs}.

	It is easy to see that a metric space which is bilipschitz
equivalent to some ${\bf R}^n$ has uniform local coordinates (relative
to ${\bf R}^n$ rather than ${\bf R}^4$, as above).  The required local
coordinates can simply be obtained from restrictions of the global
bilipschitz parameterization.  The converse is not true in general,
i.e., there are spaces which have uniform local coordinates and not be
bilipschitz equivalent to the corresponding ${\bf R}^n$.

	Examples of this are given by $({\bf R}^3 / W) \times {\bf R}$
with the kind of geometry as above.  We saw before that the singular
line $\{q\} \times {\bf R}$, which has Hausdorff dimension $1$ in our
geometry on $({\bf R}^3 / W) \times {\bf R}$, is always sent to sets
of Hausdorff dimension at least $2$ under any homeomorphism from
$({\bf R}^3 / W) \times {\bf R}$ onto ${\bf R}^4$, so that such a
homeomorphism can never be Lipschitz, or even H\"older continuous of
order greater than $1/2$.

	However, one can also get much simpler examples, by taking
snowflake spaces.  That is, one can take ${\bf R}^n$ equipped with the
metric $|x-y|^\alpha$, where $|x-y|$ is the usual metric, and $\alpha$
is a positive real number strictly less than $1$.  It is not hard to
check that this has the property of uniform local coordinates.  It is
not bilipschitz equivalent to ${\bf R}^n$, because it has Hausdorff
dimension $n/\alpha$ instead of Hausdorff dimension $n$.

	If one assumes uniform local coordinates relative to ${\bf
R}^n$ and Ahlfors-regularity of dimension $n$ (Definition \ref{ahlfors
regular}), then the matter becomes much more delicate.  Note that
Ahlfors-regularity of the correct dimension rules out snowflake
spaces, as in the previous paragraph.  It is still not sufficient,
though, because the $({\bf R}^3 / W) \times {\bf R}$ spaces with the
type of geometry as above are Ahlfors-regular of dimension $4$, in
addition to having uniform local coordinates (relative to ${\bf
R}^4$).

	With Ahlfors-regularity of dimension $n$ (and uniform local
coordinates relative to ${\bf R}^n$), there are some relatively strong
conclusions that one can make, even if one does not get something like
bilipschitz parameterizations.  Some versions of this came up in
Sections \ref{Quantitative topology, and calculus on singular spaces}
and \ref{Uniform rectifiability}, and with much less than uniform
local coordinates being needed.  We shall encounter another result
which is special to dimension $2$, later in this subsubsection.

	Dimension $1$ is more special in this regard (and as is quite
standard).  One can get bilipschitz equivalence with ${\bf R}$ from
Ahlfors-regularity of dimension $1$ and uniform local coordinates with
respect to ${\bf R}$ (and less than that).  This is not hard to do,
using arclength parameterizations.

	Let us note that in place of $({\bf R}^3 / W) \times {\bf R}$
in dimension $4$, as above, one can take Cartesian products with more
copies of ${\bf R}$ to get examples in all dimension greater than or
equal to $4$, with similar properties as in dimension $4$.  In
particular, one would get sets which are Ahlfors-regular of the
correct dimension, and which have uniform local coordinates, but not
bilipschitz coordinates, or even somewhat less than bilipschitz
coordinates.

	In dimension $3$, there is a different family of examples that
one can make, based on ``Bing doubling''.  This was mentioned earlier,
and is discussed in \cite{Seqs}.

	Notice that the uniform local coordinates property
\emph{would} imply a bilipschitz condition if the local coordinates
all came from restrictions of a single global parameterization.  In
general, the uniform local coordinates property allows the local
coordinate mappings to change as one changes locations and scales, and
this is why it allows for the possibility that there is no global
bilipschitz parameterization.  The case of $({\bf R}^3 / W) \times
{\bf R}$ provides a good example of this.  For this one can check the
earlier argument, and see how the local coordinates that are used are
not restrictions of a single global parameterization.  This is true
even though the local coordinate mappings at different locations and
scales often have a common ``model'' or behavior.  In other words, one
might often use \emph{rescalings} of a single mapping, and this can be
quite different from \emph{restrictions} of a single mapping.

	Let us emphasize that in the uniform local coordinates
condition, one does not ask that the uniform modulus of continuity be
something simple, like a Lipschitz or a H\"older condition.  Under
some conditions one might be able to \emph{derive} that, using the
uniformity of the hypothesis over all locations and scales.  This is
analogous to something that happens a lot in harmonic analysis; i.e.,
relatively weak conditions that hold uniformly over all locations and
scales often imply conditions which are apparently much stronger.  In
the present circumstances, we do have some limits on this, as in the
examples mentioned above.  The final answers are not clear, however.

	Another way to think about the uniformity over all locations
and scales in the uniform local coordinates property is that it is a
condition which implies the existence of homeomorphic coordinates even
after one ``blows up'' the space (in the Hausdorff or Gromov--Hausdorff
senses) along any sequence of locations and scales in the space.  With
the uniform local coordinates property, the local coordinates could be
``blown up'' along with the space, with the uniform bounds for the
moduli of continuity providing the equicontinuity and compactness
needed to take limits of the coordinate mappings (after passing to
suitable subsequences).

	Instead of looking at uniform local coordinates in connection
with \emph{bilipschitz} equivalence with ${\bf R}^n$, one can consider
\emph{quasisymmetric} equivalence.  Roughly speaking, a
\emph{quasisymmetric} mapping between two metric spaces is one that
approximately preserves \emph{relative} distances, in the same way
that bilipschitz mappings approximately preserve actual distances.  In
other words, if one has three points $x$, $y$, and $z$ in the domain
of such a mapping, and if $x$ is much closer to $y$ than $z$ is, then
this should also be true for their images under a
\emph{quasisymmetric} mapping, even if the actual distances between
the points might be changed a lot.  See \cite{TV} for more details and
information about quasisymmetric mappings.  Two metric spaces are
\emph{quasisymmetrically equivalent} if there is a quasisymmetric
mapping from one onto the other.  As with bilipschitz mappings,
compositions and inverses of quasisymmetric mappings remain
quasisymmetric.

	If a metric space admits a quasisymmetric parameterization
from ${\bf R}^n$, then it also satisfies the condition of uniform
local coordinates.  This is true for nearly the same reason as for
bilipschitz mappings; given a quasisymmetric mapping from ${\bf R}^n$
onto the metric space, one can get suitable local coordinates for the
space from restrictions of the global mapping to individual balls.
There is a difference between this case and that of bilipschitz
mappings, which is that one should allow some extra rescalings to
compensate for the fact that distances are not approximately
preserved.  Specifically, a ball $B$ in the metric space may be
covered in a nice way by the image under a quasisymmetric mapping of a
ball $\beta$ in ${\bf R}^n$, but then there is no reason why the radii
of $B$ and $\beta$ should be approximately the same.  For a
bilipschitz mapping, one would be able to choose $\beta$ so that it
has radius which is comparable to that of $B$.  In the quasisymmetric
case, one may not have that, but if one adds an extra rescaling on
${\bf R}^n$ (depending on the choices of balls), then one can still
get local coordinates with the kind of uniform control on the
(normalized) moduli of continuity as in the uniform local coordinates
condition.

	If a metric space has uniform local coordinates with respect
to ${\bf R}^n$, and if these coordinates come from restrictions and
then rescalings (on ${\bf R}^n$) of a single global parameterization,
as in the preceding paragraph, then that parameterization does have to
be quasisymmetric.  This is an easy consequence of the definitions,
and it is analogous to what happens in the bilipschitz case.

	If a metric space admits uniform local coordinates from some
${\bf R}^n$, it still may not be true that it admits a quasisymmetric
parameterization.  This is trickier than before, and in particular one
does not get examples simply by using snowflake metrics $|x-y|^\alpha$
on ${\bf R}^n$.  Indeed, the identity mapping on ${\bf R}^n$ \emph{is}
quasisymmetric as a mapping from ${\bf R}^n$ with the standard metric
to ${\bf R}^n$ with the snowflake metric $|x-y|^\alpha$, $0 < \alpha <
1$.

	However, there are counterexamples, going back to results of
Rickman and V\"ais\"al\"a.  That is, these are spaces which have
uniform local coordinates (and are even somewhat nicer than that), but
which do not admit quasisymmetric parameterizations.  Basically, these
spaces are Cartesian products, where the individual factors behave
nicely in their own right, and where the combination mixes different
types of geometry.  A basic example (which was the original one) is to
take a product of a snowflake with a straight line.  Quasisymmetric
mappings try to treat different directions in a uniform manner, and in
the end this does not work for parameterizations of these examples.
See Lemma 4 in \cite{Tu}, and also \cite{V-prod} and \cite{AV-prod}.

	These examples occur already in dimension $2$.  They do not
behave well in terms of measure, though, and for Ahlfors-regularity
(of dimension $2$) in particular.  This is a basic part of the story;
compare with \cite{AV-prod, Tu, V-prod}.

	As usual, dimension $1$ is special.  There are results
starting from more primitive conditions, and good characterizations
for the existence of quasisymmetric parameterizations, in fact.  See
Section 4 of \cite{TV}.

	On the other hand, in dimension $2$, there are positive
results about having global quasisymmetric parameterizations for a
given space, and with bounds.  For these results one would assume
Ahlfors-regularity of dimension $2$, and some modest additional
conditions concerning the geometry and topology, conditions which are
weaker than ``uniform local coordinates''.  See \cite{DS-trans, HK1,
Se-casscII}.  For these results, it is important that the dimension be
$2$, and not larger, because of the way that they rely on the
existence of conformal mappings.  We shall say a bit more about this
in Subsection \ref{Geometric and analytic results about the existence
of good coordinates}.

	If we go now to dimension $3$, then there are counterexamples,
even with good behavior in terms of measure.  That is, there are
examples of subsets of Euclidean spaces which are Ahlfors-regular of
dimension $3$, admit uniform local coordinates with respect to ${\bf
R}^3$, but are not quasisymmetrically-equivalent to ${\bf R}^3$.
These examples are based on the decompositions of ${\bf R}^3$
associated to ``Bing doubling'' (as in \cite{Bing1, Bing6, Dm2}),
using geometric realizations as in Subsubsection \ref{A basic
scenario}.  This is discussed in \cite{Seqs}.  The absence of a
quasisymmetric parameterization in this case is close to a result in
\cite{FS}, although the general setting in \cite{FS} is a bit
different.

	Concerning the space $({\bf R}^3 / W) \times {\bf R}$
considered before, equipped with a nice geometry as from Subsubsection
\ref{A basic scenario} and \cite{Seqs}, it is not clear (to my
knowledge) whether \emph{quasisymmetric} parameterizations from ${\bf
R}^4$ exist or not.  We saw earlier that \emph{bilipschitz} mappings
do not exist, because of the line in $({\bf R}^3 / W) \times {\bf R}$
which has to have Hausdorff dimension at least $2$ after any
homeomorphism from $({\bf R}^3 / W) \times {\bf R}$ onto ${\bf R}^4$.
These considerations of Hausdorff measure do not by themselves rule
out the existence of a quasisymmetric mapping, as they do for
bilipschitz mappings.  (Compare with \cite{V1}, for instance.)

	Similar remarks apply to double-suspensions of homology
spheres.  In particular, it is not known (to my knowledge) whether or
not quasisymmetric parameterizations exist for them.

	To summarize a bit, we have seen in this subsubsection how
decomposition spaces might be ``standard'' topologically, being
equivalent to ${\bf R}^n$ as a topological space, and still lead to
geometries which are tricky, and which have some interesting
structure.  We shall discuss another class of examples like this in
Subsubsection \ref{Examples that are even simpler topologically, ...}.
These examples are more trivial topologically, more tame
geometrically, but still nontrivial geometrically.  (They admit simple
quasisymmetric parameterizations, but no bilipschitz
parameterizations.)

	See \cite{Seqs} for some other types of geometric structure
which can arise naturally from decompositions.  For instance, cellular
decompositions whose quotients are not manifolds --- such as Bing's
dogbone space (\cite{Bing2}, Example 4 on p64-65 in \cite{Dm2}) ---
give rise to spaces which do not have ``uniform local coordinates'',
but which do have nice properties in other ways.  For instance, one
can make constructions so that local coordinates exist at all
locations and scales, and with uniform bounds on how they are
localized, but not for their moduli of continuity.  This last is
related to having Gromov--Hausdorff limits which are not manifolds, but
are topologically given as a decomposition space (with a cellular
decomposition), like Bing's dogbone space.

\subsubsection{Examples that are even simpler topologically, but still
nontrivial geometrically}
\label{Examples that are even simpler topologically, ...}

	Let us mention another class of examples, which one can also
think of in terms of decompositions (although they are ``trivial'' in
this respect).  These examples are based on ``Antoine's necklaces'',
which came up before, in Subsection \ref{Wildness and tameness
phenomena}.

	Antoine's necklaces are compact subsets of ${\bf R}^3$ which
are homeomorphic to the usual middle-thirds Cantor set in the real
line, but for which there is no global homeomorphism from ${\bf R}^3$
onto itself which maps these sets into subsets of a line.  In
dimension $2$ this does not happen, as in Chapter 13 of \cite{Moise}.

	The ``wildness'' of these sets is manifested in a simple
fundamental-group property.  Namely, the complement of these sets in
${\bf R}^3$ have nontrivial fundamental group, whereas this would not
be true if there were a global homeomorphism from ${\bf R}^3$ to
itself which would take one of these sets to a subset of the line.
This last uses the fact that these sets are totally disconnected
(i.e., to have simple-connectivity of the complement in ${\bf R}^3$ if
the set were to lie in a line).

	Antoine's necklaces are discussed in Chapter 18 of
\cite{Moise}.  See also p71ff in \cite{Dm2}.  The basic construction
for them can be described in terms of the same kind of ``rules'' as in
Subsubsection \ref{A basic scenario}.

	One starts with a solid torus $T$ in ${\bf R}^3$.  Inside this
torus one embeds some more tori, which are disjoint, but which form a
chain that is ``linked'' around the hole in the original torus.  (See
Figure 18.1 on p127 of \cite{Moise}, or Figure 9.9 on p71 of
\cite{Dm2}.)  In each of these smaller tori, one can embed another
collection of linking tori, in the same way as for the first solid
torus.

	One can repeat this indefinitely.  In the limit, one gets a Cantor
set, which is a necklace of Antoine.  

	Actually, we should be a little more precise here.  In saying
that we get a Cantor set in the limit, we are implicitly imagining
that the diameters of the solid tori are going to $0$ as one proceeds
through the generations of the construction.  This is easy to arrange,
if one uses enough tori in the basic rule (linking around the original
torus $T$).  If one uses enough tori, then one can do this in such a
way that all of the tori are \emph{similar} to the initial torus $T$,
i.e., can be given as images of $T$ by mappings which are combinations
of translations, rotations, and dilations.

	One can also consider using a smaller number of tori, and
where the embeddings of the new tori (in the original one) are allowed
to have some stretching.  (Compare with Figure 9.9 on p71 of
\cite{Dm2}.)  As an extreme case, the Whitehead continuum corresponds
essentially to the same construction as Antoine's necklaces, except
that only $1$ embedded torus is used in the linking in the original
torus $T$.  (See Figure 9.7 on p68 of \cite{Dm2}.)  For this one
definitely needs a fair amount of stretching.  In Bing doubling, one
uses two solid tori, embedded and linked around the hole in $T$
(Figure 3 on p357 of \cite{Bing1}, Figure 9.1 on p63 in \cite{Dm2}).
One again needs some stretching, but not as much for the Whitehead
continuum.

	These cases are different from the standard ones for Antoine's
necklaces, because when one iterates the basic rule in a
straightforward way, the components that one gets at the $n$th
generation do not have diameters tending to $0$.  For the Whitehead
continuum, this is simply unavoidable, and reflects the way that the
components are clasped, each one by itself, around the ``hole'' of the
original torus $T$.  In the case of Bing doubling, one can rearrange
the embeddings at later generations in such a way that the diameters
do tend to $0$, even if this might not be true for naive iterations.
This was proved by Bing, in \cite{Bing1}.  (See also p69-70 of
\cite{Dm2}, and \cite{Bing6}.)

	Let us imagine that we are using enough small solid tori in
the linking around $T$, as in standard constructions for Antoine's
necklaces, so that it is clear that the diameters of the components of
the sets obtained by repeating the process do go to $0$ (i.e., without
having to make special rearrangements, or anything like that, as in
Bing doubling).  In other words, in the limit, one gets a Cantor set
in ${\bf R}^3$, as above.  Let us call this Cantor set $A$.

	This Cantor set $A$ is wild, in the sense that its complement
in ${\bf R}^3$ is not simply-connected, and there is no global
homeomorphism from ${\bf R}^3$ to itself which takes $A$ to a subset
of a line.  At the level of decompositions, however, there is not much
going on here.

	Normally, to get a decomposition from a process like this, one
takes the connected components obtained in the limit of the process,
together with sets with one element for the rest of the points in
${\bf R}^3$ (or ${\bf R}^n$, as the case may be).  This was described
before, in Subsubsection \ref{A basic scenario}.  (See also the
discussion at the beginning of Section 9 of \cite{Dm2}, on p61,
concerning the notion of a \emph{defining sequence}.)  In the
construction that we are considering here, all of the connected
components in the end contain only one element each, because of the
way that the diameters of the components in finite stages of the
``defining sequence'' converge to $0$.

	In other words, the decomposition that occurs here is
automatically ``trivial'', consisting of one-element sets, one for
each point in ${\bf R}^3$.  Taking the quotient does not do anything,
and the decomposition space is just ${\bf R}^3$ again.

	There is nothing too complicated about this.  It is just
something to say explicitly, for the record, so to speak, to be clear
about it, especially since it is a situation to which one might
normally pay little attention, for being degenerate.  (See also the
text at the beginning of p71 of \cite{Dm2}, about this kind of
defining sequence and decomposition.)  

	While there is nothing going on at the level of the
decompositions \emph{topologically}, this is not the case
geometrically!  One can think of this in the same way as in
Subsubsection \ref{A basic scenario}, for making geometric
representations of decomposition spaces, with metrics and
self-similarity properties for them.  In the present situation, one
can also work more directly at the level of ${\bf R}^3$ itself, to get
geometries like this.

	Here is the basic point.  Imagine deforming the geometry of
${\bf R}^3$, at the level of infinitesimal measurements of distance,
as with Riemannian metrics.  In the general idea of a decomposition
space, one can shrink sets in some ${\bf R}^n$ which have nonzero
diameter to single points.  In the present setting, our basic
components already \emph{are} single points, and so there is nothing
to do to them.  However, one can still shrink the \emph{geometry}
around these points in ${\bf R}^3$.

	In technical terms, one can think of deforming the geometry of
${\bf R}^3$, by multiplying its standard Riemannian metric by a
function.  One can take this function to be positive and regular away
from the Antoine's necklace, and then vanish on the necklace itself.
For instance, one could take the function to be a positive power of
the distance to the necklace $A$, so that the Riemannian metric can be
written as follows:
\begin{equation}
\label{Riemannian metric, deformed along the Antoine's necklace A}
	ds^2 = \dist(x,A)^\alpha \, dx^2.
\end{equation}

	This type of metric is discussed in some detail in
\cite{Sebil} (although in a slightly different form).  The geometry of
${\bf R}^3$ with this kind of metric behaves a lot like standard
Euclidean geometry.  One still has basic properties like
Ahlfors-regularity of dimension $3$, and Sobolev, Poincar\'e, and
isoperimetric inequalities in the new geometry.  See \cite{Sebil}.
(This kind of deformation is special case of the broader class in
\cite{DS-sob} and \cite{Se-finn, Sebarc}.)

	However, with suitable choices of parameters, the metric space
that one gets in this way is \emph{not} bilipschitz equivalent to
${\bf R}^3$ with the standard Euclidean metric.  This is because the
shrinking of distances around the necklace $A$ can lead to $A$ having
Hausdorff dimension less than $1$ in the new metric.  On the other
hand, ${\bf R}^3 \backslash A$ is not simply-connected, and this means
that the geometry which has been constructed cannot be bilipschitz
equivalent to the standard geometry on ${\bf R}^3$, because of Theorem
\ref{hausdorff dimension and local fundamental groups} in Subsection
\ref{Manifold factors}.  See \cite{Sebil} for more information.
(Concerning the ``choice of parameters'' here, the main point is to
have enough shrinking of distances around $A$ to get the Hausdorff
dimension to be less than $1$, or something like that.  The amount of
shrinking needed depends on some of the choices involved in producing
the necklace.  By using Antoine's necklaces which are sufficiently
``thin'', it is enough to employ arbitrarily small powers $\alpha$ of
the distance to the necklace in (\ref{Riemannian metric, deformed
along the Antoine's necklace A}) to get enough shrinking of the metric
around the necklace.  In any case, one is always free to take the
power $\alpha$ to be larger.)

	Let us emphasize that in making this kind of construction, the
conclusion is that the metric space that one gets is not
bilipschitz-equivalent to ${\bf R}^3$ through \emph{any} homeomorphism
between the two spaces.  It is easy to make deformations of the
geometry so that the new metric seems to be much different from the
old one in the given coordinates, but for which this is not really the
case if one is allowed to make a change of variables.  For instance,
one might deform the standard Euclidean Riemannian metric on ${\bf
R}^3$ by multiplying it by a function that vanishes at a point, like a
positive power of the distance to that point.  Explicitly, this means
\begin{equation}
\label{deforming the Riemannian metric around a point p}
	d\widetilde{s}^2 = |x-p|^\beta \, dx^2,
\end{equation}
where $\beta > 0$.  In the standard coordinates, this metric and the
ordinary one look quite different.  However, for this particular type
of deformation (as in (\ref{deforming the Riemannian metric around a
point p})), the two metrics are bilipschitz equivalent, if one allows
a nontrivial change of variables.  Specifically, one can use changes
of variables of the form
\begin{equation}
\label{radial change of variables}
	f(x) = p + |x-p|^{\beta/2} (x-p).
\end{equation}
This is not hard to verify.

	As a more complicated version of this, one can also make
deformations of the standard geometry on ${\bf R}^3$ of the form
\begin{equation}
	d\widehat{s}^2 = \dist(x, K)^\gamma \, dx^2,
\end{equation}
where $\gamma > 0$ and $K$ is a self-similar Cantor set in ${\bf R}^3$
which is \emph{not} wild.  In this case one can again get geometries
which may look different from the usual one in the given coordinates,
but for which there are changes of variables which give a bilipschitz
equivalence with the standard metric.  Compare with Remark 5.28 on
p390 of \cite{Sebil}.

	When one makes deformations based on Antoine's necklaces, as
above, the linking that goes on can ensure that there is no
bilipschitz equivalence between ${\bf R}^3$ with the new geometry and
${\bf R}^3$ with the standard geometry.  In fact, there will not be a
homeomorphism which is Lipschitz (from the new geometry to the
standard one), without asking for bilipschitzness.  With suitable
choices of parameters, it can be impossible to have a homeomorphism
like this which is even H\"older continuous of an arbitrary exponent
$\delta > 0$, given in advance.

	Under the conditions in \cite{Sebil}, the identity mapping
itself on ${\bf R}^3$ always gives a homeomorphism which is H\"older
continuous with \emph{some} positive exponent.  In fact, it is also
\emph{quasisymmetric}, in the sense of \cite{TV}.  Thus, here one gets
examples of spaces which are quasisymmetrically equivalent to a
Euclidean space, and which are Ahlfors-regular of the correct
dimension, but which are not bilipschitz equivalent to a Euclidean
space.  (Compare with Subsubsection \ref{Quotient spaces can be
topologically standard, but geometrically tricky}.)

\subsection{Geometric and analytic results about the existence of good 
coordinates}
\label{Geometric and analytic results about the existence of good coordinates}

	In Subsection \ref{Contractable open sets}, we considered the
question of whether a nonempty contractable open set in ${\bf R}^n$ is
homeomorphic to the standard open unit ball in ${\bf R}^n$.  When $n = 2$
this is true, and it is a standard result in topology.

	One can establish this result in $2$ dimensions analytically
via the \emph{Riemann Mapping Theorem}.  This theorem gives the
existence of a \emph{conformal} mapping from the unit disk in ${\bf
R}^2$ onto any nonempty simply-connected open set in ${\bf R}^2$ which
is not all of ${\bf R}^2$.  See Chapter 6 of \cite{Ahlfors1}, for
instance.

	The Riemann Mapping Theorem is of course very important for
many aspects of analysis and geometric function theory in ${\bf R}^2
\cong {\bf C}$, but it also does a lot at a less special level.  From
it one not only obtains homeomorphisms from the unit disk in ${\bf
R}^2$ onto any nonempty simply-connected proper open subset of ${\bf
R}^2$, but one gets a way of choosing such homeomorphisms which is
fairly canonical.  In particular, Riemann mappings are unique modulo a
three real-dimensional group of automorphisms of the unit disk (which
can be avoided through suitable normalizations), and there are results
about the dependence of Riemann mappings on the domains being
parameterized.  

	By comparison, one might try to imagine doing such things
without the Riemann mapping, or in other contexts where it is not
available.  In this regard, see \cite{Hatcher1, Hatcher2, Laudenbach,
Randall-Schweitzer}, concerning related matters in higher dimensions.

	Another fact in dimension $2$ is that any smooth Riemannian
metric on the $2$-sphere ${\bf S}^2$ is conformally-equivalent to the
standard metric.  In other words, if $g$ is a smooth Riemannian metric
on ${\bf S}^2$, and if $g_0$ denotes the standard metric, then there
is a diffeomorphism from ${\bf S}^2$ onto ${\bf S}^2$ which converts
$g$ into a metric of the form $\lambda \, g_0$, where $\lambda$ is a
smooth positive function on ${\bf S}^2$.  There are also local and
other versions of this fact, but for the moment let us stick to this
formulation.

	One way to try to use this theorem is as follows.  Suppose
that one has a $2$-dimensional space which behaves roughly like a
$2$-dimensional Euclidean space (or sphere) in some ways, and one
would like to know whether it can be realized as nearly-Euclidean in
more definite ways, through a parameterization which respects the
geometry.  Let us assume for simplicity that our space is given to us
as ${\bf S}^2$ with a smooth Riemannian metric $g$, but without bounds
for the smoothness of $g$.  One can then get a conformal
diffeomorphism $f : ({\bf S}^2, g_0) \to ({\bf S}^2, g)$, as in the
result mentioned in the preceding paragraph.  A priori the behavior of
this mapping could be pretty complicated, and one might not know much
about it at definite scales.  It would be nice to have some bounds for
the behavior of $f$, in terms of simple geometric properties of $({\bf
S}^2, g)$.

	Some results of this type are given in \cite{DS-trans, HK1,
Se-casscII}.  Specifically, general conditions are given in
\cite{DS-trans, HK1, Se-casscII} under which a conformal equivalence
$f : ({\bf S}^2, g_0) \to ({\bf S}^2, g)$ actually gives a
\emph{quasisymmetric} mapping (as in \cite{TV}), with uniform bounds
for the quasisymmetry condition.  In other words, these results have
the effect of giving uniform bounds for the behavior of $f$ at any
location or scale, under suitable conditions on the initial space
$({\bf S}^2, g)$, and using the conformality of $f$.

	Another use of conformal mappings of an analogous nature is
given in \cite{MuS}.  There the assumptions on the space involved more
smoothness --- an integral condition on (principal) curvatures, for a
surface in some ${\bf R}^n$ --- and the conclusions are also stronger,
concerning bilipschitz coordinates.  This gave a new approach to
results in \cite{T1}.  See also \cite{T2}, and the recent and quite
different method in \cite{Fu-J}.

	More precisely, \cite{MuS} works with conformal mappings,
while \cite{T1, T2} and \cite{Fu-J} obtain bilipschitz coordinates by
quite different means.  In the context of \cite{DS-trans, HK1,
Se-casscII}, no other method for getting quasisymmetric or other
coordinates with geometric estimates (under similar conditions) is
known, at least to my knowledge.  

	One might keep in mind that conformality is defined in
infinitesimal terms, through the differential of $f$.  To go from
infinitesimal or very small-scale behavior to estimates at larger
scales, one in effect tries to ``integrate'' the information that one
has.

	This is a very classical subject for conformal and
quasiconformal mappings.  A priori, it is rather tricky, because one
is not given any information about the conformal factors (like the
function $\lambda$ before).  Thus one cannot ``integrate'' directly in
a conventional sense.  One of the basic methods is that of ``extremal
length'', which deals with the balance between length and area.  At
any rate, methods like these are highly nonlinear.  See
\cite{Ahlfors2, Ahlfors3, Lehto-Virtanen, V-book} for more
information.

	What would happen if one attempted analogous enterprises in
higher dimensions?  One can begin in the same manner as before.  Let
$n$ be an integer greater than or equal to $2$, and suppose (as a
basic scenario) that one has a smooth Riemannian metric $g$ on ${\bf
S}^n$.  Let $g_0$ denote the standard metric on ${\bf S}^n$.  One
might like to know that $({\bf S}^n, g)$ can be parameterized by
$({\bf S}^n, g_0)$ through a mapping with reasonable properties, and
with suitable bounds, under some (hopefully modest) geometric
conditions on $({\bf S}^n, g)$.  Here, as before, the smoothness of
$g$ should be taken in the character of an a priori assumption.  One
would seek uniform bounds that do not depend on this in a quantitative
way.  (The bounds would depend on constants in the geometric
conditions on $({\bf S}^n, g)$.)

	If one has a mapping $f : ({\bf S}^n, g_0) \to ({\bf S}^n, g)$
which is conformal, or which is quasiconformal (with a bound for its
dilatation), then \cite{HK1} provides some natural hypotheses on
$({\bf S}^n, g)$ under which one can establish that $f$ is
quasisymmetric, and with bounds.  In other words, this works for all
$n \ge 2$, and not just $n = 2$, as above.  See \cite{HK2, HK3} for
further results along these lines.

	However, when $n > 2$, there are no general results about
\emph{existence} of conformal parameterizations for a given space, or
quasiconformal parameterizations with uniform bounds for the
dilatation.  It is simply not true that arbitrary smooth Riemannian
metrics admit conformal coordinates, even locally, as they do when $n
= 2$.  Quasiconformal coordinates automatically exist for
reasonably-nice metrics, but with the quasiconformal dilatation
depending on the metric or on the size of the region being
parameterized in a strong way.  The issue would be to avoid or reduce
that.

	One can easily see that the problem is highly overdetermined,
in the following sense.  A general Riemannian metric in $n$ dimensions
is described (locally, say) by $n(n+1)/2$ real-valued functions of $n$
variables.  A conformal deformation of the standard metric is defined
by $1$ real-valued function of $n$ real variables, i.e., for the
conformal factor.  A general diffeomorphism in $n$ dimensions is
described by $n$ real-valued functions of $n$ variables.  Thus,
allowing for general changes of variables, the metrics which are
conformally-equivalent to the standard metric are described by $n+1$
real-valued functions of $n$ real variables.  When $n = 2$, this is
equal to $n (n+1) / 2$, but for $n > 2$ one has that $n (n+1)/2 > n +
1$.

	In fact, one knows that in dimension $3$ there are numerous
examples of spaces which satisfy geometric conditions analogous to
ones that work in dimension $2$, but which do not admit quasisymmetric
parameterizations.  There are also different levels of structure which
occur in dimension $3$, between basic geometric properties and having
quasisymmetric parameterizations, and which would come together in
dimension $2$.  See \cite{Seqs}.  Parts of this are reviewed or
discussed in Subsection \ref{Geometric structures for decomposition
spaces}, especially Subsubsection \ref{Quotient spaces can be
topologically standard, but geometrically tricky}.

	Thus, not only does the method based on conformal and
quasiconformal mappings not work in higher dimensions, but some of the
basic results that one might hope to get or expect simply are not
true, by examples which are pretty concrete.

	These examples can be viewed as geometrizations of classical
examples from geometric topology, and they are based on practically
the same principles.  They are not especially strange or pathological
or anything like that, but have a lot of nice properties.  They
reflect basic phenomena that occur.

	This is all pretty neat!  One has kinds of ``parallel
tracks'', with geometric topology on one side, and aspects of geometry
and analysis on the other.  A priori, these two tracks can exist
independently, even if there are ways in which each can be involved in
the other.

	Each of these two tracks has special features in low
dimensions.  This concerns the existence of homeomorphisms with
certain properties, for instance.  Each has statements and results and
machinery which make sense for the given track, and not for the other
side, even if there are also some overlaps (as with applications of
Riemann mappings).

	Each of these two tracks also starts running into trouble in
higher dimensions, and at about the same time!  We have seen a number
of instances of this by now, in this subsection, in this appendix more
broadly, and also Section \ref{parameterization problems}.  The kinds
of trouble that they encounter can be rather different a priori (such
as localized fundamental group conditions, versus behavior of partial
differential equations), even if there is again significant overlap
between them.  In particular, this concerns the existence of
homeomorphisms with good properties.

\subsubsection{Special coordinates that one might consider in other 
dimensions}
\label{Special coordinates that one might consider in other dimensions}

	We have already discussed a number of basic topological
phenomena in this appendix.  Let us now briefly consider a couple of
things that one might try in higher dimensions on the side of geometry
and analysis, in similar veins as above.

	One basic approach would be to try to find and use mappings
which minimize some kind of ``energy''.  As before, one can consider
smooth metrics on smooth manifolds (like ${\bf S}^n$), and try to get
parameterizations with uniform bounds on their behavior, under modest
conditions on the geometry of the spaces.  (One can also try to work
directly with spaces and metric that are not smooth.)

	A very standard energy functional to consider would be the
$L^2$ norm of the differential, as with harmonic mappings.  In
dimension $2$, conformal mappings can be placed in this framework.
One can also consider energy functionals based on $L^p$ norms of
differentials of mappings.  This is more complicated in terms of the
differential equations that come up, but it can have other advantages.
The choice of $p$ as the dimension $n$ has some particularly nice
features, for having the energy functional cooperate with the geometry
(and analysis).  (This is one of the ways that $n = 2$ is special; for
this one can have both $p = n$ and $p = 2$ at the same time!)  In
particular, the energy becomes invariant under conformal changes in
the metric when $p$ is equal to the dimension.

	In elasticity theory, one considers more elaborate energy
functionals as well.  For instance, these might include integral norms
of the \emph{inverse} of the Jacobian of the mapping, in addition to
$L^p$ norms of the differentials.  In other words, the functional can
try to limit both the way that the differential becomes large and
small, so that it takes into account both stretchings and
compressions.

	In any case, although there is a lot of work concerning
existence and behavior of minimizers for functionals like these, I do
not really know of results in dimensions $n \ge 3$ where they can be
used to obtain well-behaved parameterizations of spaces, with bounds,
under modest or general geometric conditions.  This is especially true
in comparison with what one can get in dimension $2$, as discussed
before.

	This is a bit unfortunate, compared with the way that the
normal form of conformal mappings can be so useful in dimension $2$.
On the other hand, perhaps someone will find ways of using such
variational problems for geometric questions like these some day, or
will find some kind of special structure connected to them.  In this
regard, one might bear in mind issues related to mappings with
\emph{branching}, as in Appendix \ref{branching}.  We shall say a bit
more about this later in this subsection.

	One might also keep in mind the existence of spaces with good
properties, but not good parameterizations, as mentioned earlier (and
discussed in Subsubsections \ref{Quotient spaces can be topologically
standard, but geometrically tricky} and \ref{Examples that are even
simpler topologically, ...}, and in \cite{Seqs, Sebil}).

	In dimension $3$, there is another kind of special structure
that one might consider.  Namely, instead of metrics which are
conformal deformations of the standard Euclidean metric, let us
consider metrics $g = g_{i,j}$ for which only the \emph{diagonal}
entries $g_{i,i}$ are nonzero.

	In this case the diagonal entries are allowed to vary
independently.  For conformal deformations of the standard Euclidean
metric, the off-diagonal entries are zero, and the diagonal entries
are all equal.

	In dimension $3$, the problem of making a change of variables
to put a given metric into diagonal form like this is ``determined'',
in the same way as for conformal deformations of Euclidean metrics in
dimension $2$.  Specifically, one can compute as follows.  A general
Riemannian metric is described by $n (n+1) / 2$ real-valued functions
of $n$ variables, which means $6$ real-valued functions of $3$ real
variables in dimension $3$.  Metrics with only diagonal nonzero
entries are defined by $3$ real-valued functions of $3$-real
variables, and changes of variables are given by $3$ real-valued
functions of $3$ real variables as well.  Thus, allowing for changes
of variables, the metrics that can be reduced to diagonal metrics can
be described by $6$ real-valued functions of $3$ real-variables, which
is the same as for the total family of Riemannian metrics in this
dimension.  In dimensions greater than or equal to $4$, this would not
work, and there would again be too many Riemannian metrics in general
compared to diagonal metrics and ways of reducing to them via changes
of variables.

	Of course this is just an informal ``dimension'' count, and
not a justification for being able to put metrics into diagonal form
in dimension $3$.  (One should also be careful that there is not
significant overlap between changes of variables and diagonal metrics,
i.e., so that there was no ``overcounting'' for the combination of
them.)  However, it does turn out that one can put metrics in diagonal
form (in dimension $3$), at least locally.  This was established in
\cite{deturck-yang} in the case of smooth metrics.  There were earlier
results in the real-analytic category.  (See \cite{deturck-yang} for
more information.)

	However, this type of ``normal form'' does not seem to be as
useful for the present type of issue as conformal parameterizations
are.  As in the case of conformal coordinates, part of the problem is
that even if one has such a normal form, one does not a priori know
anything about the behavior of the diagonal entries of the metric in
this normal form.  One would need methods of getting estimates without
this information, and only the nature of the normal form.  In the
context of conformal mappings, one has extremal lengths, conformal
capacities, and other conformal and quasiconformal invariants and
quasi-invariants.  For diagonal metrics, it is not clear what one
might do.

	A related point is that the analysis of the partial
differential equations which permits one to put smooth Riemannian
metrics in dimension $3$ into diagonal form is roughly ``hyperbolic'',
in the same way that the corresponding differential equations for
conformal coordinates in dimension $2$ are elliptic.  See
\cite{deturck-yang}.  This is closely connected to the kind of
stability that one has for conformal mappings, and the possibilities
for having estimates for them under mild or primitive geometric
conditions.

	In a way this is all ``just fair'', and nicely so.  With
diagonal metrics one does have something analogous to conformal
coordinates in dimension $3$.  On the other hand, this analogue
behaves differently in fundamental ways, including estimates.  This is
compatible with other aspects of the story as a whole, like the
topological and geometric examples that one has in dimension $3$
(where homeomorphisms may not exist, with the properties that one
might otherwise hope for).

	In any event, this illustrates how analytic and geometric
methods seem to behave rather differently in dimensions $3$ and
higher, compared to the special structure and phenomena which occur in
dimension $2$.  This is somewhat remarkable in analogy with
topological phenomena, which have similar differences between
dimensions.  With the topology there are both some crossings and
overlaps with geometry and analysis, and much that is separate or
independent.

	On the side of geometry and analysis, let us also note that
there are some other special features in low dimensions that we have
not mentioned.  As a basic example, the large amount of flexibility
that one has in making conformal mappings in dimension $2$ leads to
some possibilities in dimension $3$ that are not available in higher
dimensions.  That is, the large freedom that one has in dimension $2$
can sometimes permit one to make more limited constructions in
dimension $3$, e.g., by starting with submanifolds of dimension $2$,
and working from there (with extensions, gluings, etc.)  These
possibilities in dimension $3$ can be much more restricted than in
dimension $2$, but having them at all can be significantly more than
what happens in higher dimensions.

	We should also make clear that if one allows mappings with
\emph{branching}, as in Appendix \ref{branching}, then a number of
things can change.  Some topological difficulties could go away or be
ameliorated, as has been indicated before (and in Appendix
\ref{branching}).  For instance, the branching can unwind obstructions
or problems with localized fundamental groups (in complements of
points or other sets).  There are many basic examples of this, as in
Appendix \ref{branching}, and the constructions in \cite{HR1, HR2}.
Ideas of Sullivan \cite{Su2, Su3} are also important in this regard
(and as mentioned in Appendix \ref{branching}).

	General pictures for mappings with branching, including
existence and good behavior, have yet to be fully explored or
understood.  The Alexander argument described in Appendix
\ref{branching}, the constructions of Heinonen and Rickman \cite{HR1,
HR2}, classical work on quasiregular mappings (as in \cite{Res, Ri1}),
and the work of Sullivan \cite{Su2, Su3}, seem to indicate many
promising possibilities and directions.

	One can perhaps use variational problems in these regards
as well.

	Concerning variational problems, one might also keep in mind
the approaches of \cite{DS8, DS-q-min} (and some earlier ideas of
Morel and Solimini \cite{morel-solimini}).  For these one does not
necessarily work directly with mappings or potential parameterizations
of sets, and in particular one may allow sets themselves to be
variables in the minimization (rather than mappings between fixed
spaces).  This broader range can make it easier for the minimizations
to lead to useful conclusions about geometric structure and
complexity, under natural and modest conditions.  In particular, one
can get substantial ``partial parameterizations'', as with uniform
rectifiability conditions.

	These approaches are also nicely compatible with the trouble
that one knows can occur, related to topology and homeomorphisms (and
in geometrically moderate situations, as in Subsections \ref{Geometric
structures for decomposition spaces} and \ref{Examples that are even
simpler topologically, ...}, and \cite{Seqs, Sebil}).

	Finally, while we have mentioned a lot about the special
phenomena that can occur in dimension $2$, and what happens in higher
dimensions, we should also not forget about dimension $1$.  This is
even more special than dimension $2$.  This is a familiar theme in
geometric topology, for the ways that one can recognize and
parameterize curves.  In geometry and analysis, one can look for
parameterizations with bounds, and these are often constructible.

	A fundamental point along the lines is the ability to make
parameterizations by arclength, for curves of locally finite length.
More generally, one can use parameterizations adapted to other
measures (rather than length), when they are around.

	Arclength parameterizations provide a very robust and useful
way for obtaining parameterizations in dimension $1$ with good
behavior and bounds.  In dimension $1$, simple conditions in terms of
mass can often be immediately ``integrated'' to get well-behaved
parameterizations, in ways that are not available (or do not work
nearly as well) in higher dimensions, even in dimension $2$.

	To put the matter in more concrete terms, in dimension $1$ one
can often make parameterizations, or approximate parameterizations,
simply by ordering points in a good way.  This does not work in higher
dimensions.  Once one has the ordering, one can regularize the
geometry by parameterizing according to arclength, or some other
measure (as appropriate).

	For another version of this, in connection with quasisymmetric
mappings, see Section 4 of \cite{TV}.

	In differential-geometric language, one might say that
dimension $1$ is special for the way that one can make isometries
between spaces, through arclength parameterizations.  This no longer
works in dimension $2$, but one has conformal coordinates there.
Neither of these are generally available in higher dimensions.  In
higher dimensions one has less special structure for getting the
existence in general of well-behaved parameterizations, and then the
kinds and ranges of geometric and topological phenomena which can
exist open up in a large way.

\subsection{Nonlinear similarity: Another class of examples}
\label{Nonlinear similarity: Another class of examples}

	A very nice and concrete situation in which issues of
existence and behavior of homeomorphisms can come up is that of
``nonlinear similarity''.  Specifically, it is possible to have linear
mappings $A$, $B$ on ${\bf R}^n$ which are conjugate to each other by
homeomorphisms from ${\bf R}^n$ onto itself --- i.e., $B = h \circ A
\circ h^{-1}$, where $h$ is a homeomorphism of ${\bf R}^n$ onto itself
--- and which are \emph{not} conjugate by \emph{linear} mappings!

	Examples of this were given in \cite{C-S1}.  For related
matters, including other examples and conditions under which one can
deduce linear equivalence, see \cite{C-S-bulletin, C-S2, C-S3, C-S4,
C-S-S-W, C-S-S-W-W, HamPed1, HamPed2, HP1, HP2, KR, M-R1, M-R2, Mio,
deRham, deRham2, RW, Weinberger, Weinberger-book, Wilker}.

	Note that if one has a conjugation of linear mappings $A$ and
$B$ by a \emph{diffeomorphism} $h$ on ${\bf R}^n$, then one can derive
the existence of a linear conjugation from this.  This comes from
passing to the differential of the diffeomorphism at the origin.

	Thus, when a linear conjugation does not exist, but a
homeomorphic conjugation does, then the homeomorphism cannot be smooth
at all points in ${\bf R}^n$, or even at just the origin.  One might
wonder then about the kinds of processes and regularity that might be
entailed in the homeomorphisms that provide the conjugation.  In this
regard, see \cite{C-S3, RW, Weinberger}.

\section{Doing pretty well with spaces which may not have nice coordinates}
\label{Working on spaces which may not have nice coordinates}
\setcounter{equation}{0}

\renewcommand{\thetheorem}{D.\arabic{equation}}
\renewcommand{\theequation}{D.\arabic{equation}}

	If one has a topological or metric space (or whatever) which
has nice coordinates, then that can be pretty good.

	However, there is a lot that one can do without having
coordinates.  As in Section \ref{parameterization problems} and
Appendix \ref{more on existence and behavior of homeomorphisms}, there
are many situations in which homeomorphic coordinates might not be
available, or might be available only in irregular forms, or forms
with large complexity.  Even if piecewise-linear coordinates exist,
for instance, it may not be very nice if they have enormous
complexity, as in Section \ref{parameterization problems}.

	Or, as in Section \ref{parameterization problems}, local
coordinates might exist, but one might not have an algorithmic way to
know this.  Similarly, coordinates might exist, but it may not be so
easy to find them.

	As a brief digression, let us mention some positive results
about situations in which coordinates exist, but are not as regular as
one might like, at least not at first.  Consider the case of
topological manifolds, which admit homeomorphic local coordinates, but
for which there may not be a compatible piecewise-linear or smooth
structure.  Homeomorphic coordinates are not suitable for many basic
forms of analysis in which one might be interested, e.g., involving
differential operators.  However, there is a famous theorem of
Sullivan \cite{Su1}, to the effect that topological manifolds of
dimension $\ge 5$ admit unique quasiconformal and Lipschitz structures
(for which one then has quasiconformal and bilipschitz coordinates).
(In dimensions less than or equal to $3$, unique piecewise-linear and
smooth structures always exist for topological manifolds, by more
classical results.  In dimension $4$, quasiconformal and Lipschitz
structures may not exist for topological manifolds, or be unique.  See
\cite{DoS}.  Concerning smooth structures in dimension $4$, see
\cite{DoK, FQ}.  Some brief surveys pertaining to different structures
on manifolds, and in general dimensions, are given in Section 8 in
\cite{FQ}, and the ``Epilogue'' in \cite{Milnor-Stasheff}.)

	Thus, with Sullivan's theorem, one has the possibility of
improving the structure in a way that does make tools of analysis
feasible.  Some references related to this include \cite{CST1, CST2,
DoS, RsW, Su1, Sullivan-Teleman, Teleman1, Teleman2}.

	Some aspects of working on spaces without good coordinates
came up in Section \ref{Quantitative topology, and calculus on
singular spaces}.  One could also consider ``higher-order'' versions
of this, along the lines of differential forms.  We shall not pursue
this here, but for a clear and simple version of this, one can look at
the case of Euclidean spaces with the geometry deformed through a
metric doubling measure (as in \cite{DS-sob, Se-finn, Sebarc}).  Some
points about this are explained in \cite{Sebarc}, beginning near the
bottom of p427.  (Compare also with \cite{Sebil}, concerning the
possible behavior of Euclidean spaces with geometry deformed by metric
doubling measures.)

	In this appendix, we shall focus more on traditional objects
from algebraic topology, like homology and cohomology groups.  In this
setting, there is a lot of structure around, concerning spaces which
might be approximately like manifolds, but not quite manifolds.

	Let us begin with some basic conditions.  Let $M$ be a
topological space which is compact, Hausdorff, and metrizable.  We
shall assume that $M$ has finite topological dimension, in the sense
of \cite{HW}.  In these circumstances, this is equivalent to saying
that $M$ is homeomorphic to a subset of some ${\bf R}^n$.  (See
\cite{HW}.)

	For simplicity, let us imagine that $M$ simply \emph{is} a
compact subset of some ${\bf R}^n$.  It is also convenient to ask that
$M$ be \emph{locally contractable}.  This means (as in Subsubsection
\ref{A couple of comparisons between geometric and topological
properties}) that for each point $p \in M$, and each neighborhood $U$
of $p$ in $M$, there is a smaller neighborhood $V \subseteq U$ of $p$
in $M$ such that $V$ can be contracted to a point in $M$.

	As a class of examples, finite polyhedra are locally
contractable.  Finite polyhedra make a nice special case to consider
throughout this appendix, and we shall return to it several times.

	For another class of examples, one has cell-like quotients of
topological manifolds (and some locally contractable spaces more
generally), at least when the quotient spaces have finite topological
dimension.  See Corollary 12B on p129 of \cite{Dm2}.  (Cell-like
quotients were discussed somewhat in Appendix \ref{more on existence
and behavior of homeomorphisms}, Subsections \ref{Decomposition
spaces, 1} and \ref{Decomposition spaces, 2} in particular.  Some
concrete instances of cell-like quotients are mentioned in these
subsections, and \cite{Dm2} provides more examples and information.)

	Although we shall mostly not emphasize metric structures or
quantitative aspects in the appendix, let us mention that cell-like
quotients like these often have natural and nice geometries, as
indicated in Subsection \ref{Geometric structures for decomposition
spaces}.  These geometries are quite different from those of finite
polyhedra, but they can also have some analogous properties.  In
particular, there can be forms of self-similarity or scale-invariant
boundedness of the geometry, and these can be analogous to local
conical structure in polyhedra in their effects.  They are not as
strong or special, and they are also more flexible.  In any case,
local contractability (and conditions like local linear
contractability, in some geometric settings) is a basic property to
perhaps have.

	As a general fact about local contractability, let us note the
following.

\beginproposition
\label{compact subsets of R^n are locally contractable if and only if nhbd ret}
Let $M$ be a compact subset of some ${\bf R}^n$.  Then $M$ is locally
contractable if and only if there is a set $V \subseteq {\bf R}^n$ which
contains $M$ in its interior, and a continuous mapping $r : V \to M$ which
is a \emph{retract}, i.e., $r(w) = w$ for all $w \in V$.
\end{proposition}

	This is a fairly standard observation.  The ``if'' part is an easy
consequence of the local contractability of ${\bf R}^n$ (through linear
mappings).  I.e., to get local contractions inside of $M$, one makes
standard linear contractions in ${\bf R}^n$, which normally do not stay
inside $M$, and then one applies the retraction to keep the contractions
inside ${\bf R}^n$.

	For the converse, one can begin by defining $r$ on a discrete
and reasonably-thick set of points outside $M$, but near $M$.  For
such a point $w$, one could choose $r(w) \in M$ so that it lies as
close to $w$ as possible (among points in $M$), or is at least
approximately like this.  To fill in $r$ in the areas around these
discrete points, one can make extensions first to edges, then
$2$-dimensional faces, and so on, up to dimension $n$.  To make these
extensions, one uses local contractability of $M$.  It is also
important that the local extensions do not go to far from the
selections already made, so that $r : V \to M$ will be continuous in
the end, and this one can also get from the local contractability.

	The notion of ``Whitney decompositions'', as in Chapter VI of
\cite{St1}, is helpful for this kind of argument.  It gives a way of
decomposing ${\bf R}^n \backslash M$ into cubes with disjoint
interiors, and some other useful properties.  (In particular, this
kind of decomposition can be helpful for keeping track of bounds, if
one should wish to do so.)  One can use the vertices of these cubes
for the discrete set in the complement of $M$ mentioned above.

	See also \cite{Dm2} for a proof of Proposition \ref{compact
subsets of R^n are locally contractable if and only if nhbd ret},
especially p117ff.

	Let us return now to the general story.  Suppose that $M$ is a
compact subset of ${\bf R}^n$, and that $M$ is locally contractable.
Let $r : V \to M$ be a continuous retraction on $M$, as in Proposition
\ref{compact subsets of R^n are locally contractable if and only if
nhbd ret}.  Thus $V$ contains $M$ in its interior.  By replacing $V$
by a slightly smaller subset, if necessary, we may assume that $V$ is
compact, and in fact that it is a finite union of dyadic cubes in
${\bf R}^n$.  (A \emph{dyadic cube} in ${\bf R}^n$ is one which can be
represented as a Cartesian product of intervals $[j_i \, 2^{-k},
(j_i+1) \, 2^{-k}]$, $i = 1, 2, \ldots, n$, where the $j_i$'s and $k$
are integers.)

	This type of choice for $V$ is convenient for having nice
properties in terms of homology and cohomology.  In particular, $V$ is
then a finite complex.  The inclusion of $M$ into $V$, and the mapping
$r : V \to M$, induce mappings between the homology and cohomology of
$M$ and $V$.  If $\iota : M \to V$ denotes the mapping coming from
inclusion, then $r \circ \iota : M \to M$ is the identity mapping, and
thus it induces the identity mapping on the homology and cohomology of
$M$.  Using this, one can see that the mapping from the homology of
$M$ into the homology of $V$ induced by $\iota$ is an injection (in
addition to being a group homomorphism, as usual), and that the
mapping from the homology of $V$ to the homology of $M$ induced by $r$
is a surjection.  This follows from standard properties of homology
and mappings, as in \cite{Massey}.  Similarly, $r$ induces a mapping
from cohomology of $M$ to cohomology of $V$ which is injective, and
$\iota$ induces a mapping from cohomology of $V$ into cohomology of
$M$ which is surjective.

	This provides a simple way in which the algebraic topology of
$M$ can be ``bounded'', under the type of assumptions on $M$ that we
are making.  (There are more refined things that one can also do, but
we shall not worry about this here.)  Local contractability, and the
existence of a retraction as in Proposition \ref{compact subsets of
R^n are locally contractable if and only if nhbd ret}, are also nice
for making it clear and easy to work with continuous mappings into
$M$.  In particular, one can get continuous mappings into $M$ from
continuous mappings into $V$, when one has a retraction $r : V \to M$,
as above.  This is as opposed to standard examples like the closure of
the graph of $\sin (1 / x)$, $x \in [-1, 1] \backslash \{0\}$.  (This
set is connected but not arcwise connected.)

	Now let us consider the following stronger conditions on $M$.

\begindefinition [Generalized $k$-Manifolds]
Let $M$ be a compact subset of ${\bf R}^n$ which is locally contractable.
Then $M$ is a \emph{generalized $k$-manifold} if for every point $z \in M$,
the relative homology $H_j(M, M \backslash \{z\})$ is the same (up to
isomorphism) as the relative homology of $H_j({\bf R}^k, {\bf R}^k \backslash
\{0\})$ for each $j$.  (In other words, $H_j(M, M \backslash \{z\})$ should
be zero when $j \ne k$, and isomorphic to ${\bf Z}$ when $j=k$.)
\end{definition}

	We are implicitly working with homology defined over the
integers here, and there are analogous notions with respect to other
coefficient groups (like rational numbers, for instance).  One may
also wish to use weaker conditions than local contractability (as in
\cite{Bredon2, Wilder}).  There are other natural variations for this
concept.

	If $M$ is a finite polyhedron, then the property of being a
generalized manifold is equivalent to asking that the links of $M$ be
homology spheres (i.e., have the same homology as standard spheres, up
to isomorphism) of the right dimension.  Such a polyhedron may not be
a topological manifold (in dimensions greater than or equal to $4$),
because the codimension-$1$ links may not be simply-connected.  (This
is closely related to some of the topics of Appendix \ref{more on
existence and behavior of homeomorphisms}, and the condition
(\ref{local vanishing of pi_1 in widehat{M} backslash {q}}) in
Subsection \ref{Interlude: looking at infinity, or looking near a
point} in particular.)

	Another class of examples comes from taking quotients of
compact topological manifolds by cell-like decompositions (Subsections
\ref{Decomposition spaces, 1} and \ref{Decomposition spaces, 2}), at
least when the quotient space has finite topological dimension.  See
Corollary 1A on p191 of \cite{Dm2} (and Corollary 12B on p129 there),
and compare also with Theorem 16.33 on p389 of \cite{Bredon2}, and
\cite{Fe4}.  As in Appendix \ref{more on existence and behavior of
homeomorphisms} (especially Subsections \ref{Decomposition spaces, 1}
and \ref{Decomposition spaces, 2}) and the references therein, such
spaces are \emph{not} always topological manifolds themselves, in
dimensions $3$ and higher.  For a concrete example of this, one can
take the standard $3$-sphere, and collapse a copy of the Fox--Artin
(wild) arc to a point.  (Compare with \cite{Fe4}.)  One could also
collapse a Whitehead continuum to a point (where Whitehead continua
are as in Subsection \ref{Contractable open sets}), as in Subsection
\ref{Decomposition spaces, 1}.  As in Appendix \ref{more on existence
and behavior of homeomorphisms}, the fact that these spaces are not
manifolds can be seen in the nontriviality of localized fundamental
groups, i.e., fundamental groups of the complement of the
distinguished point (corresponding to the Fox--Artin curve or
Whitehead continuum), localized at that point.  This is analogous to
the possibility of having codimension-$1$ links in polyhedra which are
not simply-connected, as in the previous paragraph.

	As usual, dimensions $1$ and $2$ are special for generalized
manifolds, which are then topological manifolds.  See \cite{Wilder},
Theorem 16.32 on p388 of \cite{Bredon2}, and the introduction to
\cite{Fe4}.

	For more on ways that generalized manifolds can arise, see
\cite{Borel2, Bredon2, BFMW1, BFMW2, Dm2, Fe4, Weinberger-book} (and
the references therein).  A related topic is the ``recognition
problem'', for determining when a topological space is a topological
manifold.  This is also closely connected to the questions in Appendix
\ref{more on existence and behavior of homeomorphisms}.  Some
references for this include \cite{BFMW1, BFMW2, C1, C-bulletin, C2,
Dm2, E, Fe4, Weinberger-book}.

	What are some properties of generalized manifolds?  In what
ways might they be like manifolds?  In particular, how might they be
different from compact sets in ${\bf R}^n$ which are locally
contractable in general?

	A fundamental point is that Poincar\'e duality (and other
duality theorems for manifolds) also work for generalized manifolds.
See \cite{Borel1, Borel2, Bredon2, Wilder} and p277-278 of \cite{Sp}.
This is pretty good, since Poincar\'e duality is such a fundamental
aspect of manifolds.  (See \cite{BT, Bredon, Massey, Milnor-Stasheff,
Sp}, for instance.)

	A more involved fact is that rational Pontrjagin classes can
be defined for generalized manifolds.  (See the introduction to
\cite{BFMW1}.)  

	For \emph{smooth} manifolds, the definition of the Pontrjagin
classes is classical.  (See \cite{BT, Milnor-Stasheff}.)  More
precisely, one can define Pontrjagin classes for vector bundles in
general, and then apply this to the tangent bundle of a smooth
manifold to get the Pontrjagin classes of a manifold.  As
\emph{integral} cohomology classes, the Pontrjagin classes are
preserved by diffeomorphisms between smooth manifolds, but not, in
general, by homeomorphisms.  However, a famous theorem of Novikov is
that the Pontrjagin classes of smooth manifolds are preserved as
\emph{rational} cohomology classes by homeomorphisms in general.
Further developments lead to the definition of rational characteristic
classes on more general spaces.

	For finite polyhedra, there is an earlier treatment of
rational Pontrjagin classes, which goes back to work of Thom and
Rohlin and Schwarz.  See Section 20 of \cite{Milnor-Stasheff}.  More
precisely, this gives a procedure by which to define rational
Pontrjagin classes for finite polyhedra which are generalized
manifolds, and which is invariant under piecewise-linear equivalence.
(For this, the generalized-manifold condition can be given in terms of
rational coefficients for the homology groups.)  If one starts with a
smooth manifold, then there it can be converted to a piecewise-linear
manifold (unique up to equivalence) by earlier results, and the
classical rational Pontrjagin classes for the smooth manifold are the
same as the ones that are obtained by the procedure for polyhedral
spaces.  (See \cite{Milnor-Stasheff} for more information.)

	The results mentioned above indicate some of the ways that
generalized manifolds are like ordinary manifolds.  In particular,
there is a large extent to which one can work with them, including
making computations or measurements on them, in ways that are similar
to those for ordinary manifolds. 

	This is pretty good!  This is especially nice given the
troubles that can come with homeomorphisms, as in Section
\ref{parameterization problems} and Appendix \ref{more on existence
and behavior of homeomorphisms}.  I.e., homeomorphisms can be
difficult to get or have, to begin with; even if they exist, they may
necessarily be irregular, as in the case of double-suspensions of
homology spheres \cite{SS}, manifold factors (Subsection \ref{Manifold
factors}), some other classical decomposition spaces (Subsections
\ref{Decomposition spaces, 2}, and \ref{Geometric structures for
decomposition spaces}), and homeomorphisms between $4$-dimensional
manifolds (see \cite{DoK, DoS, Fr, FQ}); even if homeomorphisms exist
and are of a good regularity class, their complexity may have to be
very large, as in Section \ref{parameterization problems} and the
results in \cite{BHP}.

	In the case of finite polyhedra, the condition of being a
generalized manifold involves looking at the homology groups of the
links (as mentioned before), and this is something which behaves in a
fairly nice and stable way.  Compare with Appendix \ref{A few facts
about homology}.  By contrast, the property of being a manifold
involves the fundamental groups of the links (at least in codimension
$1$ for having topological manifolds), and this is much more
complicated.  We have run into this already, in Section
\ref{parameterization problems}.

	More generally, we have also seen in Appendix \ref{more on
existence and behavior of homeomorphisms} how conditions related to
localized fundamental groups can arise, in connection with the
existence of homeomorphisms and local coordinates.  This includes the
vanishing of $\pi_1$ in some punctured neighborhoods of points in
topological manifolds of dimension at least $3$.  With generalized
manifolds, one has certain types of special structure, but one does
not necessarily have localized $\pi_1$ conditions like these.  This
gives a lot of simplicity and stableness, as well as flexibility.

\section{Some simple facts related to homology}
\label{A few facts about homology}
\setcounter{equation}{0}

\renewcommand{\thetheorem}{E.\arabic{equation}}
\renewcommand{\theequation}{E.\arabic{equation}}

	Let $P$ be a finite polyhedron, which we shall think of as
being presented to us as a finite simplicial complex.  Fix a positive
integer $k$, and imagine that one is interested in knowing whether the
homology $H_k(P)$ of $P$ (with coefficients in ${\bf Z}$) is $0$ in
dimension $k$.  This is a question that can be decided by an
algorithm, unlike the vanishing of the fundamental group.  This fact
is standard and elementary, but not necessarily too familiar in all
quarters.  It came up in Section \ref{parameterization problems}, and,
for the sake of completeness, a proof will be described here.  (Part
of the point is to make it clear that there are no hidden surprises
that are too complicated.  We shall also try to keep the discussion
elementary and direct, with a minimum of machinery involved.)

	The first main point is that it suffices to consider
\emph{simplicial} homology of $P$, rather than something more general
and elaborate (like \emph{singular} homology).  This puts strong
limits on the type of objects with which one works.  By contrast, note
that the higher homotopy groups of a finite complex need not be
finitely-generated.  See Example 17 on p509 of \cite{Sp}.  This is true,
but quite nontrivial, for simply-connected spaces.  See Corollary 16
on p509 of \cite{Sp}.  

	Let us begin by considering the case of homology with
coefficients in the rational numbers, rather than the integers.  In
this situation the homology is a vector space over ${\bf Q}$, and this
permits the solution of the problem through means of linear algebra.
To be explicit, one can think of the vanishing of the homology in
dimension $k$ in the following terms.  One starts with the set of
$k$-dimensional \emph{chains} in $P$ (with coefficients in ${\bf Q}$),
which is the set of formal sums of oriented $k$-dimensional simplices
with coefficients in ${\bf Q}$.  This is a vector space, and
multiplication by $-1$ is identified with reversing orientations on
the simplices.  The set of $k$-dimensional \emph{cycles} then consists
of $k$-dimensional chains with ``boundary'' equal to $0$.  This can be
described by a finite set of linear equations in our vector space of
$k$-dimensional chains.  The set of $k$-dimensional \emph{boundaries}
consist of $k$-dimensional chains which are themselves boundaries of
$(k+1)$-dimensional chains.  This is the same as taking the span of
the boundaries of $(k+1)$-dimensional simplices in $P$, a finite set.
Chains which are boundaries are automatically cycles, and the
vanishing of the $k$-dimensional homology of $P$ (with coefficients in
${\bf Q}$) is equivalent to the equality of the vector space of cycles
with the vector space of boundaries.  The problem of deciding whether
this equality holds can be reduced to computations of determinants,
for instance, and one can use other elementary techniques from linear
algebra too.

	In working with integer coefficients, one has the same basic
definitions of chains, cycles, and boundaries, but now with
coefficients in ${\bf Z}$ rather than ${\bf Q}$.  The spaces of
chains, cycles, and boundaries are now abelian groups, and the
equations describing cycles and boundaries make sense in this context.
``Abelian'' is a key word here, and a crucial difference between this
and the fundamental group.  Homotopy groups $\pi_j$ in dimensions $j
\ge 2$ are always abelian too, but they do not come with such a simple
presentation.

	Let us be more explicit again.  The space of $k$-dimensional
chains in $P$, now with integer coefficients, is a free abelian group,
generated by the $k$-dimensional simplices in $P$.  One can think of
it as being realized concretely by ${\bf Z}^r$, where $r$ is the
number of $k$-dimensional simplices in $P$.

	The set of \emph{cycles} among the chains is defined by a
finite number of linear equations, as above.  That is, the cycles $z
\in {\bf Z}^r$ are determined by finitely many equations of the form
\begin{equation}
\label{equation}
	a_1 z_1 + a_2 z_2 + \cdots + a_r z_r = 0, 
			\quad z = (z_1, z_2, \ldots, z_r) \in {\bf Z}^r,
\end{equation}
where the $a_i$'s are themselves integers.  In fact, for the equations
which actually arise in this situation, the $a_i$'s are always either
$0$, $1$, or $-1$.  It will be useful, though, to allow for general
vectors $(a_1, a_2, \ldots, a_r)$ in this discussion, and our
computations will lead us to that anyway.

	We begin with an observation about sets of vectors in ${\bf
Z}^r$ defined by a single homogeneous linear equation (with integer
coefficients).

\beginlemma
\label{lemma about one equation on {bf Z}^r}
Let $r$ be a positive integer, and let $a = (a_1, a_2, \ldots, a_r)$
be a vector of integers.  Set
\begin{equation}
\label{def of C_a}
	C_a = \{z \in {\bf Z}^r : a_1 z_1 + a_2 z_2 + \cdots + a_r z_r = 0\}.
\end{equation}
Then there is a group homomorphism $\phi_a : {\bf Z}^r \to {\bf Z}^r$
such that $\phi_a({\bf Z}^r) = C_a$ and $\phi_a(z) = z$ when $z \in C_a$.
This homomorphism can be effectively constructed given $r$ and $a$.
\end{lemma}

	Let us prove Lemma \ref{lemma about one equation on {bf Z}^r}.
Let $r$ and $a$ be given, as above.  We may as well assume that the
$a$ is not the zero vector, since if it were, $C_a$ would be all of
${\bf Z}^r$, and we could take $\phi_a$ to be the identity mapping.

	We may also assume that the components $a_i$ of $a$ are not all
divisible by an integer different from $1$ or $-1$, since we can always
cancel out common factors from the $a_i$'s without changing $C_a$.

\beginsublemma
\label{b from a}
Under these conditions ($a$ is not the zero vector, and no integers
divide all of the components of $a$ except $\pm 1$), there exists a
vector $b = (b_1, b_2, \ldots, b_r)$ of integers such that
\begin{equation}
\label{a-b equation}
	a_1 b_1 + a_2 b_2 + \cdots + a_r b_r = 1.
\end{equation}
A choice of $b$ can be constructed effectively given $a$.
\end{sublemma}

	The $r=2$ version of this is solved by the well known
``Euclidean algorithm''.  It is more commonly formulated as saying
that if one has two nonzero integers $c$ and $d$, then one can find
integers $e$ and $f$ such that $c e + d f$ is equal to the (positive)
greatest common divisor of $c$ and $d$.  This can be proved using
the more elementary ``division algorithm'', to the effect that if
$m$ and $n$ are positive integers, with $m < n$, then one can write
$n$ as $jm + i$, where $j$ is a positive integer and $0 \le i < m$.

	For general $r$ one can reduce to the $r=2$ case using
induction, for instance.  This is not hard to do, and we omit the
details.  (We should perhaps say also that the $r=1$ case makes sense,
and is immediate, i.e., $a_1$ has to be $\pm 1$.)

	Let us return now to Lemma \ref{lemma about one equation on
{bf Z}^r}.  Given any $z \in {\bf Z}^r$, write $z \cdot a$ for $z_1
c_1 + z_2 c_2 + \cdots + z_r c_r$.  Define $\phi_a : {\bf Z}^r \to
{\bf Z}^r$ by putting
\begin{equation}
\label{def of phi_a}
	\phi_a(z) = z - (z \cdot a) \, b,
\end{equation}
where $b$ is chosen as in Sublemma \ref{b from a}.  This is clearly a
group homomorphism (with respect to addition).  If $z$ lies in $C_a$,
then $z \cdot a = 0$, and $\phi_a(z) = z$.  In particular,
$\phi_a({\bf Z}^r) \supseteq C_a$.  On the other hand, if $z$ is any
element of ${\bf Z}^r$, then
\begin{equation}
	\phi_a(z) \cdot a = z \cdot a - (z \cdot a) (b \cdot a),
\end{equation}
and this is equal to $0$, by (\ref{a-b equation}).  Thus $\phi_a(z)
\in C_a$ for all $z \in {\bf Z}^r$, so that $\phi_a({\bf Z}^r)
\subseteq C_a$ too.  This proves Lemma \ref{lemma about one equation
on {bf Z}^r}.

	Next we consider the situation of sets in ${\bf Z}^r$ defined
by multiple linear equations (with linear coefficients).

\beginlemma
\label{lemma about several equations on {bf Z}^r}
Let $a^1, a^2, \ldots, a^p$ be a collection of vectors in
${\bf Z}^r$, and define $C \subseteq {\bf Z}^r$ by
\begin{equation}
\label{def of C}
   C = \{z \in {\bf Z}^r : a^i \cdot z = 0 \ \hbox{for } i = 1, 2, \dots, p \}.
\end{equation}
Then there is a group homomorphism $\psi : {\bf Z}^r \to {\bf Z}^r$
such that $\psi({\bf Z}^r) = C$ and $\psi(z) = z$ when $z \in C$.
This homomorphism can be effectively constructed given $r$ and the
vectors $a^1, a^2, \ldots, a^p$.
\end{lemma}

	To prove this, we use induction, on $p$.  The point is
basically to iterate the construction of Lemma \ref{lemma about one
equation on {bf Z}^r}, but one has to be a bit careful not to disrupt
the previous work with the new additions.

	When $p = 1$, Lemma \ref{lemma about several equations on {bf
Z}^r} is the same as Lemma \ref{lemma about one equation on {bf Z}^r}.
Now suppose that we know Lemma \ref{lemma about several equations on
{bf Z}^r} for some value of $p$, and that we want to verify it for
$p+1$.

	Let $a^1, a^2, \ldots, a^{p+1}$ be a collection of vectors in
${\bf Z}^r$, and let $C^{p+1}$ denote the set of vectors $z \in {\bf
Z}^r$ which satisfy $a^i \cdot z = 0$ for $i = 1, 2, \ldots, p+1$, as
in (\ref{def of C}).  Similarly, let $C^p$ denote the set of $z \in
{\bf Z}^r$ such that $a^i \cdot z = 0$ when $1 \le i \le p$.  Our
induction hypothesis implies that there is a group homomorphism
$\psi^p : {\bf Z}^r \to {\bf Z}^r$ such that $\psi^p({\bf Z}^r) = C^p$
and $\psi^p(z) = z$ when $z \in C^p$, and which can be effectively
constructed given $r$ and $a^1, a^2, \ldots, a^p$.  We want to produce
a similar homomorphism $\psi^{p+1}$ for $C^{p+1}$.

	The basic idea is to apply Lemma \ref{lemma about one equation
on {bf Z}^r} to the vector $a^{p+1}$.  This does not quite work, and
so we first modify $a^{p+1}$ to get a vector which has practically the
same effect as $a^{p+1}$ for defining $C^{p+1}$, and which is in a
more convenient form.

	Specifically, let $\alpha$ be the vector in ${\bf Z}^r$ such
that
\begin{equation}
\label{choice of alpha}
	\alpha \cdot z = a^{p+1} \cdot \psi^p(z) 	
				\qquad\hbox{for all } z \in {\bf Z}^r.
\end{equation}
It is easy to compute $\alpha$ given $a^{p+1}$ and $\psi^p$.  One can
also describe $\alpha$ as the vector which results by applying the
\emph{transpose} of $\psi^p$ to $a^{p+1}$.

	For our purposes, the main properties of $\alpha$ are as follows.
First,
\begin{equation}
\label{first property of alpha}
	\alpha \cdot z = a^{p+1} \cdot z \qquad \hbox{when } z \in C^p.
\end{equation}
This is a consequence of (\ref{choice of alpha}) and the fact that
$\psi^p(z) = z$ when $z \in C^p$.  Second,
\begin{equation}
\label{second property of alpha}
	\alpha \cdot \psi^p(z) = \alpha \cdot z 
				\qquad \hbox{for all } z \in {\bf Z}^r.
\end{equation}
To see this, note that $\psi^p(\psi^p(z)) = \psi^p(z)$ for all $z$,
since $\psi^p$ maps ${\bf Z}^r$ into $C^p$ and $\psi^p$ is equal to
the identity on $C^p$, by construction.  This and (\ref{choice of
alpha}) give (\ref{second property of alpha}).

	Because of (\ref{first property of alpha}), we have that
\begin{equation}
\label{C^{p+1} = ..., 1}
	C^{p+1} = \{z \in C^p : \alpha \cdot z = 0\}.
\end{equation}
That is, $C^{p+1} = \{z \in C^p : a^{p+1} \cdot z = 0\}$, by the definition
of $C^{p+1}$ and $C^p$, and then (\ref{C^{p+1} = ..., 1}) follows from this
and (\ref{first property of alpha}).

	If $\alpha = 0$, then $C^{p+1} = C^p$, and we can stop here.
That is, we can use $\psi^p$ for $\psi^{p+1}$.  Thus we assume instead
that $\alpha$ is not the zero vector.

	Let $\alpha'$ denote the vector in ${\bf Z}^r$ obtained from
$\alpha$ by eliminating all common factors of the components of
$\alpha$ which are positive integers greater than $1$.  In particular,
$\alpha$ is a positive integer multiple of $\alpha'$.  From
(\ref{second property of alpha}) and (\ref{C^{p+1} = ..., 1}), we get
that
\begin{equation}
\label{identity for alpha'}
	\alpha' \cdot \psi^p(z) = \alpha' \cdot z
				\qquad \hbox{for all } z \in {\bf Z}^r
\end{equation}
and
\begin{equation}
\label{C^{p+1} = ..., 2}
	C^{p+1} = \{z \in C^p : \alpha' \cdot z = 0\}.
\end{equation} 
	
	Now let us apply Sublemma \ref{b from a} to obtain a vector
$\beta' \in {\bf Z}^r$ such that
\begin{equation}
	\alpha' \cdot \beta' = 1.
\end{equation}
Put $\beta'' = \psi^p(\beta')$.  Then
\begin{equation}
\label{alpha' cdot beta'' = 1}
	\alpha' \cdot \beta'' = 1,
\end{equation}
because of (\ref{identity for alpha'}).  Also, $\beta''$ lies in
$C^p$, since $\psi^p$ maps ${\bf Z}^r$ into $C^p$.  (It was for this
purpose that we have made the above modifications to $a^{p+1}$, i.e.,
to get $\beta'' \in C^p$.)

	Now that we have all of this, we can proceed exactly as in the
proof of Lemma \ref{lemma about one equation on {bf Z}^r}.  Specifically,
define $\phi' : {\bf Z}^r \to {\bf Z}^r$ by 
\begin{equation}
\label{def of phi'}
	\phi'(z) = z - (z \cdot \alpha') \, \beta''.
\end{equation}
A key point is that
\begin{equation}
\label{phi'(C^p) subseteq C^p}
	\phi'(C^p) \subseteq C^p,
\end{equation}
which holds since $\beta'' \in C^p$, as mentioned above.  Also,
\begin{equation}
	\phi'(z) \cdot \alpha' = 0 \qquad\hbox{for all } z \in {\bf z}^r,
\end{equation}
because of (\ref{alpha' cdot beta'' = 1}) and the definition (\ref{def
of phi'}) of $\phi'(z)$.  From this and (\ref{phi'(C^p) subseteq C^p})
it follows that
\begin{equation}
\label{phi'(C^p) subseteq C^{p+1}}
	\phi'(C^p) \subseteq C^{p+1},
\end{equation}
using also (\ref{C^{p+1} = ..., 2}).

	On the other hand, 
\begin{equation}
\label{phi'(z) = z whenever ...}
	\phi'(z) = z \qquad\hbox{whenever } z \in {\bf Z}^r, 
						\ z \cdot \alpha' = 0
\end{equation}
(by the definition (\ref{def of phi'}) of $\phi'$).  In particular,
$\phi'(z) = z$ when $z \in C^{p+1}$ (using (\ref{C^{p+1} = ..., 2})
again).  Combining this with (\ref{phi'(C^p) subseteq C^{p+1}}), we
get that
\begin{equation}
\label{phi'(C^p) = C^{p+1}}
	\phi'(C^p) = C^{p+1}.
\end{equation}

	Now we are almost finished with the proof.  Define $\psi^{p+1}
: {\bf Z}^r \to {\bf Z}^r$ by $\psi^{p+1} = \phi' \circ \psi^p$.  This
defines a group homomorphism (with respect to the usual addition of
vectors), since $\psi^p$ and $\phi'$ are group homomorphisms.  We
also have that
\begin{equation}
	\psi^{p+1}({\bf Z}^r) = C^{p+1};
\end{equation}
this follows from (\ref{phi'(C^p) = C^{p+1}}) and the fact that
$\psi^p({\bf Z}^r) = C^p$, which was part of our ``induction
hypothesis'' on $\psi^p$.

	Let us verify that
\begin{equation}
\label{psi^{p+1}(z) = z whenever ...}
	\psi^{p+1}(z) = z \qquad\hbox{whenever } z \in C^{p+1}.
\end{equation}
If $z \in C^{p+1}$, then $z \in C^p$ in particular.  This implies that
$\psi^p(z) = z$, again by our ``induction hypothesis'' for $\psi^p$.
Thus $\psi^{p+1}(z) = \phi'(z)$ in this case.  We also have that
$\phi'(z) = z$ if $z \in C^{p+1}$, because of (\ref{phi'(z) = z
whenever ...})  (and (\ref{C^{p+1} = ..., 2})).  This gives
(\ref{psi^{p+1}(z) = z whenever ...}), as desired.

	It is easy to see from this construction that $\psi^{p+1}$ is
effectively computable from the knowledge of $r$ and $a^1, a^2,
\ldots, a^{p+1}$, given the corresponding assertion for $\psi^p$.
Thus $\psi^{p+1}$ has all the required properties.  This completes the
proof of Lemma \ref{lemma about several equations on {bf Z}^r}.

	Now let us return to the original question, about determining
whether the integral homology of a given finite complex vanishes in a
given dimension.  We can put this into a purely algebraic form, as
follows.  Suppose that a positive integer $r$ is given, as well as two
finite collections of vectors in ${\bf Z}^r$, $a^1, a^2, \ldots, a^p$,
and $d^1, d^2, \ldots, d^q$.  Given this data, the problem asks,
\begin{eqnarray}
\label{integer question, first version} 
	&& \hbox{is it true that for every $z \in {\bf Z}^r$ which satisfies }
		a^i \cdot z = 0 					\\
        && \hbox{for all } i = 1, 2, \ldots, p,  \hbox{ there exist } 
            \lambda_1, \lambda_2, \ldots, \lambda_q \in {\bf Z}	\nonumber \\
        && \hbox{such that }
            z = \lambda_1 d^1 + \lambda_2 \, d^2 + \cdots + \lambda_q \, d^q?
                                                         \nonumber
\end{eqnarray}
The question of vanishing of homology is of this form, and so an
effective procedure for determining an answer of ``yes'' or ``no'' to
(\ref{integer question, first version}) also provides a way to decide
whether the homology of a given finite complex vanishes in a given
dimension.

	Let us use Lemma \ref{lemma about several equations on {bf
Z}^r} to reduce (\ref{integer question, first version}) to a simpler
problem, as follows.  Let $r \in {\bf Z}_+$ and $d^1, d^2, \ldots, d^q
\in {\bf Z}^r$ be given, and also another vector $z \in {\bf Z}^r$.
Given this data, the new problem asks,
\begin{eqnarray}
\label{integer question, second version}
	&& \hbox{is it true that there exist }
             \lambda_1, \lambda_2, \ldots, \lambda_q \in {\bf Z}    \\
        && \hbox{such that }
     z = \lambda_1 d^1 + \lambda_2 \, d^2 + \cdots + \lambda_q \, d^q?
                                                           \nonumber
\end{eqnarray}
To show that (\ref{integer question, first version}) can be reduced to
(\ref{integer question, second version}), let us first mention an
auxiliary observation.

	Suppose that vectors $a^1, a^2, \ldots, a^p$ in ${\bf Z}^r$
are given.  Consider the set
\begin{equation}
\label{C = { w in {bf Z}^r : ....}}
        C = \{ w \in {\bf Z}^r : a^i \cdot w = 0 \hbox{ for }  
                                                i = 1, 2, \ldots, p \}.
\end{equation}
Then there is a finite collection of vectors $z^j \in C$ which
generate $C$ (as an abelian group), and which can be obtained through
an effective procedure (given $r$ and $a^1, a^2, \ldots, a^p$).  This
follows from Lemma \ref{lemma about several equations on {bf Z}^r}.
Specifically, for the $z^j$'s one can take $\psi^p(e^j)$, $1 \le j \le
r$, where $\psi^p : {\bf Z}^r \to {\bf Z}^r$ is the homomorphism
provided by Lemma \ref{lemma about several equations on {bf Z}^r}, and
$e^j$ is the $j$th standard basis vector in ${\bf Z}^r$, i.e., the
vector whose $j$th component is $1$ and whose other components are
$0$.  From Lemma \ref{lemma about several equations on {bf Z}^r}, we
know that $\psi^p$ maps ${\bf Z}^r$ onto $C$, and that $\psi^p$ can be
effectively produced given $r$ and $a^1, a^2, \ldots, a^p$.  This
implies that $z^j = \psi^p(e^j)$, $1 \le j \le r$, generate $C$ and
can be effectively produced as well.  (One might analyze this further
to reduce the number of vectors in the generating set, but this is not
needed for the present purposes.)

	Thus, in order to determine the answer to the question in
(\ref{integer question, first version}), one can first apply this
observation to produce a finite set of generators $z^j$ for $C$ as in
(\ref{C = { w in {bf Z}^r : ....}}).  An answer of ``yes'' for the
question in (\ref{integer question, first version}) is then equivalent
to having an answer of ``yes'' for the question in (\ref{integer
question, second version}) for each $z^j$ (in the role of $z$ in
(\ref{integer question, second version})).  In this way, we see that
an algorithm for deciding the answer to (\ref{integer question,
second version}) gives rise to an algorithm for determining the
answer to (\ref{integer question, first version}).

	Now let us consider (\ref{integer question, second version}).
Imagine that $r, q \in {\bf Z}^r$ and $d^\ell \in {\bf Z}^r$, $1 \le \ell
\le q$, are given.  We want to know if there is a vector $\lambda \in
{\bf Z}^q$, $\lambda = (\lambda_1, \lambda_2, \ldots, \lambda_q)$,
such that
\begin{equation}
\label{z = ..., 1}
	z = \lambda_1 d^1 + \lambda_2 \, d^2 + \cdots + \lambda_q \, d^q.
\end{equation}
Let us rewrite this as 
\begin{equation}
\label{z = ..., 2}
   z_i = \lambda_1 d^1_i + \lambda_2 \, d^2_i + \cdots + \lambda_q \, d^q_i
                                 \qquad\hbox{for } i = 1, 2, \ldots, r,
\end{equation}
where $z_i$, $d^\ell_i$, $1 \le i \le r$, denote the $i$th components
of $z$, $d^\ell$, respectively.

	Define vectors $\delta^i \in {\bf Z}^q$ by $\delta^i_j = d^j_i$
for $j = 1, 2, \ldots, q$.  With this use of ``transpose'' we can rewrite
(\ref{z = ..., 2}) as
\begin{equation}
\label{z = ..., 3}
	z_i = \lambda \cdot \delta^i \qquad\hbox{for } i = 1, 2, \ldots, r.
\end{equation}
Here ``$\cdot$'' denotes the usual dot product for vectors, although
now for vectors in ${\bf Z}^q$, rather than ${\bf Z}^r$, as before.

	With (\ref{z = ..., 3}), we are in a somewhat similar
situation as we have considered before, but with inhomogeneous
equations rather than homogeneous ones.  Note that for the
inhomogeneous equations there can be issues of torsion, i.e., it may be
that no solution $\lambda \in {\bf Z}^q$ to (\ref{z = ..., 3}) exists,
but a solution does exist if one replaces $z_i$ with $n z_i$ for some
positive integer $n$.

	We want to have an effective procedure for deciding when a
vector $\lambda \in {\bf Z}^q$ exists which provides a solution to
(\ref{z = ..., 3}).  To do this, we shall try to systematically reduce
the number of equations involved.  Consider first the equation
with $i=1$, i.e.,
\begin{equation}
\label{equation with i=1}
	z_1 = \lambda \cdot \delta^1.
\end{equation}
If there is no $\lambda \in {\bf Z}^q$ which satisfies this single
equation, then there is no solution for the system (\ref{z = ..., 3})
either.  If there are solutions to this equation, then we can try to
analyze the remaining equations on the set of $\lambda$'s which
satisfy this equation.

	In fact it is easy to say exactly when there is a $\lambda \in
{\bf Z}^q$ which satisfies (\ref{equation with i=1}).  A necessary
condition is that $z_1$ be divisible by all nonzero integers that
divide each component $\delta^1_j$, $1 \le j \le q$, of $\delta^1$.
(If $z_1$ or some $\delta^1_j$ is $0$, then it is divisible by all
integers.)  This necessary condition is also sufficient, because of
Sublemma \ref{b from a}.  The validity or not of this condition can be
determined effectively, and, when the condition holds, a particular
solution $\widetilde{\lambda}$ of (\ref{equation with i=1}) can be
produced effectively from the knowledge of $z_1$ and $\delta^1$,
because of Sublemma \ref{b from a}.

	If no solution to (\ref{equation with i=1}) exists in ${\bf
Z}^q$, then one can simply stop, as the answer to the question of the
existence of a solution to the system (\ref{z = ..., 3}) is then known
to be ``no''.  Let us suppose therefore that there is at least one
solution to (\ref{equation with i=1}).  As in the preceding paragraph,
this means that there is a solution $\widetilde{\lambda}$ which can be
effectively computed from the data.

	Set $L_1 = \{\lambda \in {\bf Z}^q : z_1 = \lambda \cdot
\delta^1\}$.  The existence of a solution $\lambda \in {\bf Z}^q$ to
the original system of equations in (\ref{z = ..., 3}) is equivalent
to the existence of a $\lambda \in L_1$ which satisfies 
\begin{equation}
\label{smaller system}
	z_i = \lambda \cdot \delta^i  \qquad\hbox{for } i = 2, \ldots, r.
\end{equation}
Thus we have reduced the number of equations involved, at the cost of
having a possibly more complicated set of $\lambda$'s as admissible
competitors.

	However, we can rewrite $L_1$ as
\begin{equation}
\label{L_1 revisited}
	L_1 = \{\lambda \in {\bf Z}^q : 
                  (\lambda - \widetilde{\lambda}) \cdot \delta^1 = 0\}.
\end{equation}
Set $L_1' = \{\tau \in {\bf Z}^r : \tau \cdot \delta^1 = 0\}$.  We
can reformulate the question of whether there is a $\lambda \in L_1$
such that (\ref{smaller system}) holds as asking whether there is a
$\tau \in L_1'$ such that
\begin{equation}
\label{smaller system, with tau now}
	z_i + \widetilde{\lambda} \cdot \delta^i = \tau \cdot \delta^i
				\qquad\hbox{for } i= 2, \ldots, r.
\end{equation} 

	In other words, we can make a change of variables and modify
the equations slightly so that the set $L_1'$ in which we look for
solutions is defined by a \emph{homogeneous} equation.  This permits
us to apply Lemma \ref{lemma about one equation on {bf Z}^r} (with $r$
replaced by $q$) to get a homomorphism $\phi : {\bf Z}^q \to {\bf
Z}^q$ such that $\phi({\bf Z}^q) = L_1'$.  Using this, our question
now becomes the following:
\begin{eqnarray}
\label{smaller system, no constraint}
	&& \hbox{does there exist } \xi \in {\bf Z}^q \hbox{ such that } \\
        && z_i + \widetilde{\lambda} \cdot \delta^i = \phi(\xi) \cdot \delta^i
		\hbox{ for } i = 2, \ldots, r?			\nonumber
\end{eqnarray}
This is equivalent to the earlier question, because the set of vectors
$\tau \in {\bf Z}^q$ which lie in $L_1'$ is the same as the set of
vectors of the form $\phi(\xi)$, where $\xi$ is allowed to be any element
of ${\bf Z}^q$.

	From Lemma \ref{lemma about one equation on {bf Z}^r} we know
that $\phi$ can be effectively computed from the knowledge of $q$ and
$\delta^1$.  By making straightforward substitutions, we can rewrite
the equations in (\ref{smaller system, no constraint}) as
\begin{equation}
	\widehat{z}_i = \xi \cdot \widehat{\delta}^i  
			\qquad\hbox{for } i = 2, \ldots, r,
\end{equation}
where the integers $\widehat{z_i}$ and vectors $\widehat{\delta}^i$
can be computed in terms of the original $z_i$'s and $\delta^i$'s,
$\widetilde{\lambda}$, and $\phi$.  In particular, they can be
computed effectively in terms of the original data in the problem.

	Thus we are back to the same kind of problem as we started with,
asking about the existence of a vector in ${\bf Z}^q$ which satisfies
a family of inhomogeneous linear equations.  However, now we have 
reduced the number of equations by $1$.  By repeating the process, we
can reduce to the case of a single inhomogeneous equation, which we
know how to solve (as we saw before).

	This shows that there is an effective method to determine the
answer to the question in (\ref{integer question, second version}).
As indicated earlier, this also gives a method for deciding the answer
to the question in (\ref{integer question, first version}), and to
the original problem about vanishing of homology in a finite complex.

	Let us briefly mention a cruder and more naive approach to
(\ref{integer question, second version}).  If the answer to
(\ref{integer question, second version}) is ``yes'', then it means
that there do exist integers $\lambda_1, \lambda_2, \ldots, \lambda_q$
which satisfy $z = \lambda_1 d^1 + \cdots + \lambda_q \, d^q$.  To
look for an answer of ``yes'' to (\ref{integer question, second
version}), one can simply start searching among all vectors $\lambda =
(\lambda_1, \lambda_2, \ldots, \lambda_q)$ in ${\bf Z}^q$, stopping
one finds a $\lambda$ which satisfies the equation above.

	If the answer to the question in (\ref{integer question,
second version}) is ``no'', then this search will not produce an
answer in a finite amount of time.  However, in this case one can make
a ``dual'' search to find a reason for the vector $z$ not to be in the
subgroup of ${\bf Z}^r$ generated by $d^1, d^2, \ldots, d^q$, a reason
which can also be found in finite time, when it exists.  Specifically,
$z$ does not lie in the subgroup of ${\bf Z}^r$ generated by $d^1,
d^2, \ldots, d^q$ if and only if there is a homomorphism $\sigma :
{\bf Z}^r \to {\bf Q}/{\bf Z}$ such that $\sigma(d^j) = 0$ for each
$j$ but $\sigma(z) \ne 0$.  We shall explain why this is true in a
moment, but first let us notice how this ``reason'' for an answer of
``no'' does fit our purpose.

	A homomorphism $\sigma : {\bf Z}^r \to {\bf Q}/{\bf Z}$ can be
described by $r$ elements of ${\bf Q}/{\bf Z}$, i.e., the values of
$\sigma$ on the $r$ standard basis vectors in ${\bf Z}^r$.  Any
elements of ${\bf Q}/{\bf Z}$ can be used here, and elements of ${\bf
Q}/{\bf Z}$ can be described in finite terms, i.e., by pairs of
integers.  The conditions $\sigma(d^j) = 0$, $1 \le j \le q$, and
$\sigma(z) \ne 0$, can be verified in finite time in a straightforward
manner.  Thus, if a $\sigma : {\bf Z}^r \to {\bf Q}/{\bf Z}$ exists
with these properties, then it can be found in a finite amount of
time, through an exhaustive search.

	If no such $\sigma$ exists, then this exhaustive search will
not stop in finite time.  However, the exhaustive search for the
$\lambda$'s will stop in a finite time in this case.  Thus one can run
the two searches in parallel, and stop whenever one of them stops.
One of the two searches will always stop in a finite amount of time,
and thereby give an answer of ``yes'' or ``no'' to the original
question (about whether $z$ lies in the subgroup of ${\bf Z}^r$
generated by $d^1, d^2, \ldots, d^q$).

	Although naive, this argument fits nicely with what happens
for the problem of deciding whether the fundamental group of a finite
complex is trivial.  When the answer is ``yes'', one can find this out
in finite time, again through exhaustive searches.  In algebraic
terms, one searches for realizations of the generators of the
fundamental group as trivial words, using the relations for the
fundamental group which can be read off from the given complex.  In
general there is no finite test for the nontriviality of words,
however, and indeed the original question is not algorithmically
decidable.

	In some cases, one might have extra information which does
allow for effective tests for answers of ``no'', and abelian groups
are a very special instance of this.

	Let us come back now to the assertion above, that $z \in {\bf
Z}^r$ does not lie in the subgroup generated by $d^1, d^2, \ldots, d^q
\in {\bf Z}^r$ if and only if there is a homomorphism $\sigma : {\bf
Z}^r \to {\bf Q}/{\bf Z}$ such that $\sigma(d^j) = 0$ for $j = 1, 2,
\ldots, q$ and $\sigma(z) \ne 0$.  Of course this is closely analogous
to familiar statements about vectors in a vector space and linear
mappings into the ground field.

	The ``if'' part of the statement above is immediate, and so
it suffices to consider the ``only if'' part.  Thus we assume that
$z$ does not lie in the subgroup generated by the $d^j$'s, and we
want to find a homomorphism $\sigma : {\bf Z}^r \to {\bf Q}/{\bf Z}$
with the required properties.

	We begin by setting $\sigma$ to be $0$ on the subgroup
generated by $d^1, d^2, \ldots, d^q$.  For $\sigma(z)$ we have to be
slightly careful.  If there is a positive integer $n$ such that $n z$
lies in the subgroup generated by $d^1, d^2, \ldots, d^q$, then we
need to choose $\sigma(z)$ so that $n \sigma(z) = 0$.  If no such $n$
exists, take $\sigma(z)$ to be the element of ${\bf Q}/{\bf Z}$
corresponding to $1/2 \in {\bf Q}$.  If such an $n$ does exist, let
$n_0$ be the smallest positive integer with that property.  Thus $n_0
> 1$, since $z$ itself does not lie in the subgroup generated by the
$d^j$'s.  In this case we take $\sigma(z)$ to be the element of ${\bf
Q}/{\bf Z}$ which corresponds to $1/n_0$.  This element is not $0$,
since $n_0 > 1$, but $\sigma(n_0 z)$ is then $0$ in ${\bf Q}/{\bf Z}$.

	We now extend $\sigma$ to the subgroup generated by $z$ and
the $d^j$'s, in the obvious way (so that $\sigma$ is a homomorphism).
One should be a bit careful here too, i.e., that this can be done in a
consistent manner, so that $\sigma$ really is well-defined on the
subgroup generated by $z$ and the $d^j$'s.  This comes down to the
fact that if $n$ is an integer such that $n z$ lies in the subgroup
generated by $d^1, d^2, \ldots, d^q$, then $n$ should be divisible by
$n_0$, which ensures that $\sigma(n z)$ is equal to $0$ in ${\bf
Q}/{\bf Z}$.  These things are not hard to check.

	Now we simply want to extend $\sigma$ to all of ${\bf Z}^r$,
in such a way that it is still a homomorphism into ${\bf Q}/{\bf Z}$.
This is not difficult to do; for a general assertion along these
lines, see Theorem 4.2 on p312 of \cite{Massey}.  The main point is
that ${\bf Q}/{\bf Z}$ is \emph{divisible}, which means that for each
element $x$ of ${\bf Q}/{\bf Z}$ and each nonzero integer $m$, there
is a $y \in {\bf Q}/{\bf Z}$ such that $my = x$.  In order to extend
$\sigma$ to all of ${\bf Z}^r$, one can extend it to new elements one
at a time, and to the subgroups that they generate together with the
subgroup of ${\bf Z}^r$ on which $\sigma$ is already defined.  The
divisibility property of ${\bf Q}/{\bf Z}$ guarantees that there are
always values available in ${\bf Q}/{\bf Z}$ by which to make
well-defined extensions.  That is, if $\sigma$ is already defined on
some subgroup $H$ of ${\bf Z}^r$, and $w$ is an element in ${\bf Z}^r$
not in $H$, then there is always a point in ${\bf Q}/{\bf Z}$ to use
as the value of $\sigma$ at $w$.  This is not a problem if $n w$ does
not lie in $H$ for any nonzero integer $n$ --- in which case one could
just as well take $\sigma(w) = 0$ --- but if $n w \in H$ for some
nonzero $n$, then one has to choose $\sigma(w)$ so that $n \sigma(w)$
is equal to $\sigma(n w)$, where the latter is already been determined
by the definition of $\sigma$ on $H$.

	By repeating this process, one can eventually extend $\sigma$
so that it becomes a homomorphism from all of ${\bf Z}^r$ into ${\bf
Q}/{\bf Z}$.  For instance, one can apply this process to the standard
basis vectors in ${\bf Z}^r$, at least when they are not already
included in the subgroup of ${\bf Z}^r$ on which $\sigma$ has already
been defined (at the given stage of the construction).  In the end one
obtains a homomorphism $\sigma : {\bf Z}^r \to {\bf Q}/{\bf Z}$ such
that $\sigma(d^j) = 0$ for $1 \le j \le q$ and $\sigma(z) \ne 0$, as
desired.

\end{document}